\documentclass[a4paper,12pt]{amsart}

\pdfoutput=1 

\usepackage[top=3cm,bottom=3cm,outer=3cm,inner=2cm,marginpar=2.45cm]{geometry}
\usepackage[abbrev]{amsrefs}

\usepackage[french,ngerman,english]{babel}
\usepackage[utf8]{inputenc}
\usepackage[T1]{fontenc}

\usepackage[all,pdf]{xy}

\usepackage{enumitem}
\usepackage{mathtools}  
\usepackage[usenames,dvipsnames]{xcolor}
\usepackage{calc} 

\babeltags{de = ngerman}
\babeltags{fr = french}

\usepackage[no-mario]{marionotations}

\usepackage[Smaller]{cancel} 

\usepackage[final,colorlinks=true]{hyperref}
\usepackage{breakurl}  

\usepackage{subfiles}

\newcommand{\termRcolor}{NavyBlue}
\newcommand{\termQcolor}{ForestGreen}
\newcommand{\termNcolor}{VioletRed}
\newcommand{\termNpcolor}{violet}

\newenvironment{termR}{\color{\termRcolor}}{}

\newenvironment{termN}{\color{\termNcolor}}{}

\newcommand{\alert}[2][RoyalBlue]{{\color{#1}#2}}

\newcommand{\defaultDimension}{n}
\newcommand{\setDefaultDimension}[1]{\renewcommand{\defaultDimension}{#1}}

\newcommand{\defaultMetric}{\clomega}
\newcommand{\setDefaultMetric}[1]{\renewcommand{\defaultMetric}{#1}}

\newcommand{\defaultlclocus}{D}

\newcommand{\defaultpsi}{\psi_D}

\NewDocumentCommand{\logKX}{
  t{M} 
}{K_X \otimes D \otimes F \IfBooleanT{#1}{\otimes M}}

\NewDocumentCommand{\vphilist}{
  D||{\vphi}           
  t{F}                 
  t{M}                 
  d()                  
  D<>{\defaultMetric}  
}{\IfBooleanTF{#2}{\vphi_F}{#1} \IfBooleanT{#3}{+\vphi_M}, \IfNoValueF{#4}{(#4),} #5}

\NewDocumentCommand{\Ltwo}{ 
  D//{\bullet,\bullet}   
  O{X}                   
  s                      
  m                      
}{L^{#1}_{(2)}\paren{\IfBooleanF{#3}{#2;} #4}}

\NewDocumentCommand{\Ltwosp}{
  D//{q}               
  t{M}                 
  D||{\vphi}           
  o 
  D<>{\defaultMetric}  
}{\Ltwo/n,#1/*{D \otimes F \IfBooleanT{#2}{\otimes M}}_{#3, \IfNoValueF{#4}{#4,} #5}}

\NewDocumentCommand{\Harm}{ 
  D//{q}               
  D||{\vphi}           
  o 
  D<>{\defaultMetric}  
}{\mathcal{H}^{n,#1}_{#2, \IfNoValueF{#3}{#3,} #4}}

\newcommand{\spHbase}{\mathcal{H}}
\NewDocumentCommand{\spH}{ 
  D//{q}         
  t{M}           
  s              
  O{\sigma}      
  t{|}           
  D<>{\sigma-1}  
  d()            
}{\spHbase^{\defaultDimension,#1}\IfNoValueTF{#7}{
    \begingroup%
    \newcommand{\upidl}{\IfBooleanTF{#3}{
        \mtidlof{\vphi_{F \IfBooleanT{#2}{\otimes M}}}
      }{\aidlof|#4|{\vphi_{F \IfBooleanT{#2}{\otimes M}}}}
    }%
    \paren{\IfBooleanT{#2}{M\otimes}
      \IfBooleanTF{#5}{
        \frac{\upidl}{\aidlof|#6|{\vphi_{F \IfBooleanT{#2}{\otimes M}}}}
      }{\upidl}}
    \endgroup%
  }{\paren{\IfBooleanT{#2}{M\otimes}#7}}}

\DeclareMathOperator{\lc}{lc} 
\NewDocumentCommand{\lcc}{ 
  O{\sigma}                
  D<>{X}                   
  D(){\defaultlclocus}
}{\lc_{#2}^{#1}\paren{#3}}

\NewDocumentCommand{\lcS}{  
  s            
  D//{\defaultlclocus}       
  D||{\sigma}  
  O{p}         
}{\mathtt{\IfBooleanT{#1}{\rs} #2}^{#3}_{#4}}

\NewDocumentCommand{\PRes}{ 
  O{}      
  d()      
}{\mathcal R_{#1}\IfNoValueF{#2}{\paren{#2}}}

\newcommand{\defidlof}[1]{\mathcal{I}_{#1}}  
\NewDocumentCommand{\mtidlof}{   
  O{}      
  D<>{#1}  
  m        
}{\multidl_{#2}\paren{#3}} 

\NewDocumentCommand{\residlof}{  
  D||{\sigma}   
  d<>           
  m             
  s             
  O{RoyalBlue}  
}{\sheaf R_{\IfNoValueTF{#2}{}{#2,} #1}\IfBooleanT{#4}{^{{\color{#5}\infty}}}\paren{#3}}

\NewDocumentCommand{\Adjidlof}{
  D||{\sigma}  
  O{X}         
  D<>{\defaultlclocus}       
  m            
}{\operatorname{\mathit{Adj}}^{#1}_{\paren{#2,#3}}\paren{#4}}

\NewDocumentCommand{\aidlof}{
  D||{\sigma}     
  d<>             
  m               
  D//{}           
  O{\defaultpsi}  
  s               
}{\sheaf{J}_{\!\IfNoValueTF{#2}{}{#2,} #1}\IfBooleanF{#6}{\paren{#3; #4#5}}}

\NewDocumentCommand{\lcV}{ 
  D||{\sigma}      
  D//{\vphi_L}     
  d()              
  O{\defaultpsi}   
}{\:d\operatorname{lcv}^{#1}_{\IfNoValueF{#3}{#3,}#2}\left[#4\right]}

\newcommand{\dvol}{\:d\vol}

\newcommand{\RTFsym}{\mathfrak{F}} 
\NewDocumentCommand{\RTF}{ 
  s    
  G{\RTFsym} 
  d//  
  o    
  >{\SplitArgument{1}{,}} d<> 
  d||  
  d()  
  o    
}{%
  \begingroup%
    \newif\ifboolup%
    \booluptrue%
    \IfNoValueT{#4}{\IfNoValueT{#5}{\IfNoValueT{#6}{\boolupfalse}}}%
    \IfNoValueT{#7}{\boolupfalse}%
    \newcommand{\srptstr}{\cramped{{}^{\IfNoValueF{#4}{#4}\IfNoValueF{#5}{\inner#5}\IfNoValueF{#6}{\abs{#6}^2}}%
      \ifboolup _
      \fi{\ifboolup\displaystyle\fi\IfNoValueF{#7}{\paren{#7}}\IfNoValueF{#8}{%
          \ifboolup {\scriptstyle #8} \else _{#8} \fi%
        }}}}%
    \ifboolup%
      \IfBooleanTF{#1}{
        \smash[t]{
          \IfNoValueF{#3}{{}^{#3}}#2\raisebox{\depthof{$\srptstr$} * \real{0.3}}{$\srptstr$}%
        }%
      }{\IfNoValueF{#3}{{}^{#3}}#2\raisebox{\depthof{$\srptstr$} * \real{0.3}}{$\srptstr$}}%
    \else%
      \IfNoValueF{#3}{{}^{#3}}#2\srptstr%
    \fi%
  \endgroup%
} 

\def\RTI{\RTF{\mathfrak{I}}}

\NewDocumentCommand{\mtlog}{O{e} d() D||{\defaultpsi}}{\log\!#1^{\paren{#2}}\abs{#3}}
\NewDocumentCommand{\slog}{O{e} D||{\defaultpsi}}{\log\abs{#1 #2}}
\NewDocumentCommand{\dlog}{O{e} D||{\defaultpsi}}{\mtlog[#1](2)|#2|}

\NewDocumentCommand{\logpole}{ 
  D||{\defaultpsi}     
  D//{\sigma}   
  O{e}          
  D<>{1+\eps}   
  s             
}{\abs{#1}^{#2} \IfBooleanTF{#5}{\slog[#3]|#1|}{\paren{\slog[#3]|#1|}^{#4}}}

\DeclareFontFamily{OMX}{MnSymbolE}{}
\DeclareSymbolFont{MnLargeSymbols}{OMX}{MnSymbolE}{m}{n}
\SetSymbolFont{MnLargeSymbols}{bold}{OMX}{MnSymbolE}{b}{n}
\DeclareFontShape{OMX}{MnSymbolE}{m}{n}{
    <-6>  MnSymbolE5
   <6-7>  MnSymbolE6
   <7-8>  MnSymbolE7
   <8-9>  MnSymbolE8
   <9-10> MnSymbolE9
  <10-12> MnSymbolE10
  <12->   MnSymbolE12
}{}
\DeclareFontShape{OMX}{MnSymbolE}{b}{n}{
    <-6>  MnSymbolE-Bold5
   <6-7>  MnSymbolE-Bold6
   <7-8>  MnSymbolE-Bold7
   <8-9>  MnSymbolE-Bold8
   <9-10> MnSymbolE-Bold9
  <10-12> MnSymbolE-Bold10
  <12->   MnSymbolE-Bold12
}{}
\DeclareMathDelimiter{\llangle}{\mathopen}%
{MnLargeSymbols}{'164}{MnLargeSymbols}{'164}
\DeclareMathDelimiter{\rrangle}{\mathclose}%
{MnLargeSymbols}{'171}{MnLargeSymbols}{'171}

\newcommand{\iinner}[2]{\left\llangle#1,#2\right\rrangle}

\NewDocumentCommand{\idxup}{ 
  m                  
  O{\defaultMetric}  
  s                  
}{\paren{#1}^{\mathrlap{\!\IfBooleanTF{#3}{\smash[t]{#2}}{#2}}\makebox[\maxof{\widthof{$#2$}-\widthof{$\!\omega$}}{0pt}]{}}}

\newcommand{\dbadj}{\dbar^{\smash{\mathrlap{*}\;\:}}}

\NewDocumentCommand{\dep}{t{;} d<> O{\nu} m}{#4\IfBooleanTF{#1}{_}{^}{\IfNoValueF{#2}{#2\:}(#3)}}

\NewDocumentCommand{\sm}{s m}{{#2}\IfBooleanTF{#1}{_}{^}\text{sm}}

\newcommand{\cvr}[1]{\mathfrak{#1}} 
\NewDocumentCommand{\rs}{ 
  s  
  m  
}{\IfBooleanTF{#1}{\smash[t]{\widetilde{#2}}}{\widetilde{#2}}}

\NewDocumentCommand{\clt}{m}{\widetilde{#1}} 
\NewDocumentCommand{\clomega}{O{\omega}}{{\clt{#1}}} 

\newcommand{\BK}{\text{(BK)}}
\newcommand{\tBK}{\text{(tBK)}}
\DeclareMathOperator{\mlc}{mlc} 

\newcommand{\sect}[1][s]{\mathtt{#1}} 

\NewDocumentCommand{\cbn}{  
  D//{\sigma_V}
  O{\sigma}
}{\mathfrak{C}^{#1}_{#2}} 
\NewDocumentCommand{\Iset}{  
  D//{V}
  O{\sigma}
}{I^{#1}_{#2}} 

\newtheorem{prop}{Proposition}[subsection]
\newtheorem{lemma}[prop]{Lemma}
\newtheorem{thm}[prop]{Theorem}
\newtheorem{cor}[prop]{Corollary}

\newtheorem{conjecture}[prop]{Conjecture}

\theoremstyle{remark}
\newtheorem{remark}[prop]{Remark}

\theoremstyle{definition}
\newtheorem{definition}[prop]{Definition}

\newtheorem{notation}[prop]{Notation}

\numberwithin{equation}{subsection}

\allowdisplaybreaks  

\begin{document}

\newcommand{\titlestr}{%
  On an injectivity theorem for log-canonical pairs \\ with analytic
  adjoint ideal sheaves%
}

\newcommand{\shorttitlestr}{%
  On an injectivity theorem for lc pairs with adjoint ideal sheaves%
}

\newcommand{\MCname}{Tsz On Mario Chan}
\newcommand{\MCnameshort}{Mario Chan}
\newcommand{\MCemail}{mariochan@pusan.ac.kr}

\newcommand{\YJname}{Young-Jun Choi}
\newcommand{\YJnameshort}{Young-Jun Choi}
\newcommand{\YJemail}{youngjun.choi@pusan.ac.kr}

\newcommand{\addressstr}{%
  Dept.~of Mathematics, Pusan National
  University, Busan 46241, South Korea%
}

\newcommand{\subjclassstr}[1][,]{%
  32J25 (primary)#1  
  32Q15#1   
  14B05 (secondary)
}

\newcommand{\keywordstr}[1][,]{%
  $L^2$ injectivity#1
  adjoint ideal sheaf#1
  multiplier ideal sheaf#1
  log-canonical centres%
}

\newcommand{\dedicatorystr}{%
  In the memory of Prof.~Jean-Pierre Demailly%
}

\newcommand{\thankstr}{%
  The first author would like to thank Shin-ichi Matsumura for drawing
  his attention to Fujino's conjecture.
  The authors started their collaboration while they were respectively
  being a visitor and a postdoc of Jean-Pierre Demailly in
  \textfr{Institut Fourier}, whose immense influence on and
  encouragement to both authors can never be overstated.
  This work was supported by the National Research Foundation (NRF) of
  Korea grant funded by the Korea government (No.~2018R1C1B3005963 and
  No.~2023R1A2C1007227).
}


\title[\shorttitlestr]{\titlestr}
 
\author[\MCnameshort]{\MCname}
\email{\MCemail}

\author{\YJname}
\email{\YJemail}
\address{\addressstr}

\thanks{\thankstr}
 
\subjclass[2020]{\subjclassstr}


\dedicatory{\dedicatorystr}

\begin{abstract}

As an application of the residue functions corresponding to the
lc-measures developed by the authors, the proof of the
injectivity theorem on compact K\"ahler manifolds for plt pairs by
Matsumura is improved in this article to allow multiplier ideal
sheaves of plurisubharmonic functions with neat analytic singularities
in the coefficients of the relevant cohomology groups.
With the use of a refined version of analytic adjoint ideal sheaves, a
plan towards a solution to the generalised version of Fujino's conjecture
(i.e.~an injectivity theorem on compact K\"ahler manifolds for lc
pairs with multiplier ideal sheaves) is laid down and, in addition to
the result for plt pairs, a proof for lc pairs in dimension $2$, which
is also an improvement of Matsumura's result, is given.


\end{abstract} 

\date{\today} 

\maketitle



\section{Introduction}



The injectivity theorem was first formulated by \textfr{Kollár} in
\cite{Kollar_injectivity}*{Thm.~2.2} (originally as a means for
proving the torsion-freeness of the higher direct images of the
canonical sheaf under a proper morphism) in the algebraic setting
which can be viewed as a generalisation of the celebrated Kodaira
vanishing theorem (see, for example,
\cite{Esnault&Viehweg_book}*{Cor.~5.2} or
\cite{Lazarsfeld_book-I}*{Remark 4.3.8}).
It was generalised to the setting on compact \textde{Kähler} manifolds
by Enoki in \cite{Enoki} using harmonic theory.

The theorem is further generalised to log-canonical (lc) pairs in the
algebraic setting via the theory of (mixed) Hodge structures (see, for
example, \cite{Esnault&Viehweg_book}*{\S 5}, \cite{Fujino_log-MMP}*{\S
  5} and \cite{Ambro_injectivity}).
On the transcendental side, the latest results that the authors are
aware of are those of Fujino (\cite{Fujino_injectivity}), Matsumura
(\cite{Matsumura_injectivity}) and Gongyo--Matsumura
(\cite{Gongyo&Matsumura}), who obtain the injectivity theorem
in the setting on compact \textde{Kähler} manifolds for Kawamata
log-terminal (klt) pairs, and that of Matsumura
(\cite{Matsumura_injectivity-lc}) for purely log-terminal (plt) pairs,
which are sub-cases of lc pairs.
All these results make use of $L^2$ theory.
This is the starting point of the current research.
Readers are referred to \cite{Fujino_injectivity},
\cite{Matsumura_injectivity} and \cite{Matsumura_injectivity-lc} for
more references on the development of the injectivity theorem.
See also \cite{Matsumura_injectivity-Kaehler} for the development of
the injectivity theorem in the relative setting.
Moreover, readers are also referred to \cite{Kollar_Sing-of-MMP}*{Def.~2.8}
for the precise definitions of the various notions of singularities,
including klt, plt, dlt and lc, in algebraic geometry.

In addition to the fact that this topic is an interesting research
problem, another initiative of the study is to illustrate the use of 
the computation of residue functions corresponding to lc-measures
studied in \cite{Chan&Choi_ext-with-lcv-codim-1} and
\cite{Chan_on-L2-ext-with-lc-measures}, and also the results on the
analytic adjoint ideal sheaves studied in
\cite{Chan_adjoint-ideal-nas}.

Readers who would like to skip the background and motivation of the
statements in this work may go directly to
Sections \ref{sec:main-result-and-strategy} and
\ref{sec:towards-gen-Fujino-conj} for the main results of this
article.

\subsection{Fujino's conjecture and Matsumura's results}
\label{sec:fujino-conj-and-Matsumura-results}

Let $X$ be a compact \textde{Kähler} manifold of dimension $n$, $D$ a
\emph{reduced divisor} on $X$ with simple normal crossings (snc) (see,
for example, \cite{Kollar_Sing-of-MMP}*{Def.~1.7}), and $F$ a
holomorphic line bundle on $X$.
Under the snc assumption on $D$, the lc centres of $(X,D)$ are simply
the irreducible components (with reduced structure) of any
intersections of irreducible components of $D$ (see, for example,
\cite{Kollar_Sing-of-MMP}*{Def.~4.15} for the precise definition of lc
centres; see also \cite{Chan&Choi_ext-with-lcv-codim-1}*{Def.~1.4.1}
for the authors' attempt to generalise to the case when $D$ may not be
a divisor but the zero locus of certain multiplier ideal sheaf).
Fujino's conjecture on the injectivity theorem can be stated as
follows.
\begin{conjecture}[\cite{Fujino_survey}*{Conj.~2.21}] \label{conj:Fujino-conj}
  Suppose that $F$ is semi-positive and there exists a holomorphic
  section $s$ of $F^{\otimes m}$ on $X$ for some positive integer $m$
  such that $s$ does not vanish identically on any lc centres of
  $(X,D)$.
  Then, the multiplication map induced by $\otimes s$,
  \begin{equation*}
    \cohgp q[X]{K_X \otimes D\otimes F} \xrightarrow{\;\otimes s\;}
    \cohgp q[X]{K_X \otimes D\otimes F^{\otimes (m+1)}} \; ,
  \end{equation*}
  is injective for every $q \geq 0$.
\end{conjecture}
In the algebraic setting, i.e.~$X$ being a smooth projective manifold,
the assumption on $F$ is replaced by that $F$ being \emph{semi-ample}.
Note that semi-ampleness implies semi-positivity (see, for example,
\cite{Fujino_injectivity}*{Lemma 1.6}).
The conjecture is then known to hold true in the algebraic setting via
the theory of mixed Hodge structures (see \cite{Fujino_log-MMP}*{\S 5
  and \S 6}).

Let $M$ be another holomorphic line bundle on the compact
\textde{Kähler} manifold $X$ and let $e^{-\vphi_F}$ and $e^{-\vphi_M}$
be hermitian metrics on $F$ and $M$ respectively. 
Matsumura proves the conjecture for the case where $(X,D)$ being plt,
i.e.~$D$ is a smooth subvariety consisting of disjoint irreducible
components.
\begin{thm}[\cite{Matsumura_injectivity-lc}*{Thm.~1.3 and
    Cor.~1.4}] \label{thm:Matsumura-plt}
  Suppose that $(X,D)$ is plt, and suppose that $\vphi_F$ and
  $\vphi_M$ are smooth such that their curvature forms satisfy
  \begin{equation*}
    0 \leq \ibddbar \vphi_M \leq C \ibddbar \vphi_F \quad\text{ on } X
  \end{equation*}
  for some constant $C > 0$ (so $F$ is semi-positive in particular).
  Let $s$ be a holomorphic section of $M$ on $X$ such that $s$ does
  not vanish identically on any lc centres of $(X,D)$.
  Then, the multiplication map induced by $\otimes s$,
  \begin{equation*}
    \cohgp q[X]{K_X \otimes D\otimes F} \xrightarrow{\;\otimes s\;}
    \cohgp q[X]{K_X \otimes D\otimes F\otimes M} \; ,
  \end{equation*}
  is injective for every $q \geq 0$.

  In particular, the conclusion holds true when $M =F^{\otimes m}$ for
  some positive integer $m$ (with $\vphi_M := m\vphi_F$),
  i.e.~Conjecture \ref{conj:Fujino-conj} holds true for any plt pairs
  $(X,D)$.
\end{thm}
Let $\phi_D$ be a potential (of the curvature of the metric
$e^{-\phi_D}$) on $D$ (see Notation
\ref{notation:potential-definition}) induced from canonical sections
of irreducible components of $D$ (see Notation
\ref{notation:potentials}).
The proof of the theorem in \cite{Matsumura_injectivity-lc} is
proceeded by reducing the original question to the questions on the
injectivity of the composition of maps
\begin{equation*}
  \xymatrix{
    {\cohgp q[X]{K_X\otimes D\otimes F\otimes \mtidlof{\phi_D}}}
    \ar[r]^-{\iota_0} \ar@/^1.7pc/[rr]|*+{\mu_0}
    & {\cohgp q[X]{K_X \otimes D\otimes F}} \ar[r]^-{\otimes s}
    & {\cohgp q[X]{K_X \otimes D\otimes F\otimes M}}
  }
\end{equation*}
(more precisely, it is to check whether $\ker\mu_0 =\ker\iota_0$)
and the injectivity of the map
\begin{equation*}
  \nu_1 \colon \cohgp q[D]{K_X\otimes D\otimes F\otimes
    \frac{\holo_X}{\mtidlof{\phi_D}}}
  \xrightarrow{\otimes \res s_{D}}
  \cohgp q[D]{K_X\otimes D\otimes F\otimes M\otimes
    \frac{\holo_X}{\mtidlof{\phi_D}}}
\end{equation*}
(see \S \ref{sec:towards-gen-Fujino-conj}), where $\mtidlof{\phi_D} =
\mtidlof[X]{\phi_D}$ is the multiplier ideal sheaf of $\phi_D$ on $X$
and the map $\iota_0$ that $\mu_0$ factors through is induced by the
inclusion $\mtidlof{\phi_D} \subset \holo_X$.
Since $(X,D)$ is plt, the map $\nu_1$ can be decomposed into a direct
sum of homomorphisms between cohomology groups on irreducible
components of $D$.
The injectivity of $\nu_1$ is thus a consequence of the injectivity
theorem of Enoki (\cite{Enoki}).
The main focus of \cite{Matsumura_injectivity-lc} is to show that
$\ker\mu_0 =\ker\iota_0$ (see \cite{Matsumura_injectivity-lc}*{Thm.~1.6}).

When the potentials $\vphi_F$ and $\vphi_M$ are allowed to be singular, there
is the following result for pseudo-effective line bundles.

\begin{thm}[\cite{Gongyo&Matsumura}*{Thm.~1.3}; see also
  \cite{Matsumura_injectivity}*{Thm.~1.3}] \label{thm:inj-thm-klt-Matsumura}
  Suppose that $D$ is any effective $\fieldR$-divisor.
  Let $\phi_D$ be a potential on $D$ induced from canonical
  sections of irreducible components of $D$ (see Notation
  \ref{notation:potentials}).
  Suppose also that
  \begin{equation*}
    \vphi_F :=a \vphi_M +\phi_D \quad\text{ for some number } a > 0
  \end{equation*} 
  and that $\vphi_M$ (and thus $\vphi_F$) is plurisubharmonic (psh)
  locally everywhere in $X$, i.e.
  \begin{equation*}
    \ibddbar\vphi_M \geq 0 \quad(\text{and thus } \ibddbar\vphi_F \geq
    0 ) \quad\text{ on } X \; .
  \end{equation*}
  Let $s$ be a non-zero holomorphic section of $M$ on $X$ such that
  $\sup_X \abs s_{\vphi_M}^2 < \infty$.
  Then, the multiplication map induced by $\otimes s$,
  \begin{equation*}
    \cohgp q[X]{K_X \otimes F \otimes \mtidlof{\vphi_F}} \xrightarrow{\;\otimes s\;}
    \cohgp q[X]{K_X \otimes F\otimes M \otimes
      \mtidlof{\vphi_F+\vphi_M}} \; ,
  \end{equation*}
  is injective for every $q \geq 0$, where $\mtidlof{\vphi}$ is the
  multiplier ideal sheaf of the potential $\vphi$.
\end{thm}
Recall that $(X,D)$ is a klt pair if, under the assumption that $D$ is
an snc $\fieldR$-divisor, the coefficient of every irreducible
component of $D$ is $< 1$. 
In this case, $\phi_D$ has only klt singularities,
i.e.~$\mtidlof{\phi_D} =\holo_X$.
When $\vphi_M$ is smooth, the above theorem can be viewed as a version
of Conjecture \ref{conj:Fujino-conj} (with $F^{\otimes m}$ replaced by
$M$) for klt pairs $(X,D)$. 

\subsection{The first result and strategy of proof}
\label{sec:main-result-and-strategy}

The goal of this research is to prove Conjecture
\ref{conj:Fujino-conj} while allowing certain multiplier ideal sheaves
in the coefficients of the cohomology groups.
This article is the first step in this direction.
By revising the proofs in \cite{Matsumura_injectivity} and
\cite{Matsumura_injectivity-lc} into the one which is, in the authors'
point of view, more favourable to the study of lc pairs, the following
generalisation of Theorem \ref{thm:Matsumura-plt}, or more precisely,
the generalisation of \cite{Matsumura_injectivity-lc}*{Thm.~1.6} (a
statement on the map $\mu_0$ in Section
\ref{sec:fujino-conj-and-Matsumura-results}), is obtained.

\begin{thm} 
  \label{thm:main-result}
  Suppose that $(X,\omega)$ is a compact \textde{Kähler} manifold and
  $D$ a reduced divisor with snc such that $(X,D)$ is lc.
  Let $\vphi_F$ and $\vphi_M$ be potentials on $F$ and $M$
  respectively such that
  \begin{itemize}
  \item $\ibddbar\vphi_F \geq 0$ and $-C\omega \leq \ibddbar\vphi_M
    \leq C\ibddbar\vphi_F$ on $X$ in the sense of currents for some
    constant $C > 0$,
  \item $\vphi_F$ and $\vphi_M$ have only neat analytic singularities,
  \item the polar sets $P_F :=\vphi_F^{-1}(-\infty)$ and $P_M
    :=\vphi_M^{-1}(-\infty)$ both contain no lc centres of $(X,D)$, and
  \item both $P_F$ and $P_M$ are divisors and $P_F\cup P_M \cup D$ has
    only snc.
  \end{itemize}
  Also let $\phi_D$ be a potential defined by a canonical section of
  $D$.
  Suppose that there exists a non-trivial holomorphic section $s \in
  \cohgp 0[X]{M}$ such that $\sup_X \abs{s}_{\vphi_M}^2 < \infty$.
  Then, given the commutative diagram
  \begin{equation*}
    \xymatrix@R=0.6cm{
      {\cohgp q[X]{K_X\otimes D\otimes F\otimes \mtidlof{\vphi_F+\phi_D}}}
      \ar[d]_-{\iota_0} \ar@{}@<8ex>[d]|*+{\circlearrowleft} \ar[dr]^-{\mu_0}
      \\
      {\cohgp q[X]{K_X\otimes D\otimes F\otimes \mtidlof{\vphi_F}}}
      \ar[r]^-{\otimes s}
      &
      {\cohgp q[X]{K_X\otimes D\otimes F\otimes M\otimes
        \mtidlof{\vphi_F+\vphi_M}} \; ,}
    }
  \end{equation*}
  in which $\iota_0$ is induced from the inclusion
  $\mtidlof{\vphi_F+\phi_D} \subset \mtidlof{\vphi_F}$,
  one has $\ker\mu_0 =\ker\iota_0$ for every $q \geq 0$.
\end{thm}
Together with the injectivity of the corresponding map $\nu_1$ in
Section \ref{sec:fujino-conj-and-Matsumura-results} (a
consequence of Theorem \ref{thm:inj-thm-klt-Matsumura}), a
statement slightly more general than Conjecture \ref{conj:Fujino-conj}
in the plt case can be proved.
See Corollary \ref{cor:gen-Fujino-conj-plt} for details.

Although Theorem \ref{thm:main-result} is only a slight improvement
($\vphi_F$ and $\vphi_M$ are allowed to be singular but only for neat
analytic singularities) to the corresponding statement in
\cite{Matsumura_injectivity-lc}*{Thm.~1.6 or Thm.~3.9}, a different
proof from that in \cite{Matsumura_injectivity-lc} is presented here.
While both proofs follow the same spirit of arguments of Enoki in
\cite{Enoki}*{\S 2} (\emph{in
  view of the Dolbeault isomorphism, consider a harmonic form $u$
  representing a class in the domain cohomology group which is in
  $\ker\mu_0$ and also in the orthogonal complement of $\ker\iota_0$,
  then argue via the $L^2$ theory and Bochner--Kodaira--Nakano identity
  to show that $u=0$ under the positivity assumption on $F$ and $M$}),
the two proofs differ in the handling of the non-integrable lc
singularities in the potentials (namely, $\phi_D$).

Assume that $\vphi_F$ and $\vphi_M$ are smooth for the moment.
Let $\sect_D$ be a canonical section of $D$ such that $\phi_D =
\log\abs{\sect_D}^2$ and let $\sm\vphi_D$ be a smooth potential on $D$.
In the proof of \cite{Matsumura_injectivity-lc}*{Thm.~3.9}, $\phi_D$
is smoothed to
\begin{equation*}
  \dep[\eps]\phi_D :=\log\paren{\abs{\sect_D}_{\sm\vphi_D}^2
    +\eps} +\sm\vphi_D \quad\paren{\in \smooth_X}\; .
\end{equation*}
In view of the $L^2$ Dolbeault isomorphism (which is named as de Rham--Weil
isomorphism in \cite{Matsumura_injectivity} and
\cite{Matsumura_injectivity-lc} \footnote{\label{fn:L2-Dolbeault-name}%
  The name ``de Rham--Weil isomorphism'' is used in
  \cite{Demailly}*{Ch.~IV, \S 6} to mean more generally the
  isomorphisms between the cohomology of a sheaf and the cohomology of
  an acyclic resolution of the sheaf.
  The so named isomorphism in \cite{Matsumura_injectivity} and
  \cite{Matsumura_injectivity-lc} stands for, more specifically, the
  isomorphisms between the \v Cech cohomology of a multiplier ideal
  sheaf and the $\dbar$-cohomology computed from the associated
  Dolbeault complex of locally $L^2$ forms (with respect to some $L^2$
  norm with a possibly singular weight).
  This \emph{latter} type of isomorphisms, while not named in
  \cite{Takegoshi_higher-direct-images}*{Prop.~4.6},
  \cite{Ohsawa_book}*{Thm.~4.13} and \cite{Fujino_injectivity}*{Claim
    1} when it is stated or proved, is named as ``\emph{Leray
    isomorphism}'' in \cite{Siu_non-Kaehler}*{\S 2}.
  The version of such isomorphism by Fujino in
  \cite{Fujino_injectivity}*{Claim 1} allows the involving $L^2$ norm
  to be the one induced from a singular quasi-psh potential which is
  smooth on a Zariski open set, while the one by Matsumura in
  \cite{Matsumura_injectivity}*{Prop.~5.5} allows the quasi-psh
  potential to have arbitrary singularities and also allows a more
  flexible choice of the \textde{Kähler} form
  (\cite{Matsumura_injectivity-lc}*{Prop.~2.8}).
  In this paper, the version in \cite{Matsumura_injectivity} is used.
  Although it may be more proper to attribute the isomorphism to
  Fujino and Matsumura, the authors incline to name it as
  ``\emph{$L^2$ Dolbeault isomorphism}'' ($L^2$ version of the Dolbeault
  isomorphism), which seems to be more suggestive and self-explanatory.
}),
let $[u]$ be a cohomology class in $\ker\mu_0$ and let
$u$ be the \emph{harmonic} $D\otimes F$-valued $(n,q)$-form with respect to
the (global) $L^2$ norm $\norm\cdot_{\vphi_F+\phi_D}$ induced from
$\vphi_F+\phi_D$ which represents the class $[u]$.
Then, $s u =\dbar v$ for some $D\otimes F\otimes M$-valued
$(n,q-1)$-form $v$ which is $L^2$ in $\norm\cdot_{\vphi_F+\vphi_M+\sm\vphi_D}$
(but not clear whether it is $L^2$ in
$\norm\cdot_{\vphi_F+\vphi_M\alert{+\phi_D}}$).
The proof is based on the inequality
\begin{equation*}
  \norm{su}^2 \xleftarrow{\eps \tendsto 0^+}
  \norm{su}_{(\eps)}^2 =\iinner{su}{\dbar v}_{(\eps)}
  =\iinner{\dfadj_{(\eps)}(su)}{v}_{(\eps)}
  \leq \norm{\dfadj_{(\eps)}(su)}_{(\eps)} \: \norm{v}_{(\eps)} \; ,
\end{equation*}
where $\norm\cdot =\norm\cdot_{\phi_D
+\vphi_F+\vphi_M}$ while $\norm\cdot_{(\eps)}$ and
$\iinner\cdot\cdot_{(\eps)}$ are the norm and inner product obtained
after $\phi_D$ is smoothed to $\dep[\eps]\phi_D$, and $\dfadj_{(\eps)}$ is
the corresponding formal adjoint of $\dbar$.
In order to show that the right-hand-side converges to $0$ as $\eps
\tendsto 0^+$, the rate of divergence of the integral of
$e^{-\dep[\eps]\phi_D}$ has to be controlled so that
\begin{equation*} \tag{$*$} \label{eq:desired-estimate-for-smoothened-metric}
  \int_V e^{-\dep[\eps]\phi_D} \:d\vol_V = o\paren{\frac 1\eps}
  \quad\text{(little-o notation)}
\end{equation*}
for any local open set $V$ in $X$ as $\eps \tendsto 0^+$.
When $(X,D)$ is plt, it is easy to show that 
\begin{equation*}
  \int_V e^{-\dep[\eps]\phi_D} \:d\vol_V = \BigO\paren{\abs{\log\eps}}
  \quad\text{(Big-O notation)}
\end{equation*}
as $\eps \tendsto 0^+$, which gives the required estimate
\eqref{eq:desired-estimate-for-smoothened-metric} (see 
\cite{Matsumura_injectivity-lc}*{Lemma 3.11, Prop.~3.12 and
  Prop.~3.14}).\footnote{
  \cite{Matsumura_injectivity-lc}*{Lemma 3.11} holds
  only in the plt case.
  Indeed, on a neighbourhood $V$ such that $V\cap D =\set{r_1r_2=0}$,
  where $r_1$ and $r_2$ are the radial components of the polar
  coordinates such that $(r_1,r_2) \in [0,1)^2$ on $V$, one has
  \begin{equation*}
    \int_V e^{-\dep[\eps]\phi_D} \dvol_V
    \sim
    \int_{[0,1)^2} \frac{dr_1^2 dr_2^2}{r_1^2 r_2^2 +\eps}
    \geq
    \int_{[0,1)^2} \frac{dr_1^2 dr_2^2}{\paren{r_1^2 +\sqrt\eps}
      \paren{r_2^2 +\sqrt\eps}}
    =\BigO\paren{\abs{\log\eps}^2} \; .
  \end{equation*}
  Nevertheless, one can still obtain
  \eqref{eq:desired-estimate-for-smoothened-metric} by a simple
  adjustment, namely, when $(X,D)$ is lc but not plt and when $D
  =\sum_{i \in I} D_i$ such that each $D_i$ is irreducible and has a
  canonical section $\sect_{D_i}$ and a smooth potential $\sm\vphi_{D_i}$, set
  \begin{equation*}
    \dep[\eps]\phi_D
    :=\sum_{i \in I} \log\paren{\abs{\sect_{D_i}}_{\sm\vphi_{D_i}}^2
      +\eps} +\sm\vphi_D \; .
  \end{equation*}
  The variables can then be separated when estimating the integral of
  $e^{-\dep[\eps]\phi_D}$ and thus a direct computation yields 
  \begin{equation*}
    \int_V e^{-\dep[\eps]\phi_D} \dvol_V = \BigO\paren{\abs{\log\eps}^{\sigma_V}}
  \end{equation*}
  when $V \cap D =\set{z_1 \dotsm z_{\sigma_V} = 0}$ for some $\sigma_V
  \leq n$ (where $(z_1,\dots,z_n)$ is a holomorphic coordinate system on
  $V$).
  The estimate \eqref{eq:desired-estimate-for-smoothened-metric} then follows.}
The proof of \cite{Matsumura_injectivity-lc}*{Thm.~1.6} relies on this
estimate.

In the proof presented in this article, instead of smoothing out
the lc singularities on the potential $\phi_D$, a sequence of smooth
cut-off functions $\seq{\theta_\eps}_{\eps > 0}$ vanishing identically
on some neighbourhoods of $D$ and converging to the identity map on
$X$ as $\eps \tendsto 0^+$ is considered such that
\begin{equation*}
  \norm{su}^2 \xleftarrow{\eps \tendsto 0^+}
  \iinner{su}{\theta_\eps \dbar v}
  =\iinner{\dfadj(su)}{\theta_\eps v} -\iinner{su}{\dbar\theta_\eps
    \wedge v} \; ,
\end{equation*}
where $\dfadj$ is the formal adjoint of $\dbar$ with respect to
the potential $\phi_D+\vphi_F+\vphi_M$ (which is denoted by
$\dfadj_{\vphi_M}$ in latter sections).
It can be shown that $\dfadj\paren{su} =0$ (see Corollary
\ref{cor:dfadj_M-su=0}; also compare with the result
$\norm{\dfadj_{(\eps)}\paren{su}}_{(\eps)}^2
=\BigO\paren{\eps\abs{\log\eps}}$ in the plt case in
\cite{Matsumura_injectivity-lc}*{\S 3.2}), so it suffices to estimate
the inner product on the far right-hand-side in order to show that $u
= 0$.
A fundamental trick at play is that, although $e^{-\phi_D}$ is
non-integrable around $D$, using the computation of the
residue functions associated to lc-measures studied in
\cite{Chan&Choi_ext-with-lcv-codim-1} and
\cite{Chan_on-L2-ext-with-lc-measures} (or simply via a direct
computation), one has 
\begin{equation*}
  \eps\int_V \frac{e^{-\phi_D}}{\abs{\psi_D}^{\sigma +\eps}} \:d\vol_V
  =\BigO(1) \quad \text{ as } \eps \tendsto 0^+
  \quad\paren{\psi_D := \phi_D -\sm\vphi_D \leq -1}
\end{equation*}
on any local open set $V \subset X$ when $\sigma \geq \sigma_{\mlc}$,
where $\sigma_{\mlc}$ is the codimension of the minimal lc centres
(mlc) of $(X,D)$ (see Theorem \ref{thm:residue-fcts-and-norms} or
\cite{Chan&Choi_ext-with-lcv-codim-1}*{Prop.~2.2.1}; note also that
the integral diverges when $\sigma < \sigma_{\mlc}$).
It turns out that, with a careful analysis on the properties possessed
by $u$, in order to prove Theorem \ref{thm:main-result} for any 
values of $\sigma_{\mlc} \geq 1$ (i.e.~no matter whether $(X,D)$ is
plt or not), it suffices to put $\abs{\psi_D}^{\alert{1+\eps}}$ into the
denominators (via a suitable choice of the cut-off functions
$\theta_\eps$) of the integrand of the inner product (see Steps
\ref{item:results-from-harmonic-u-and-BK},
\ref{item:Takegoshi-argument} and
\ref{item:residue-of-final-inner-prod} of the outline of the proof of
Theorem \ref{thm:main-result} in \S \ref{sec:outline-of-pf}).

In order to prove that $u=0$, it is necessary to assume
that $u$ is sitting inside the orthogonal complement
$\paren{\ker\iota_0}^\perp$ of $\ker\iota_0$ (note that $\ker\iota_0
\neq 0$ for $q=1$, for example, when $X$ is an elliptic curve,
$F=\holo_X$ and $D$ is an effective divisor of $\deg D =1$ with the
realisation that $D \otimes \mtidlof{\phi_D} \isom \holo_X$).
An argument of Takegoshi
(see \cite{Takegoshi_cohomology-nef-line-bdl}*{Prop.~3.8} or
\cite{Matsumura_injectivity-lc}*{Prop.~3.13}; see also Step
\ref{item:Takegoshi-argument} in \S\ref{sec:outline-of-pf}) is needed
to make use of such assumption, which requires $u$ to be smooth
around the lc locus $D$ (see Remark
\ref{rem:reason-for-smooth-omega-along-D} for details).
Indeed, to compute the above inner product using the computation
of residue functions in \cite{Chan&Choi_ext-with-lcv-codim-1} and
\cite{Chan_on-L2-ext-with-lc-measures}, $u$ is
also required to be smooth around the lc locus $D$.
Using the refined version of the hard Lefschetz theorem of Matsumura
(see Theorem \ref{thm:refined-hard-Lefschetz} or
\cite{Matsumura_injectivity-lc}*{Thm.~3.3}), this can be guaranteed
when the \textde{Kähler} metric $\omega$ is smooth around $D$.
As a result, even though the metric $e^{-\phi_D}$ on $D$ is singular,
one has to keep using a \textde{Kähler} metric $\omega$ which is
incomplete on $X \setminus D$ when making use of the (twisted)
Bochner--Kodaira formula, and thus extra care is needed (see
\S\ref{sec:BK-formulas}, Proposition \ref{prop:nabla-u=0_curv-u=0} and
Corollary \ref{cor:dfadj_M-su=0}).

When $\vphi_F$ and $\vphi_M$ are not smooth but have neat analytic
singularities as in the assumption of Theorem \ref{thm:main-result},
one would expect that the arguments employed in the smooth case should
still hold true since the singularities on $\vphi_F$ and $\vphi_M$
along $P_F \cup P_M$ and the lc locus $D$ are ``separated''.
In practice, a suitably chosen complete \textde{Kähler} metric
$\clomega$ on $X^\circ := X \setminus \paren{P_F \cup P_M}$ is
considered.
The corresponding harmonic forms $u$ (denoted by $\clt u$ in latter
sections) may not be smooth along $P_F \cup P_M$, but their
singularities can be determined (see Proposition
\ref{prop:regularity-of-clt-u}) and are not interfering with the
computations around the lc locus $D$, thanks to Fubini's theorem.
The argument of Takegoshi is also adjusted to adapt to such
situation (see Step \ref{item:Takegoshi-argument} in
\S\ref{sec:outline-of-pf}).

In the following sections, $\omega$ is used to mean a fixed
(smooth) \textde{Kähler} form on $X$ and $\clomega$ a chosen complete
\textde{Kähler} form on $X^\circ$.
The harmonic forms with respect to $\omega$ and $\clomega$ in the same
class $[u]$ discussed above are denoted by $u$ and $\clt u$
respectively.

\subsection{Towards Fujino's conjecture and its generalisation}
\label{sec:towards-gen-Fujino-conj}

Already in the proof of the injectivity theorem for plt pairs in
\cite{Matsumura_injectivity-lc} involves arguments of restriction
of the relevant cohomology classes to the lc centres of $(X,D)$.
It is therefore natural to incorporate the corresponding adjoint ideal
sheaves and their residue exact sequences into the potential solution
of Fujino's conjecture.
The analytic adjoint ideal sheaves studied in
\cite{Chan_adjoint-ideal-nas} is introduced below for that purpose.

For any integer $\sigma= 1, \dots, n$, let $\lcc$ be the \emph{union of
  lc centres of $(X,D)$ of codimension $\sigma$} (or \emph{union of
  $\sigma$-lc centres} for short)
and $\defidlof{\lcc}$ be its defining ideal sheaf in $\holo_X$.
If $\sigma_{\mlc}$ is the codimension of the mlc of $(X,D)$, set
$\defidlof{\lcc} := \holo_X$ for all $\sigma > \sigma_{\mlc}$.
Let $L$ denote, in this section, either the line bundle $F$ or $F
\otimes M$ and let $\vphi_L$ be either the potential $\vphi_F$ or
$\vphi_F +\vphi_M$ accordingly.
Notice that the family $\seq{\mtidlof{\vphi_L+m\psi_D}}_{m \geq 0}$ of
multiplier ideal sheaves on $X$ has a jumping number $m=1$ as seen
from the assumptions on (the singularities of) $\vphi_L$ and $\psi_D$
($:= \phi_D -\sm\vphi_D \leq -1$) in Theorem \ref{thm:main-result}.
In \cite{Chan_adjoint-ideal-nas}, the first author introduces the
following version of analytic adjoint ideal sheaves.
\begin{definition}[\cite{Chan_adjoint-ideal-nas}*{Def.~1.2.1}] \label{def:adjoint-ideal-sheaves}
  Given any integer $\sigma \geq 0$ and a family
  $\set{\mtidlof{\vphi_L+m\psi_D}}_{m\in[0,1]}$ with a jumping number
  at $m=1$, the \emph{(analytic) adjoint ideal sheaf $\aidlof{\vphi_L}
    := \aidlof<X>{\vphi_L}$ of index $\sigma$} of $(X,\vphi_L,\psi_D)$
  is the sheaf associated to the presheaf over $X$ given by 
  \begin{equation*}
    \bigcap_{\eps > 0} \mtidlof{\vphi_L+\psi_D +\log\paren{\logpole}}\paren{V}
  \end{equation*}
  for every open set $V \subset X$.
  Then, its stalk at each $x \in X$ can be described as
  \begin{equation*}
    \aidlof{\vphi_L}_x
    =\setd{f\in \holo_{X,x}}{\exists~\text{open set } V_x \ni x \: , \;
      \forall~\eps > 0 \: , \; \frac{\abs f^2
        e^{-\vphi_L-\psi_D}}{\logpole} \in L^1\paren{V_x} } \; . 
  \end{equation*}
\end{definition}
According to \cite{Chan_adjoint-ideal-nas}*{Thm.~1.2.3}, under the
assumption that $\vphi_L$ and $\vphi_L+\psi_D$ have only neat analytic
singularities with snc, one has
\begin{equation*}
  \aidlof{\vphi_L} =\mtidlof{\vphi_L} \cdot \defidlof{\lcc[\sigma+1]}
\end{equation*}
for all integers $\sigma \geq 0$, which fit into the chain of natural
inclusions 
\begin{equation*}
  \mtidlof{\vphi_L+\phi_D} =\aidlof|0|{\vphi_L}
  \subset \aidlof|1|{\vphi_L} \subset \dotsm \subset
  \aidlof|\sigma_{\mlc}|{\vphi_L} =\mtidlof{\vphi_L} \; .
\end{equation*}
Moreover, since $\vphi_L^{-1}(-\infty)$ contains no lc centres of
$(X,D)$, the analytic adjoint ideal sheaves fit into the residue short
exact sequence
\begin{equation} \label{eq:short-ext-seq-of-ideals}
  \xymatrix{
    0 \ar[r]
    & {K_X \otimes \aidlof|\sigma-1|{\vphi_L}} \ar[r]
    & {K_X \otimes \aidlof{\vphi_L}} \ar[r]^-{\Res}
    & {K_X \otimes \residlof{\vphi_L}} \ar[r]
    & 0
  } \; ,
\end{equation}
where $\residlof{\vphi_L}$ is a coherent sheaf supported on $\lcc$
such that, on an open set $V$ with $\lcc \cap V = \bigcup_{p\in \Iset}
\lcS$, where $\lcS$'s are the $\sigma$-lc centres in $V$ indexed
by $p \in \Iset$, one has
\begin{equation} \label{eq:residl-definition}
  K_X \otimes \residlof{\vphi_L}(V)
  =\prod_{p\in \Iset} K_{\lcS} \otimes \res{
    D^{\otimes (-1)}}_{\lcS} 
  \otimes \mtidlof[\lcS]{\res{\vphi_L}_{\lcS}}\paren{\lcS} 
\end{equation}
(see \cite{Chan_adjoint-ideal-nas}*{\S 4.2} for the precise
construction of $\residlof{\vphi_L}$).
For every $f \in K_X \otimes \aidlof{\vphi_L}(V)$, the component of $\Res(f)$ on
$\lcS$ is given by
\begin{equation*}
  \PRes[\lcS](\frac{f}{\sect_D}) \; ,
\end{equation*}
where $\sect_D$ is the canonical section of $D$ such that $\phi_D
=\log\abs{\sect_D}^2$, $\PRes[\lcS]$ is the \textfr{Poincaré} residue
map corresponding to the restriction from $X$ to $\lcS$ (see
\cite{Kollar_Sing-of-MMP}*{Def.~4.1 and para.~4.18}; see also Section
\ref{sec:residue-functions}).
Readers are referred to \cite{Chan_adjoint-ideal-nas} for the
comparison between the analytic adjoint ideal sheaves introduced above
and the version studied in
\cite{Guenancia}\nocite{Guenancia_AdjIdl-Erratum} and
\cite{KimDano-adjIdl}, as well as the algebraic versions studied in
\cite{Ein&Lazarsfeld_adjIdl}, \cite{Hacon&Mckernan} and
\cite{Ein-Popa}.

For the sake of convenience, for any sheaf $\sheaf F$ on $X$, set
\begin{equation*}
  \spH(\sheaf F) :=\cohgp q[X]{K_X\otimes D\otimes
                                   F\otimes \sheaf F}
\end{equation*}
for any integer $q =0, \dots, n$.
From the residue short exact sequence \eqref{eq:short-ext-seq-of-ideals}
and the multiplication map (where $\vphi_{F \otimes M} :=\vphi_F +\vphi_M$)
\begin{equation*}
  K_X\otimes D\otimes F\otimes \aidlof{\vphi_F}
  \xrightarrow{\otimes s\;}
  K_X\otimes D\otimes F\otimes M\otimes \aidlof{\vphi_{F \otimes M}} \; ,
\end{equation*}
one obtains the following commutative diagram of cohomology groups:

\begin{equation} \label{eq:commut-diagram_Fujino-conj}
  \renewcommand{\objectstyle}{\displaystyle}
  \begin{aligned}
    \xymatrix@R=0.65cm@+0.25cm{
      {\vdots} \ar[d]
      & {\vdots}
      \ar[d]
      & {\vdots} \ar[d] \\
      {\spH[\sigma-1]} \ar[d] \ar@{=}[r]
      & {\spH[\sigma-1]}
      \ar[r]
      \ar[d]_-{\iota_{\sigma-1}}
      \ar[dr]|-*+{\mu_{\sigma-1}}
      & {\spH M[\sigma-1]} \ar[d] \\
      {\spH} \ar[d] \ar[r]^-{\iota_{\sigma}}
      \ar@/_1.69pc/[rr]|(.6)*+{\mu_{\sigma}}
      & {\spH*} \ar[d]|(.38)*+<3pt>{ } \ar[r]
      & {\spH M*} \ar[d] \\
      {\spH(\residlof{\vphi_F})} \ar[d] \ar[r]^-{\tau_\sigma}
      \ar@/_1.7pc/[rr]+<-39pt,-15pt>|(.63)*+{\nu_\sigma}
      & {\spH*|}
      \ar[d]|(.52)*+<3pt>{}
      \ar[r]
      & {\spH M*|} \ar[d] \\
      {\vdots} & {\vdots} & {\vdots} }
  \end{aligned}
\end{equation} 
Note that the columns are all exact.
The middle horizontal map on the left-hand-side is induced from the
natural inclusion $\aidlof{\vphi_F} \subset \mtidlof{\vphi_F}$, while
the horizontal maps on the right-hand-side are induced from the
multiplication map $\otimes s$.
Each homomorphism of $\mu_\sigma$'s and $\nu_\sigma$'s is the
composition of the maps on the corresponding row.

Through a simple diagram-chasing, one sees that, for each $\sigma \geq
1$, if the homomorphisms $\mu_{\sigma-1}$ and $\nu_\sigma$ satisfy
$\ker\mu_{\sigma-1} =\ker\iota_{\sigma-1}$ and $\ker\nu_\sigma
=\ker\tau_\sigma$ respectively, then it follows that $\ker\mu_{\sigma}
=\ker\iota_\sigma$.
One then obtains the following theorem via induction.
\begin{thm} \label{thm:induction-on-Fujino-conj}
  If one has $\ker\mu_0 =\ker\iota_0$ and $\ker\nu_{\sigma}
  =\ker\tau_{\sigma}$ for $\sigma =1, \dots, \sigma_{\mlc}$, then
  $\mu_{\sigma_{\mlc}}$ is injective (as $\iota_{\sigma_{\mlc}}$ is the
  identity map).
  In particular, Fujino's conjecture, which concerns about the
  situation when $\vphi_F$ and $\vphi_M$ are smooth and $M=F^{\otimes
    m}$ for some integer $m \geq 1$, holds true under the given
  assumptions.
\end{thm}

Suppose $(X,D)$ is plt and suppose that $\vphi_F$ and
$\vphi_M$ both have only neat analytic singularities.
The following corollary to Theorem \ref{thm:main-result} can then be
stated and proved.
\begin{cor}[cf.~\cite{Matsumura_injectivity-lc}*{Thm.~3.16}] \label{cor:gen-Fujino-conj-plt}
  Suppose that $(X,D)$ is a plt pair.
  Suppose that $X$, $D$, $\vphi_F$, $\vphi_M$ and $s$ satisfy all the assumptions
  in Theorem \ref{thm:main-result}.
  Assume further that the section $s\in \cohgp 0[X]{M}$
  does not vanish identically on any lc centres of $(X,D)$.
  Then, the multiplication map $\mu_1$, that is,
  \begin{equation*}
    \cohgp q[X]{K_X \otimes D\otimes F\otimes
      \mtidlof{\vphi_F}} 
    \xrightarrow{\;\otimes s\;}
    \cohgp q[X]{K_X \otimes D\otimes F \otimes M \otimes
      \mtidlof[X]{\vphi_F +\vphi_M}} \; ,
  \end{equation*}
  is injective for any integer $q \geq 0$.
  (Put $M:=F^{\otimes m}$ and $\vphi_M :=m \vphi_F$ and assume
  that $\vphi_F$ is smooth when Fujino's conjecture is concerned.)
\end{cor}

\begin{proof}
  The pair $(X,D)$ being plt means that $\sigma_{\mlc} = 1$.
  By Theorem \ref{thm:induction-on-Fujino-conj}, it suffices to show
  that $\ker\mu_0=\ker\iota_0$ and $\ker\nu_1=\ker\tau_1$.
  The equality $\ker\mu_0=\ker\iota_0$ is guaranteed by Theorem
  \ref{thm:main-result}.
  Notice that, as $\sigma_{\mlc} =1$, the homomorphism $\tau_1$ is the
  identity map.
  The goal is therefore to prove that $\nu_1$ is injective.

  Write $D = \sum_{i \in I} D_i$, where $D_i$'s are the mutually
  disjoint irreducible components of $D$. 
  Then, it follows from \eqref{eq:residl-definition} (with $\sigma =1$ and
  $\setd{D_i \cap V }{ i \in I \text{ s.t. } D_i \cap V \neq \emptyset
  }  =\setd{\lcS|1|}{p \in \Iset[1]}$) that the homomorphism
  $\nu_1$ is reduced to
  \begin{multline*}
    \nu_1 \colon
    \bigoplus_{i \in I} \cohgp q[D_i]{
      K_{D_i} \otimes \res F_{D_i} \otimes
      \mtidlof[D_i]{\res{\vphi_F}_{D_i}}
    } \\
    \longrightarrow~ \bigoplus_{i \in I} \cohgp q[D_i]{
      K_{D_i} \otimes \parres{F \otimes M}_{D_i} \otimes
      \mtidlof[D_i]{\parres{\vphi_F +\vphi_M}_{D_i}} 
    } \; ,
  \end{multline*}
  which maps the $i$-th summand to the $i$-th summand via the
  multiplication $\otimes \res s_{D_i}$.
  Write the homomorphism on the $i$-th summand as $\nu_{1,i}$, and thus
  $\nu_1 = \bigoplus_{i \in I} \nu_{1,i}$.
  Note that $\res s_{D_i}$ is non-trivial and $\res{\vphi_F}_{D_i}$ is
  psh for each $i \in I$.
  When $M =F^{\otimes m}$ and $\vphi_M =m\vphi_F$, 
  each $\nu_{1,i}$ is injective by Theorem
  \ref{thm:inj-thm-klt-Matsumura} (putting $D=0$ and $a=\frac 1m$ in
  the theorem; notice that each $D_i$ is a compact \textde{Kähler}
  manifold).
  For a more general pair $(M,\vphi_M)$ which satisfies the given
  assumptions in Theorem \ref{thm:main-result}, following the proof of
  \cite{Matsumura_injectivity}*{Thm.~1.3} or the arguments given in
  Section \ref{sec:outline-of-pf} under the current setup
  (i.e.~$\vphi_F$ and $\vphi_M$ having only neat analytic
  singularities with snc), it is easy to see that the injectivity of
  $\nu_{1,i}$ for each $D_i$ (or, more precisely, for each pair $(D_i,
  0)$) still holds true.
  In any case, this implies that $\nu_1$ itself is injective.
\end{proof}

\begin{remark} \label{rem:consequence-inj-thm-klt}
  When $(X,D)$ is an lc pair (which need not be plt) and
  $\lcc =\bigcup_{p \in I} \lcS$, where $\lcS$'s are the $\sigma$-lc
  centres, the description in the proof above implies that, if $s$
  does not vanish identically on any $\sigma$-lc centres $\lcS$, the
  multiplication map
  \begin{equation*}
    \spH(\residlof{\vphi_F}) \xrightarrow{\otimes s}
    \spH M(\residlof{\vphi_F +\vphi_M}) \; ,
  \end{equation*}
  which can be rewritten as
  \begin{multline*}
    \bigoplus_{p \in I} \cohgp q[\lcS]{
      K_{\lcS} \otimes \res F_{\lcS} \otimes
      \mtidlof[\lcS]{\res{\vphi_F}_{\lcS}}
    } \\
    \xrightarrow{\;\;\otimes s\;\;}~ \bigoplus_{p \in I} \cohgp q[\lcS]{
      K_{\lcS} \otimes \parres{F \otimes M}_{\lcS} \otimes
      \mtidlof[\lcS]{\parres{\vphi_F +\vphi_M}_{\lcS}} 
    } 
  \end{multline*}
  according to \eqref{eq:residl-definition}, in which the $p$-th summand is
  mapped to the $p$-th summand via $\otimes \res s_{\lcS}$, is indeed
  \emph{injective}.
\end{remark}

Corollary \ref{cor:gen-Fujino-conj-plt} is reduced to Theorem
\ref{thm:Matsumura-plt} of Matsumura (with a slightly relaxed
assumption on $\ibddbar\vphi_M$) when $\vphi_F$ and $\vphi_M$ are
smooth.
The corresponding statement for general lc pairs $(X,D)$ is a
generalisation of Fujino's conjecture (Conjecture
\ref{conj:Fujino-conj}).
Theorem \ref{thm:induction-on-Fujino-conj}, together with Theorem
\ref{thm:main-result}, guarantees that such generalised conjecture is
solved once it is shown that $\ker\nu_\sigma =\ker\tau_{\sigma}$ for
all $\sigma = 1,\dots, \sigma_{\mlc}$.
Even without deeper analysis of the adjoint ideal sheaves, one can
already solve the generalised version of Fujino's conjecture when
$\dim_\fieldC X =2$.
The same result for $M = F^{\otimes m}$ with smooth $\vphi_F$ and
$\vphi_M = m \vphi_F$ is obtained by Matsumura in
\cite{Matsumura_rel-vanishing-w-nd}*{Thm.~1.4}.

\begin{thm}[cf.~\cite{Matsumura_rel-vanishing-w-nd}*{Thm.~1.4}] \label{thm:Fujino-conj-lc-dim-2}
  Suppose that $X$, $D$, $\vphi_F$, $\vphi_M$ and $s$ satisfy all the
  assumptions in Theorem \ref{thm:main-result} (so, in particular,
  $(X,D)$ is an lc pair which need not be plt) and suppose also that
  $\dim_\fieldC X = 2$.
  Assume further that the section $s\in \cohgp 0[X]{M}$
  does not vanish identically on any lc centres of $(X,D)$.
  Then, the homomorphism
  \begin{equation*}
    \cohgp q[X]{K_X \otimes D\otimes F\otimes
      \mtidlof{\vphi_F}} 
    \xrightarrow{\;\otimes s\;}
    \cohgp q[X]{K_X \otimes D\otimes F\otimes M \otimes
      \mtidlof[X]{\vphi_F +\vphi_M}} \; ,
  \end{equation*}
  is injective for any integer $q \geq 0$.
\end{thm}

{
  \NewDocumentCommand{\umu}{ 
    m              
    D//{\sigma-1}  
    d()            
  }{{}^{\IfNoValueF{#3}{#3}}\upsilon^{#1}_{#2}}

  \begin{proof}
    \NewDocumentCommand{\aidlquo}{m O{\sigma-1}}{\frac{\aidlof|#1|{}*}{\aidlof|#2|{}*}}

    Set $\vphi_{F\otimes M} :=\vphi_F +\vphi_M$ for convenience and let
    \begin{equation*}
      \umu{\sigma'} := \umu{\sigma'}(q)
      \:\colon \spH[\sigma']| \xrightarrow{\otimes s} \spH M[\sigma']|
    \end{equation*}
    for any integers $\sigma, \sigma'$ and $q$ such that $1
    \leq \sigma \leq \sigma' \leq \sigma_{\mlc}$ and $q \geq 0$.
    Then $\nu_{\sigma} = \umu{\sigma_{\mlc}} \circ \tau_{\sigma}$ for
    all $\sigma =1,\dots, \sigma_{\mlc}$.
    Moreover, the discussion in Remark \ref{rem:consequence-inj-thm-klt}
    (or the injectivity theorem for the case where $D = 0$)
    implies that \emph{$\umu{\sigma}(q)$ is injective for all $\sigma
      =1,\dots, \sigma_{\mlc}$ and $q \geq 0$} (note that
    $\residlof{\vphi_L} \isom
    \frac{\aidlof{\vphi_L}}{\aidlof|\sigma-1|{\vphi_L}}$ for $L= F$ or
    $F \otimes M$).
    According to Theorem \ref{thm:induction-on-Fujino-conj} and given
    Theorem \ref{thm:main-result}, the claim in this theorem is proved
    when one shows that $\ker\nu_\sigma =\ker\tau_\sigma$ for
    $\sigma=1,\dots,\sigma_{\mlc}$.
    It therefore suffices to show that $\umu{\sigma_{\mlc}}$ is
    injective for $\sigma=1,\dots,\sigma_{\mlc}$.

    {
      \setDefaultDimension{2}

      When $\dim_\fieldC X =2$, the codimension $\sigma_{\mlc}$ of the mlc
      of $(X,D)$ can take only values $1$ or $2$.
      The case where $\sigma_{\mlc} =1$ is handled in Corollary
      \ref{cor:gen-Fujino-conj-plt}.
      Assume $\sigma_{\mlc} =2$ in what follows.

      It is known that $\umu{2}$ is injective for $\sigma =2$.
      It remains to check the injectivity of $\umu{2}/0/$ in view of
      Theorem \ref{thm:induction-on-Fujino-conj}.
      Considering the short exact sequence
      \begin{equation*}
        \renewcommand{\objectstyle}{\displaystyle}
        \xymatrix@R=0.5cm{
          0 \ar[r]
          &
          {\residlof|1|{\vphi_L}} \ar[r]
          &
          {\frac{\aidlof|2|{\vphi_L}}{\aidlof|0|{\vphi_L}}} \ar[r]
          &
          {\residlof|2|{\vphi_L}} \ar[r]
          &
          0
        }
      \end{equation*}
      for $L = F$ or $F \otimes M$ (obtained from
      \eqref{eq:short-ext-seq-of-ideals}),
      one obtains a commutative diagram
      
\begin{equation} \label{eq:commut-diagram_Fujino-conj-dim-2}
  \renewcommand{\objectstyle}{\displaystyle}
  \setDefaultDimension{2}
  \begin{aligned}
    \xymatrix@R=0.65cm@C=0.65cm{
      {0} \ar[r]
      & {\spH/q-1/(\residlof|2|{\vphi_F})}
      \ar[d]^{\delta} \ar[r]^-{\alert{\umu{2}/1/}} 
      & {\alert{\spH/q-1/ M(\residlof|2|{\vphi_{F\otimes M}})}} \ar[d]^{\delta'}
      \\
      {0} \ar[r]
      & {\spH(\residlof|1|{\vphi_F})} \ar[d] \ar[r]^-{\umu{1}/0/}
      & {\spH M(\residlof|1|{\vphi_{F\otimes M}})} \ar[d]
      \\
      & {\spH[2]|<0>} \ar[d] \ar[r]^-{\umu{2}/0/}
      & {\spH M[2]|<0>} \ar[d]
      \\
      {0} \ar[r]
      & {\spH(\residlof|2|{\vphi_F})} \ar[r]^-{\umu{2}/1/}
      & {\spH M(\residlof|2|{\vphi_{F\otimes M}})} }
  \end{aligned}
\end{equation} 
      where all columns and rows are exact.

      Notice that $\residlof|2|{\vphi_{F\otimes M}}$ is supported on
      $\lcc[2]$, which has dimension $0$ and is thus a finite set of
      points in $X$.
      Therefore,
      \begin{equation*}
        \alert{\spH/q-1/M(\residlof|2|{\vphi_{F\otimes M}})} = 0
        \quad\text{ for } q \neq 1
      \end{equation*}
      and a diagram-chasing shows that
      $\umu{2}/0/(q)$ is injective for all $q \neq 1$.

      To see that $\umu{2}/0/(q)$ is injective for $q=1$, notice that
      the map
      \begin{equation*}
        \alert{\umu{2}/1/(0)} \colon \spH/0/(\residlof|2|{\vphi_F})
        \xrightarrow{\;\otimes s\;}
        \alert{\spH/0/M(\residlof|2|{\vphi_{F\otimes M}})}
      \end{equation*}
      is an \emph{isomorphism}, as $s$ is non-zero at every point of the
      finite set $\lcc[2]$ by assumption.
      The surjectivity of $\alert{\umu{2}/1/(0)}$ makes it possible to
      show that $\umu{2}/0/(1)$ is injective via again a diagram-chasing.
      This completes the proof.
    }
  \end{proof}

  \begin{remark}
    \setDefaultDimension{2}
    
    In \cite{Matsumura_rel-vanishing-w-nd}*{Thm.~1.4} (in which
    $M=F^{\otimes m}$, $\vphi_M=m \vphi_F$ and $\vphi_F$ is smooth),
    the assumption on $s \in \cohgp 0[X]{F^{\otimes m}}$ is more
    relaxed than that in Theorem \ref{thm:Fujino-conj-lc-dim-2} in the
    sense that $s$ is required not to vanish identically only on every
    component of $(X,D)$.
    In other words, $s$ may vanish on some of the $2$-lc centres of
    $(X,D)$.
    The proof in Theorem \ref{thm:Fujino-conj-lc-dim-2} can be adjusted
    to recover also this case.
    Following the proof above, one only has to verify the injectivity of
    $\umu{2}/0/(q)$ for the case $q = 1$.
    (Note that $\alert{\umu{2}/1/(0)}$ may not be injective under the
    weakened assumption.)
    Suppose $\lcc[1] =\bigcup_{i \in I_1} D_i$ and $\lcc[2]
    =\bigcup_{p \in I_2} \lcS|2|$ (where each $\lcS|2|$ is actually a
    point).
    For every $p \in I_2$, there are $i_p,j_p \in I_1$ such that $\lcS|2|
    \in D_{i_p} \cap D_{j_p}$ and the connecting morphism $\delta$ in the
    diagram \eqref{eq:commut-diagram_Fujino-conj-dim-2} maps the
    summand in $\spH/0/(\residlof|2|{\vphi_F})$ corresponding to $p
    \in I_2$ (see Remark \ref{rem:consequence-inj-thm-klt}) into the
    sum of the summands of $\spH/1/(\residlof|1|{\vphi_F})$
    corresponding to $i_p , j_p \in I_1$.
    The same is true also for the other connecting morphism $\delta'$.
    If $s$ (the global holomorphic section of $F^{\otimes m}$) vanishes
    at the point $\lcS|2|$ (but not vanishing identically on either
    $D_{i_p}$ or $D_{j_p}$ by assumption), it follows that $\deg
    \res{F^{\otimes m}}_{D_{i_p}} > 0$, thus $\res{F^{\otimes
        m}}_{D_{i_p}}$ as well as $\res F_{D_{i_p}}$ is ample on the
    curve $D_{i_p}$ (true also for $j_p$ in place of $i_p$).
    The summand $\cohgp 1[D_{i_p}]{K_{D_{i_p}} \otimes
      \res F_{D_{i_p}}}$ in $\spH/1/(\residlof|1|{\vphi_F})$ therefore
    vanishes, and the same holds true for the corresponding summand in
    $\spH/1/(F^{\otimes m} \otimes \residlof|1|{(m+1)\vphi_F})$.
    (If $\mtidlof<D_{i_p}>{\res{(m+1)\vphi_F}_{D_{i_p}}}$ is non-trivial,
    one may, for example, put an extra assumption that the numerical
    dimension of $(\res F_{D_{i_p}}, \res{\vphi_F}_{D_{i_p}})$ should
    satisfy $\operatorname{nd}(\res F_{D_{i_p}},
    \res{\vphi_F}_{D_{i_p}}) = 1$ and apply the vanishing theorem of
    Cao in \cite{Cao_vanishing-cpt-Kahler}.)
    Let $J_2$ be the subset of $I_2$ which contains all $p \in I_2$
    such that $s$ does \emph{not} vanish at $\lcS|2|$ and let
    $\res{\spH/0/(\residlof|2|{\vphi_F})}_{J_2}$
    (resp.~$\alert{\res{\spH/0/(F^{\otimes m} \otimes
      \residlof|2|{(m+1)\vphi_F})}_{J_2}}$) be the sum of summands in
    $\spH/0/(\residlof|2|{\vphi_F})$ (resp.~$\alert{\spH/0/(F^{\otimes m} \otimes
      \residlof|2|{(m+1)\vphi_F})}$) corresponding to all $p \in J_2$.
    The vanishing result above implies that, after replacing the first row
    of the diagram \eqref{eq:commut-diagram_Fujino-conj-dim-2} by
    \begin{equation*}
      \res{\spH/0/(\residlof|2|{\vphi_F})}_{J_2}
      \xrightarrow{\;\alert{\res{\umu{2}/1/}_{J_2}}\;}
      \alert{\res{\spH/0/(F^{\otimes m} \otimes
          \residlof|2|{(m+1)\vphi_F})}_{J_2}} \; ,
    \end{equation*}
    the two columns of the diagram are still exact.
    Since $\alert{\res{\umu{2}/1/}_{J_2}}$ is now an
    \emph{isomorphism}, a diagram-chasing as in the proof of Theorem
    \ref{thm:Fujino-conj-lc-dim-2} then guarantees that $\umu{2}/0/(1)$ is
    injective.
  \end{remark}

}

The general case will be discussed in subsequent papers.

\subsection{Restrictions on the singularities of $\vphi_F$ and
  $\vphi_M$}
\label{sec:restriction-on-singularities}

It is natural to ask whether the above results can be generalised to
the setting where (quasi-)psh potentials $\vphi_F$ and $\vphi_M$
with more general singularities are allowed. 
There are two apparent constraints on the singularities of the
potentials as seen from the current exposition.

The first one comes from the refined hard Lefschetz theorem of
Matsumura (see Theorem \ref{thm:refined-hard-Lefschetz} or
\cite{Matsumura_injectivity-lc}*{Thm.~3.3}), in which $\vphi_F$ (or
$\vphi_F+\phi_D$) is required to be smooth on some Zariski open set
in $X$, although there is no restriction on its singularities on the
complement.

The other one comes from the use of the adjoint ideal sheaves introduced
in \cite{Chan_adjoint-ideal-nas}.
In that paper, all involving potentials are assumed to have neat
analytic singularities.
It is expected that the regularities of the potentials can be
relaxed,\footnote{
  In the version of adjoint ideal sheaves studied by Guenancia
  \cite{Guenancia} and Dano Kim \cite{KimDano-adjIdl}, the involving
  potential $\vphi_L$ is assumed such that $e^{\vphi_L}$ is locally
  \textde{Hölder} continuous, and the corresponding adjoint ideal sheaf
  then satisfies a residue short exact sequence similar to
  \eqref{eq:short-ext-seq-of-ideals}, thus being coherent, at least in
  the case where $(X,D)$ is plt.
  There exists a psh potential such that their adjoint ideal sheaf
  does not fit in the residue short exact sequence though (see
  \cite{Guenancia}*{Remark 2.17}).
} although the singularities on the potentials may still not be
arbitrary if one insists in the current definition of the adjoint
ideal sheaves and requires them to satisfy the residue short exact
sequences \eqref{eq:short-ext-seq-of-ideals} with the quotient sheaves
having some decent description as in \eqref{eq:residl-definition}. 

In order to allow more general singularities on $\vphi_F$ and
$\vphi_M$ in the injectivity theorem, one should first relax the
requirements on their regularities from the adjoint ideal sheaves.

\mmark{
}{Discussion on (q)-psh potentials $\vphi_F$ and $\vphi_M$ with more
  general singularities? \alert{Done.}}


\subsection{Organisation of the article}

This paper is organised as follows.

Preliminaries are given in Section \ref{sec:preliminaries}.
Sections \ref{sec:notation} and \ref{sec:setup} explain some less
commonly used notations as well as the basic setup and assumptions
used in this article. 
The $L^2$ Dolbeault isomorphism is stated in Section
\ref{sec:L2-Dolbeault-isom}, also for the purpose of fixing notation.
In view of the use of \textde{Kähler} metrics which are incomplete on
$X \setminus D$ or $X^\circ \setminus D$, justification of the
well-definedness of the formal adjoint of $\dbar$, which has
singularities along $D$, is provided in Section \ref{sec:BK-formulas}.
The (twisted) Bochner--Kodaira formulas are also stated there.
In Section \ref{sec:refined-hard-Lefschetz-thm}, the refinement
of the hard Lefschetz theorem proved in
\cite{Matsumura_injectivity-lc}*{Thm.~3.3}, with a minor adjustment
for the present use, is stated and a sketch of proof is provided.
The computation on the residue functions corresponding to
$\sigma$-lc-measures, with relaxed regularity assumptions compared to
the statements in \cite{Chan&Choi_ext-with-lcv-codim-1} and
\cite{Chan_on-L2-ext-with-lc-measures}, is given in full in Section
\ref{sec:residue-functions}.

Section \ref{sec:proof-main-result} is devoted to the proof of Theorem
\ref{thm:main-result}.
An outline of the proof is given in Section \ref{sec:outline-of-pf},
which provides the essential arguments and leaves the technical
details to latter sections.
For the sake of clarity, the technical details under the assumption
that both $\vphi_F$ and $\vphi_M$ are smooth are first presented in
Section \ref{sec:pf-smooth-vphi_FM}.
The necessary adjustments for the singular case are then presented in
Section \ref{sec:pf-singular-vphi_FM}.


\section{Preliminaries}
\label{sec:preliminaries}

\subsection{Notation}
\label{sec:notation}


In this paper, the following notations are used throughout.

\begin{notation}
  Set $\ibar := \ibardefn \;$. \ibarfootnote
\end{notation}

\begin{notation} \label{notation:potential-definition}
  Each potential $\vphi$ (of the curvature of
  a metric) on a holomorphic line bundle $L$ in the following
  represents a collection of local functions
  $\set{\vphi_\gamma}_\gamma$ with respect to some fixed local
  coordinates and trivialisation of $L$ on each open set $V_\gamma$ in
  a fixed open cover $\set{V_\gamma}_\gamma$ of $X$.  The functions
  are related by the rule
  $\vphi_\gamma = \vphi_{\gamma'} + 2\Re h_{\gamma \gamma'}$ on
  $V_\gamma \cap V_{\gamma'}$ where $e^{h_{\gamma \gamma'}}$ is a
  (holomorphic) transition function of $L$ on
  $V_\gamma \cap V_{\gamma'}$ (such that
  $s_\gamma = s_{\gamma'}e^{h_{\gamma \gamma'}}$, where $s_\gamma$ and
  $s_{\gamma'}$ are the local representatives of a section $s$ of $L$
  under the trivialisations on $V_\gamma$ and $V_{\gamma'}$
  respectively).
  Inequalities between potentials is meant to be the inequalities
  under the chosen trivialisations over open sets in the fixed open
  cover $\set{V_\gamma}_\gamma$.
\end{notation}

\begin{notation} \label{notation:potentials}
  For any prime (Cartier) divisor $E$, let
  \begin{itemize}
  \item $\phi_E := \log\abs{\sect_E}^2$, representing the collection
    $\set{\log\abs{\sect_{E,\gamma}}^2}_{\gamma}$, denote a potential (of
    the curvature of the metric) on the line bundle associated to $E$
    given by the collection of local representations
    $\set{\sect_{E,\gamma}}_{\gamma}$ of some canonical section $\sect_E$
    (thus $\phi_E$ is uniquely defined up to an additive constant);
    
  \item $\sm\vphi_E$ denote a smooth potential on the line
    bundle associated to $E$;
    
    
  \item $\psi_E := \phi_E - \sm\vphi_E$, which is a global function
    on $X$, when both $\phi_E$ and $\sm\vphi_E$ are fixed.
  \end{itemize}
  All the above definitions are extended to any $\fieldR$-divisor $E$
  by linearity.
  For notational convenience, the notations for a $\fieldR$-divisor
  and its associated $\fieldR$-line bundle are used interchangeably.
  The notation of a line bundle is also abused to mean its associated
  invertible sheaf.
\end{notation}

\begin{notation} \label{notation:norm-of-nq-form}
  For any $(n,0)$-form (or $K_X$-valued section) $f$ where
  $n=\dim_\fieldC X$, define $\abs f^2 := c_n f \wedge \conj f$, where
  $c_n := (-1)^{\frac{n(n-1)}{2}}\paren{\pi\ibar}^n$.
  For any hermitian metric $\omega =\pi\ibar \sum_{1\leq j,k\leq
    n} h_{j\conj k} \:dz^j \wedge d\conj{z^k}$ on $X$, set $d\vol_{X,\omega} :=
  \frac{\omega^{\wedge n}}{n!}$.
  When $f$ is an $(n,q)$-form with $q \geq 1$ ($n=\dim_\fieldC X$), its pointwise (squared)
  norm with respect to $\omega$ for the $(0,q)$-directions is written as $\abs
  f_{\omega}^2$ (which can be viewed as a non-negative $(n,n)$-form or
  $K_X \otimes \conj{K_X}$-valued section).
  The same convention for the symbol ``$\abs f_\omega^2$'' applies also
  to the setup where $X$ is replaced by its submanifolds.
  When the Hodge $*$-operator with respect to $\omega$ is involved, 
  let $\abs{*_\omega f}_\omega^2$ denote the \emph{function-valued}
  pointwise (squared) norm with respect to $\omega$ (with differential
  forms of $*_\omega f$ in all directions being contracted). 
\end{notation}

\begin{notation} \label{notation:ineq-up-to-constants}
  For any two non-negative functions $u$ and $v$,
  write $u \lesssim v$ (equivalently, $v \gtrsim u$) to mean that there
  exists some constant $C > 0$ such that $u \leq C v$, and $u
  \sim v$ to mean that both $u \lesssim v$ and $u \gtrsim v$ hold
  true.
\end{notation}


\subsection{Basic setup}
\label{sec:setup}


Let $(X,\omega)$ be a compact \textde{Kähler} manifold of complex
dimension $n$ equipped with a \textde{Kähler} metric $\omega$ on $X$.
Let $\mtidlof{\vphi} := \mtidlof[X]{\vphi}$
be the multiplier ideal sheaf of the potential $\vphi$ on $X$ given at
each $x \in X$ by
\begin{equation*}
  \mtidlof{\vphi}_x := \mtidlof[X]{\vphi}_x
  :=\setd{f \in \holo_{X,x}}{
    \begin{aligned}
      &f \text{ is defined on a coord.~neighbourhood } V_x \ni x \vphantom{f^{f^f}} \\
      &\text{and }\int_{V_x} \abs f^2 e^{-\vphi} \dvol_{V_x} < +\infty
    \end{aligned}
  } \; .
\end{equation*}
A potential $\vphi$ is said to have \emph{Kawamata log-terminal (klt)
  singularities} on $X$ if $\mtidlof{\vphi} =\holo_X$ on $X$, and
\emph{log-canonical (lc) singularities} on $X$ if
$\mtidlof{(1-\eps)\vphi} =\holo_X$ for all $\eps > 0$ on $X$.

Throughout this paper,
the following data are assumed:
\begin{enumerate}[itemsep=8pt]
\item $D$ is a \emph{reduced divisor} on $X$ with \emph{simple normal
    crossings (snc)} such that $(X,D)$ is log-smooth and log-canonical
  (lc), and it is endowed with a potential $\phi_D$ defined from a
  canonical section of $D$ (see Notations
  \ref{notation:potential-definition} and \ref{notation:potentials})
  and a smooth potential $\sm\vphi_D$ such that the global function
  \begin{equation*}
    \psi_D := \phi_D -\sm\vphi_D \leq -1 \quad\text{ on } X \; ;
  \end{equation*}
  
\item $(F,e^{-\vphi_F})$ and $(M,e^{-\vphi_M})$ are holomorphic line
  bundles on $X$ equipped with singular hermitian metrics such that
  \begin{itemize}
  \item $\vphi_F$ is \emph{plurisubharmonic (psh)} and $\vphi_M$ is
    \emph{quasi-plurisubharmonic (quasi-psh)} such that the curvature
    of $\vphi_M$ is dominated by some multiple of the curvature of
    $F$, i.e.~for some constant $C >0$, one has
    \begin{equation*}
      \ibddbar\vphi_F \geq 0 \quad\text{ and }
      -C \omega \leq \ibddbar\vphi_M \leq C \ibddbar\vphi_F
      \quad\text{ on }X 
    \end{equation*}
    in the sense of $(1,1)$-currents,
    
  \item $\vphi_F$ and $\vphi_M$ both have at worst \emph{neat analytic
      singularities}, i.e.~they are locally of the form (under the
    assumption that they are both quasi-psh)
    \begin{equation*}
      c \log\paren{\sum_{j=1}^N \abs{g_{j}}^2} \mod \smooth
      \quad\text{ for some $c \in \fieldR_{\geq 0}$ and $g_{j} \in
        \holo_X \;$} \; ,
    \end{equation*}

  \item the polar sets $P_F :=\vphi_F^{-1}(-\infty)$ and $P_M
    :=\vphi_M^{-1}(-\infty)$ do not contain any \emph{lc centres of
    $(X,D)$} (i.e.~irreducible components of any intersections of
  irreducible components of $D$ in $X$),
    
  \item 
      the polar sets $P_F$ and $P_M$ are assumed to be divisors in $X$
      and $P_F \cup P_M \cup D$ has only snc; 
    
  \end{itemize}

\item $s \in \cohgp 0[X]{M}$ is a \emph{non-trivial} global
  holomorphic section of $M$ on $X$ such that
  \begin{equation*}
    \sup_X \abs s_{\vphi_M}^2 < \infty \; ;
  \end{equation*}

\item \label{item:setup-clomega}
 $X^\circ := X \setminus \paren{P_F\cup P_M}$ is a complete
  \textde{Kähler} manifold (as $P_F\cup P_M$ is an analytic set, see
  \cite{Demailly_complete-Kahler}*{Thm.~1.5}) which is equipped with a
  \emph{complete \textde{Kähler} form $\clomega$} given by
  \begin{align*}
    \clomega
    &:=2\omega +\ibddbar\frac 1{\log\abs{\ell\psi_{P_F\cup P_M}}} \\
    &=
      \begin{aligned}[t]
        &2\omega + \frac{ \ibddbar{\psi_{P_F\cup P_M}} }{
          \abs{\psi_{P_F\cup P_M}} \paren{\log\abs{\ell{\psi_{P_F\cup
                  P_M}}}}^{2}
        } \\
        &~+\paren{1 +\frac{2}{\log\abs{\ell{\psi_{P_F\cup P_M}}}}}
        \frac{
          \ibar \diff{\psi_{P_F\cup P_M}} \wedge \dbar{\psi_{P_F\cup P_M}}
        }{ \abs{\psi_{P_F\cup P_M}}^{2}
          \paren{\log\abs{\ell{\psi_{P_F\cup P_M}}}}^2
        } \; ,
      \end{aligned}
  \end{align*}
  where $P_F\cup P_M$ is viewed as a reduced divisor and
  $\psi_{P_F\cup P_M} := \phi_{P_F\cup P_M} -\sm\vphi_{P_F\cup P_M}
  \leq -1$ is defined as in Notation \ref{notation:potentials},
  therefore satisfying $\ibddbar\psi_{P_F\cup P_M} \gtrsim -\omega$ on $X$,
  and $\ell \gg e$ is some constant such that
  \begin{equation*}
    \omega + \frac{
      \ibddbar{\psi_{P_F\cup P_M}} } {\abs{\psi_{P_F\cup P_M}}
      \paren{\log\abs{\ell{\psi_{P_F\cup P_M}}}}^{2}
    } \geq 0 \; ,
  \end{equation*}
  thus having $\clomega \geq \omega$ and $\clomega \geq
  \ibar\diff\paren{\log\paren{e\log\abs{\ell\psi_{P_F\cup P_M}}}} \wedge
  \dbar\paren{\log\paren{e\log\abs{\ell\psi_{P_F\cup P_M}}}}$ on $X$.
  (Note also that the local potential of $\clomega$ can be chosen to
  be \emph{locally bounded in $X$}, i.e.~for any $p\in X$, there exist
  a neighbourhood $U \ni p$ and a bounded function $\Phi$ on $U$ such
  that $\clomega = \ibddbar \Phi$ on $U \setminus \paren{P_F\cup
    P_M}$.)
\end{enumerate}

Call an open set $V \subset X$ as an \emph{admissible open set with respect
  to the data $(\vphi_F,\vphi_M,\psi_D)$ (or simply $(\vphi_F,\psi_D)$
  when the expression of $\vphi_M$ is not under concern) in the
  holomorphic coordinate system $(z_1,\dots,z_n)$} if $V$, sitting
inside a coordinate chart with a holomorphic coordinate system
$(z_1,\dots, z_n)$ on which $F$ and $M$ are trivialised, is
biholomorphic to a polydisc centred at the origin such that 
\begin{equation*} 
  \begin{gathered}
    D\cap V =\set{z_1 \dotsm z_{\sigma_V} =0} \quad\text{ for some
      integer }
    \sigma_V \leq n \; , \\
    \res{\psi_D}_V =\underbrace{\sum_{j=1}^{\sigma_V}
      \log\abs{z_j}^2}_{=~\res{\phi_D}_V}
    ~ -\ \res{\sm\vphi_D}_V
    \quad\text{ and }\quad
    \res{\vphi_\bullet}_V = \smashoperator{\sum_{k=\sigma_V+1}^n} b_{\bullet,k}
    \log\abs{z_k}^2 +\beta_\bullet \;\;\text{ for } \bullet= F, M \; ,
  \end{gathered}
\end{equation*} 
where, after shrinking $V$ if necessary, 
\begin{itemize}
\item $\sup_V \log\abs{z_j}^2 < 0$ for $j=1,\dots, n \;$,
  
\item $\beta_F$ and $\beta_M$ are smooth functions such that
  $\sup_V\beta_\bullet \leq 0$,

\item $b_{F,k}$'s and $b_{M,k}$'s are constants such that
  $b_{\bullet,k} \geq 0$ for $k=\sigma_V+1,\dots,n \;$ (as $\vphi_F$
  and $\vphi_M$ are both quasi-psh), and

\item $\sup_V r_j \fdiff{r_j} \psi_D =2 -\inf_V r_j \fdiff{r_j} \sm\vphi_D >
  0$ for $j=1,\dots, \sigma_V$, where $r_j = \abs{z_j}$ is the radial
  component of the polar coordinates.
\end{itemize} 
Such an open set $V$ is simply called \emph{admissible} if the data
$(\vphi_F,\vphi_M,\psi_D)$ (or $(\vphi_F,\psi_D)$) are understood and
some holomorphic coordinate system satisfying the above criteria is
assumed.
Note that such admissible open sets are the kind of open sets on which
the computations in \cite{Chan&Choi_ext-with-lcv-codim-1}*{\S 2.2} and
\cite{Chan_on-L2-ext-with-lc-measures}*{\S 2.2} are valid.
The family of all such admissible open sets forms a basis of the
topology of $X$.


\subsection{$L^2$ Dolbeault isomorphism}
\label{sec:L2-Dolbeault-isom}


Put
\begin{equation*}
  \vphi := \vphi_F +\phi_D \; .
\end{equation*}
The $L^2$ Dolbeault isomorphism (see
\cite{Matsumura_injectivity}*{Prop.~5.5} and
\cite{Matsumura_injectivity-lc}*{Prop.~2.8} for a proof; see also
footnote \ref{fn:L2-Dolbeault-name} on page
\pageref{fn:L2-Dolbeault-name}) asserts that
cohomology classes on $X$ with coefficients in $K_X\otimes D\otimes F \otimes
\mtidlof{\vphi}$ can be represented by the $\dbar$-closed weighted-$L^2$
$(n,q)$-forms in the Hilbert space $\Ltwosp := \Ltwo/n,q/[X^\circ]{D\otimes
  F}_{\vphi,\clomega}$ (resp.~$\Ltwosp<\omega> := \Ltwo/n,q/[X]{D\otimes
  F}_{\vphi,\omega}$) for any $q \geq 0$, which is equipped with the
$L^2$ norm
\begin{equation*}
  \norm\zeta_{\vphi,\clomega}^2
  := \int_{X^\circ} \abs\zeta_{\vphi,\clomega}^2
  := \int_{X^\circ} \abs\zeta_{\clomega}^2 \:e^{-\vphi}
  \qquad \paren{\text{resp.~} \norm\zeta_{\vphi,\omega}^2
    := \int_{X} \abs\zeta_{\vphi,\omega}^2 
  }\; .
\end{equation*}
Note that, in the notation above, the subscript ``$\clomega$''
indicates the contraction only of the $(0,q)$-forms in $\zeta$ via
$\clomega$ and thus $\abs\zeta_{\vphi,\clomega}^2$ is a global real
$(n,n)$-form which can be integrated without any further need of
metrics on $K_X$ (see Notation \ref{notation:norm-of-nq-form}; the
same for the subscript ``$\omega$'').

One has the decomposition
\begin{equation*}
  \Ltwosp = \Harm \oplus \cl{\paren{\im\dbar}}_{\vphilist} \oplus \cl{\paren{\im\dbadj}}_{\vphilist}
  = \Harm \oplus \paren{\im\dbar}_{\vphilist} \oplus \paren{\im\dbadj}_{\vphilist} \; ,
\end{equation*}\mmark{}{Consider adding metrics into subscript place
  of $\im\dbar$. \alert{Done.}}%
where $\dbadj$ is the Hilbert space adjoint of $\dbar$ with respect to
the inner product $\iinner{\cdot}{\cdot}_{\vphi,\clomega}$, the spaces
$\paren{\im\dbar}_{\vphilist}$ and $\paren{\im\dbadj}_{\vphilist}$ are the images of the operators (with
$\cl{\paren{\im\dbar}}_{\vphilist}$ and $\cl{\paren{\im\dbadj}}_{\vphilist}$ being their closures in
$\Ltwosp$), and $\Harm$ is the space of harmonic forms (with respect
to $\dbar$ and $\dbadj$) in $\Ltwosp$.
Since $X$ is compact, both $\paren{\im\dbar}_{\vphilist}$ and $\paren{\im\dbadj}_{\vphilist}$ are closed
subspaces of $\Ltwosp$ (see, for example,
\cite{Matsumura_injectivity}*{Prop.~5.8}).
%
From the inclusion $\mtidlof{\vphi_F+\phi_D} \subset
\mtidlof{\vphi_F}$ and the $L^2$ Dolbeault isomorphism, the following
commutative diagram follows:
\begin{equation*}
  \xymatrix@R=0.65cm{
    {\cohgp q[X]{K_X\otimes D\otimes F \otimes
        \mtidlof{\vphi_F +\phi_D}}} \ar[r]
    & {\cohgp q[X]{K_X\otimes D\otimes F \otimes \mtidlof{\vphi_F}}} \\
    {\mathllap{\Harm} \isom
      \dfrac{\paren{\ker\dbar}_{\vphilist}}{\paren{\im\dbar}_{\vphilist}}} 
    \ar[u]^-\isom \ar[r] 
    & **[r]
    {\dfrac{
        \paren{\ker\dbar}_{\vphilist|\vphi_F+\sm\vphi_D|}
      }{
        \paren{\im\dbar}_{\vphilist|\vphi_F+\sm\vphi_D|}
      } 
      \quad .} \ar[u]_-\isom
  }
\end{equation*}

The above still holds true when $\clomega$ is replaced by $\omega$
(see \cite{Matsumura_injectivity-lc}*{Prop.~2.8}).
Indeed, the $L^2$ Dolbeault isomorphism implies that $\Harm<\omega>
\isom \Harm$, although they may not be the same as subsets of
$\Ltwosp$.
This fact is not needed in this paper.
It is stated here just for completeness.



\subsection{Adjoints of $\dbar$ and Bochner--Kodaira formulas}
\label{sec:BK-formulas}


In this section, let $L$ be a holomorphic line bundle on $X$
which represents either $D\otimes F$ or $D\otimes F\otimes M$ in the main
body of this article.
Suppose that there are two \emph{quasi-psh} potentials $\vphi$ and
$\clt\vphi$ on $L$ over $X$ such that $\clt\vphi$ is smooth on
$X^\circ$ and 
\begin{equation*}
  \vphi  = \clt\vphi +\psi_D
\end{equation*}
(so $\clt\vphi$ is $\vphi_F+\sm\vphi_D \:(+\vphi_M)$ while $\vphi$ is
$\vphi_F+\phi_D \:(+\vphi_M)$ in the main body of this article).
Note that $\vphi$ is singular on $X^\circ$.

Let $\dfadj_{\clt\vphi}$ be the formal adjoint of $\dbar$ with respect
to the inner product $\iinner\cdot\cdot_{\vphilist|\clt\vphi|}$
(i.e.~the adjoint of $\dbar$ on compactly supported smooth forms
$\smform*/\bullet,\bullet/\paren{X^\circ;L}$ with respect to
$\iinner\cdot\cdot_{\vphilist|\clt\vphi|}$ whose domain is extended in the
sense of currents).
Define the operator $\dfadj :=\dfadj_{\vphi}$ via the formula
\begin{equation*}
  \dfadj\zeta :=\dfadj_\vphi \zeta
  :=e^{\psi_D} \dfadj_{\clt\vphi} \paren{e^{-\psi_D} \zeta}
  =\dfadj_{\clt\vphi}\zeta +\idxup{\diff\psi_D} \ctrt
  \zeta
  \quad\text{ for } \zeta \in
  \Ltwo[X^\circ]{L}_{\vphilist} \; ,
\end{equation*}
where $\idxup{\diff\psi_D} \ctrt \cdot $ denotes the adjoint of $\dbar
\psi_D \wedge \cdot $ with respect to $\iinner{\cdot}{\cdot}_{\clomega}$
on $X^\circ$.
Call $\dfadj$ as the \emph{formal adjoint of $\dbar$ with respect to
  $\iinner\cdot\cdot_{\vphilist}$ (on $X^\circ$)}, even though it is a
priori the adjoint of $\dbar$ only on
$\smform*/\bullet,\bullet/\paren{X^\circ \setminus D ; L}$ with
respect to $\iinner\cdot\cdot_{\vphilist}$.
Indeed, it is easy to see (for example, by using the smooth
cut-off function $\theta_\eps$ described in Section
\ref{sec:outline-of-pf}) that $\smform*/\bullet,\bullet/\paren{X^\circ
  \setminus D ; L}$ is also dense in $\Ltwo[X^\circ]{L}_{\vphilist}$.
An argument for proving that $\dfadj\zeta$ is well-defined as a
current on $X^\circ$ for all $\zeta \in \Ltwo[X^\circ]{L}_{\vphilist}$
and that $\dfadj$ induces a densely defined closed operator on
$\Ltwo[X^\circ]{L}_{\vphilist}$ into
$\Ltwo/\bullet,\bullet-1/[X^\circ]{L}_{\vphilist}$ can be found in
\cite{WuXiaojun_hard-Lefschetz}*{Remark 1} (in which indeed the claim
for the $(1,0)$-Chern connection $\nabla^{(1,0)}_\vphi$ on
$\Ltwo/\bullet,0/[X^\circ]{L}_{\vphilist}$ is proved instead and the
argument works even for the case when $\psi_D$, and thus $\vphi$, has
arbitrary singularities).
As $\psi_D$ has only analytic singularities, a simple proof for the
operator $\dfadj$ on $(n,q)$-forms is sketched below for the
convenience of readers.

\begin{lemma} \label{lem:formal-adjoint-densely-defined}
  For every $\zeta \in \Ltwo/n,q/[X^\circ]{L}_{\vphilist}$,
  $\dfadj\zeta$ is a well-defined current on $X^\circ$.
  Moreover,
  $\dfadj \colon \Ltwo/n,q/[X^\circ]{L}_{\vphilist} \birat
  \Ltwo/n,q-1/[X^\circ]{L}_{\vphilist}$ is a densely defined operator
  with closed graph.
\end{lemma}
\begin{proof}[Sketch of proof]

  First it is to show that $\dfadj\zeta$ is in $\Lloc$ (unweighted) on
  $X^\circ$ for any $\zeta \in \Ltwo/n,q/[X^\circ]{L}_{\vphilist}$,
  hence a current.
  It suffices to check that for $\idxup{\diff\psi_D} \ctrt \zeta$.
  As $\psi_D$ has only analytic singularities, it is easy to check that
  $\abs{\diff\psi_D}_{\clomega}^2 \:e^{\psi_D} \in \Lloc(X^\circ)$,
  since, on any admissible open set $V \Subset X^\circ$, one has
  \begin{equation*}
    \abs{\diff\psi_D}_{\clomega}^2 \:e^{\psi_D}
    \sim \abs{\sum_{j=1}^{\sigma_V} \frac{dz_j}{z_j}
      -\diff\sm\vphi_D}_{\clomega}^2 \:\abs{z_1 \dotsm z_{\sigma_V}}^2
  \end{equation*}
  (see \cite{WuXiaojun_hard-Lefschetz}*{\S 2} for the treatment when
  $\psi_D$ having general singularities).
  Then, for any relatively open set $V \Subset X^\circ$, via the
  Cauchy--Schwarz inequality,
  one obtains
  \begin{equation*}
    \int_V \abs{\idxup{\diff\psi_D} \ctrt
      \zeta}_{\vphilist|\alert{\clt\vphi}|} \dvol_{X^\circ,\clomega}
    \leq \paren{\int_V \abs{\diff\psi_D}_{\clomega}^2 \:e^{\alert{\psi_D}}
      \dvol_{X^\circ, \clomega} }^{\frac 12}
    \paren{\int_V \abs{\zeta}_{\vphilist|\alert{\vphi}|}^2}^{\frac 12}
    \; ,
  \end{equation*}
  which implies that $\idxup{\diff\psi_D} \ctrt \zeta$ is in $\Lloc$
  on $X^\circ$, as desired.

  For the remaining claims, $\dfadj$ is densely defined since it is
  well-defined on $\smform*/n,q/\paren{X^\circ \setminus D;L}$.
  Let $\Dom\dfadj$ be the maximal domain of $\dfadj$ in
  $\Ltwo/n,q/[X^\circ]{L}_{\vphilist}$.
  To see that it has closed graph, take a sequence
  $\seq{\zeta_\nu}_{\nu\in\Nnum} \subset \Dom\dfadj \subset
  \Ltwo/n,q/[X^\circ]{L}_{\vphilist}$ such that both sequences
  $\seq{\zeta_\nu}_{\nu\in\Nnum}$ and
  $\seq{\dfadj\zeta_\nu}_{\nu\in\Nnum}$ converge in their respective
  $L^2$ spaces.
  The identity
  \begin{equation*}
    \iinner{\dfadj \zeta_\nu}{\xi}_{\vphilist}
    =\iinner{\zeta_\nu}{\dbar \xi}_{\vphilist}
    \quad\text{ for any }\;
    \xi \in \smform*/n,q-1/\paren{X^\circ \setminus D ; L}
    \;\;\text{ and } \nu \in \Nnum
  \end{equation*}
  concludes that $\lim_{\nu \tendsto \infty} \dfadj \zeta_\nu =\dfadj
  \paren{\lim_{\nu \tendsto \infty} \zeta_\nu}$. 
\end{proof}

As $\smform*/n,q-1/\paren{X^\circ \setminus D; L}$ is dense in
$\Ltwo/n,q-1/[X^\circ]{L}_{\vphilist}$ and is contained in the domain $\Dom\dbar
\subset \Ltwo/n,q-1/[X^\circ]{L}_{\vphilist}$ of $\dbar \colon
\Ltwo/n,q-1/[X^\circ]{L}_{\vphilist} \birat
\Ltwo/n,q/[X^\circ]{L}_{\vphilist}$, it follows easily from the standard
argument that
\begin{equation*}
  \Dom\dbadj \subset \Dom\dfadj \quad\text{ and }\quad
  \dbadj = \dfadj \;\;\text{ on } \Dom\dbadj \; ,
\end{equation*}
where $\Dom \dfadj \subset \Ltwo/n,q/[X^\circ]{L}_{\vphilist}$ is the
maximal domain of the operator $\dfadj$ and $\Dom \dbadj \subset
\Ltwosp$ is the domain of the Hilbert space adjoint $\dbadj$ of
$\dbar$.

The following lemma is the twisted Bochner--Kodaira formula for
$K_X\otimes L$-valued $(0,q)$-forms used in
\cite{Chan&Choi_ext-with-lcv-codim-1} with a slightly different choice
of auxiliary functions, which can be derived from \cite{Siu}*{\S 1.3}
or \cite{McNeal&Varolin_adjunction}*{Eq.~(8)}.
Recall that $\clomega$ is a \textde{Kähler} metric.
Let $\nabla_{\vphi} =\nabla_{\vphi}^{(1,0)} +\nabla^{(0,1)}$ be the
covariant differential operator (with the decomposition according to
$(1,0)$- and $(0,1)$-types) induced from the Chern connection.
Moreover, $\idxup{\ibar\Theta}\paren{\zeta, \zeta}_{\vphilist}$
denotes, for any real $(1,1)$-form $\ibar\Theta$ (usually in the form
of $\ibddbar\vphi$) and any $K_X\otimes L$-valued $(0,q)$-form
$\zeta$, the trace of the contraction between $\Theta$ and $e^{-\vphi}
\zeta \wedge \conj\zeta$ with respect to $\clomega$ on $X^\circ$ (in
the convention such that $\idxup{\ibar\Theta}
\ptinner{\zeta}{\zeta}_{\vphilist} \geq 0$ whenever $\ibar\Theta \geq
0$).
To be more precise, in local coordinates, it is given as
\begin{equation*}
  \idxup{\ibar\Theta} \paren{\zeta, \zeta}_{\vphilist}
  =\sum_{j,k} \sideset{}{'}\sum_{J_{q-1}} \Theta_{j\conj k}
  \zeta^j_{\:\conj J_{q-1}} \conj{\zeta^{k J_{q-1}}} \:e^{-\vphi} \; ,
\end{equation*}
where
\begin{gather*}
  \Theta = \sum_{j,k} \Theta_{j\conj k} \:dz^j \wedge d\conj{z^k} \; , \quad
  \zeta = \sideset{}{'}\sum_{J_q} \zeta_{\conj J_q} d\conj{z^{J_q}}
  \; , \\
  J_q = (j_1, j_2, \dots, j_q) \quad\text{and}\quad
  \sideset{}{'}\sum_{J_q} := \sum_{j_1 < j_2 < \dots < j_q} 
\end{gather*}
and the indices in the components of $\zeta$ are raised via
$\clomega$.
Under this convention, the inequality $\ibar\Theta \geq -C \clomega$
for some constant $C > 0$ implies that
$\idxup{\ibar\Theta}\ptinner{\zeta}{\zeta}_{\clomega} \geq -\pi q C
\abs\zeta_{\clomega}^2$ for any $(0,q)$-form $\zeta$ (as $\clomega =
\pi\ibar \sum_{j,k} h_{j\conj k} \:dz^j \wedge d\conj{z^k}$).

\begin{lemma} \label{lem:BK-formulas}
  Set, for any $\eps > 0$,
  \begin{equation*}
    \eta_\eps := \abs{\psi_D}^{1-\eps} \quad\text{ on } X \; .
  \end{equation*}
  Then, for any given number $\eps > 0$, the twisted Bochner--Kodaira
  formula, which is referred to as $\tBK_{\eps,\vphilist}$, becomes
  \begin{align*}
    &~
      \int_{X^\circ} \abs{\dbar\zeta}_{\vphilist}^2 \eta_\eps
      +\int_{X^\circ} \abs{\dfadj\zeta}_{\vphilist}^2 \eta_\eps
    \\ \tag*{$\tBK_{\eps,\vphilist}$}
    =&\!\!\!
       \begin{aligned}[t]
         &~\int_{X^\circ} \abs{\nabla^{(0,1)} \zeta}_{\vphilist}^2
         \eta_\eps +\int_{X^\circ}
         \idxup{\ibddbar\vphi
           +\frac{1-\eps
             }{\abs{\psi_D}} \ibddbar\psi_D }
         \ptinner\zeta\zeta_{\vphilist} \eta_\eps
         \\
         &
         \begin{termR}
           -2\paren{1-\eps} \Re
           \int_{X^\circ}
           \inner{\dfadj\zeta}{
             \frac{
               \idxup{\diff\psi_D} \ctrt \zeta
             }{
               \abs{\psi_D}
             }
           }_{\mathrlap{\vphilist}} \;\;\eta_\eps
         \end{termR}
         +
         \begin{termN}
           \eps
           \int_{X^\circ} \frac{1-\eps }{\abs{\psi_D}^2}
           \abs{\idxup{\diff\psi_D} \ctrt \zeta}_{\vphilist}^2
           \eta_\eps
         \end{termN} 
       \end{aligned}
  \end{align*}
  for any compactly supported $K_X\otimes L$-valued smooth $(0,q)$-forms
  $\zeta \in \smform*/0,q/\paren{X^\circ \alert{\setminus D}; K_X\otimes L} $ on
  $X^\circ \setminus D$.

  For the ease of reference, the untwisted Bochner--Kodaira formula is
  referred to as $\BK_{\vphilist}$ ($=\tBK_{1,\vphilist}$),
  which is given by
  \begin{equation*} \tag*{$\BK_{\vphilist}$}
    \int_{X^\circ} \abs{\dbar\zeta}_{\vphilist}^2 
    +\int_{X^\circ} \abs{\dfadj\zeta}_{\vphilist}^2
    =\int_{X^\circ} \abs{\nabla^{(0,1)} \zeta}_{\vphilist}^2
    +\int_{X^\circ}
    \idxup{\ibddbar\vphi }
    \ptinner\zeta\zeta_{\vphilist} 
  \end{equation*}
  for any $\zeta \in \smform*/0,q/\paren{X^\circ \alert{\setminus D};
    K_X\otimes L} $. 
\end{lemma}

\begin{proof}
  Notice that, as $\zeta$ is chosen to be compactly supported on
  $X^\circ \setminus D$, on which $\vphi$ and $\eta_\eps
  =\abs{\psi_D}^{1-\eps}$ are smooth, the classical Bochner--Kodaira
  formula in \cite{Siu}*{\S 1.3} or
  \cite{McNeal&Varolin_adjunction}*{Eq.~(2.2)} is applicable.
  From there, it follows that
  \begin{align*}
    &~\int_{X^\circ} \abs{\dbar\zeta}_{\vphilist}^2
      \eta_\eps +\int_{X^\circ}
      \abs{\dfadj\zeta}_{\vphilist}^2 \eta_\eps
      -
      \begin{termR}
        2\Re \int_{X^\circ} \inner{\dfadj\zeta }{
          \idxup{\diff\log\eta_\eps} \ctrt
          \zeta}_{\vphilist} \eta_\eps
      \end{termR}
      +
      \begin{termN}
        \int_{X^\circ} \abs{\idxup{\diff\log\eta_\eps}
          \ctrt \zeta}_{\vphilist}^2 \eta_\eps
      \end{termN}
    \\
    =&~
       \int_{X^\circ} \abs{\nabla^{(0,1)} \zeta}_{\vphilist}^2
       \eta_\eps
       +\int_{X^\circ} \idxup{\ibddbar\vphi
       \begin{termN}
         -\ibddbar\log\eta_\eps
       \end{termN}
       } \ptinner\zeta\zeta_{\vphilist} \eta_\eps
       \; .
  \end{align*}
  A direct computation with the choice of $\eta_\eps$ then yields the
  desired formula.
\end{proof}

\begin{remark} \label{remark:BK-BKN-compare}
  \newcommand{\wtzeta}{\widetilde \zeta}
  For the comparison with the Bochner--Kodaira--Nakano formula used in
  \cite{Matsumura_injectivity}*{Prop.~2.4} or
  \cite{Matsumura_injectivity-lc}*{Prop.~2.5}, which is derived from
  the commutator identities, let $\wtzeta$ be $\zeta$ but treated as
  an $L$-valued $(n,q)$-form.
  Using a local computation (in a normal coordinate system at a given
  point of $X^\circ$ such that $\clomega$ and $\ibar\Theta$ are
  simultaneously diagonalised and the first derivatives of the
  coefficients of $\clomega$ all vanish at that point; see, for
  example, \cite{Demailly}*{Ch.~VII, (3.2)}) and noting that $D'^* =
  -*_{\clomega} \dbar \: *_{\clomega}$, one sees that
  \begin{equation*}
    \inner{\ibar\Theta \Lambda_{\clomega} \wtzeta}{ \wtzeta}_{\vphilist}
    = \idxup{\ibar\Theta} \paren{\zeta, \zeta}_{\vphilist}
    \quad\text{ and }\quad
    \abs{D'^{*}\wtzeta}_{\vphilist}^2
    =\abs{*_{\clomega} \dbar *_{\clomega} \wtzeta}_{\vphilist}^2
    =\abs{\nabla^{(0,1)}\zeta}_{\vphilist}^2 \; .
  \end{equation*}
  The equalities on the right-hand-side can be seen more transparently
  if one notices that, in a normal coordinate system at an arbitrary
  point in $X^\circ$,
  \begin{gather*}
    \abs{\nabla^{(0,1)} \zeta}_{\vphilist}^2
    =\sum_j \sideset{}{'}\sum_{J_q} \abs{\paren{\nabla^{(0,1)}
        \zeta}_{\conj j \:\conj J_q}}^2 e^{-\vphi}
    =\sum_j \sideset{}{'}\sum_{J_q} \abs{\diff_{\conj j} \zeta_{\conj
        J_q}}^2 e^{-\vphi} \quad\text{ and } \\
    \begin{aligned}
      \abs{D'^{*}\wtzeta}_{\vphilist}^2
      &=\abs{*_{\clomega} \dbar *_{\clomega} \wtzeta}_{\vphilist}^2
      =\abs{\dbar *_{\clomega} \wtzeta}_{\vphilist}^2
      =\sum_j \sideset{}{'}\sum_{I_{n-q}} \abs{\paren{\dbar *_{\clomega}
          \wtzeta}_{I_{n-q} \conj j}}^2 e^{-\vphi} \\
      &=\sum_j \sideset{}{'}\sum_{I_{n-q}} \abs{\diff_{\conj j} \paren{*_{\clomega}
          \wtzeta}_{I_{n-q}}}^2 e^{-\vphi}
      =\sum_j \sideset{}{'}\sum_{\alert{J_q}} \abs{\diff_{\conj j}
        \wtzeta_{I_{n} \alert{\conj J_q}}}^2 e^{-\vphi}
      =\sum_j \sideset{}{'}\sum_{\alert{J_q}} \abs{\diff_{\conj j}
        \zeta_{\alert{\conj J_q}}}^2 e^{-\vphi} \; ,
    \end{aligned}
  \end{gather*}
  in which $I_{n-q} =(i_1,\dots,i_{n-q})$ and $J_q =(j_1,\dots, j_q)$
  are multi-indices such that $I_{n-q}$ is complementary to $J_q$ in
  the sense that $\set{ i_1,\dots, i_{n-q}, j_1, \dots,
    j_q } = \set{1,\dots, n}$.
  Moreover, $\zeta_{\conj J_q} =\wtzeta_{I_n \conj J_q}$.
  In the rest of this article, unless stated otherwise, $\wtzeta$ and
  $\zeta$ are identified with each other.\footnote{
    The authors' preference of the notation used in Lemma
    \ref{lem:BK-formulas} over the one used in
    \cite{Matsumura_injectivity}*{Prop.~2.4} and by many others is due
    to its better reflection of its hermitian nature of the integral
    of $\ibar\Theta$ and the ease of incorporating inequality on
    $\ibar\Theta$ in mental calculations.
    Moreover, the Chern connection $\nabla^{(0,1)}$ of type $(0,1)$
    is more apparently independent of $\vphi$ when compared to $D'^*$,
    a fact used throughout the paper. 
    Otherwise, the choice is made simply out of the habit and taste of
    the first author.
  }
\end{remark}

The pointwise formula of the Laplacian
\begin{equation} \label{eq:pointwise-BK}
  \paren{\dbar \dfadj_{\clt\vphi} +\dfadj_{\clt\vphi} \dbar}\zeta
  =-\nabla_{\clt\vphi}^{(1,0)} \cdot \nabla^{(0,1)} \zeta +\ibddbar\clt\vphi
  \Lambda_{\clomega} \zeta 
\end{equation}
for any $\zeta \in \smform/0,q/\paren{X^\circ;K_X \otimes L}$, which leads to
$\BK_{\vphilist|\clt\vphi|}$, is also stated here for later use.
Note that $-\nabla_{\clt\vphi}^{(1,0)} \cdot $ is the formal adjoint
of $\nabla^{(0,1)}$ with respect to
$\iinner\cdot\cdot_{\vphilist|\clt\vphi|}$, which is given by
\begin{equation*}
  \paren{-\nabla_{\clt\vphi}^{(1,0)} \cdot \nabla^{(0,1)}
    \zeta}_{\conj J_q}
  =-\sum_j \paren{\nabla^{(1,0)} -\paren{\diff\clt\vphi} }^{\conj j}
  \:\nabla_{\conj j} \zeta_{\conj J_q}
\end{equation*}
for any $\zeta \in \smform/0,q/\paren{X^\circ;K_X \otimes L}$ in local
coordinates, in which the index $j$ is raised via $\clomega$ and
$\nabla^{(1,0)}$ is the $(1,0)$-Chern connection on sections of $K_X$
with respect to $\clomega$.


\subsection{Refined hard Lefschetz theorem}
\label{sec:refined-hard-Lefschetz-thm}


In \cite{Matsumura_injectivity-lc}*{\S 3.1}, Matsumura proves a
refinement of the hard Lefschetz theorem with multiplier ideal
sheaves, which indeed provides a preimage for each image of the
surjective map
\begin{equation*}
  \omega^q \wedge \cdot \colon
  \cohgp 0[X]{\holoform_X^{n-q} \otimes D\otimes F\otimes \mtidlof{\vphi}} \to
  \cohgp q[X]{\logKX \otimes \mtidlof{\vphi}} \; ,
\end{equation*}
where $\vphi$ is a \emph{psh} potential on $D\otimes F$ over $X$ which
has arbitrary singularities on a Zariski closed subset while being
smooth on the complement.\footnote{
  See \cite{WuXiaojun_hard-Lefschetz}*{Thm.~2} for a related statement
  which allows more general singularities on $\vphi$.
}  
The precise statement (with a slight alteration) under the current
setting is stated as follows.

\begin{thm}[\cite{Matsumura_injectivity-lc}*{Thm.~3.3}] \label{thm:refined-hard-Lefschetz}
  Under the setup and notation given in Section \ref{sec:setup} and
  \ref{sec:L2-Dolbeault-isom}, for any harmonic $\clt u \in \Harm
  \isom \cohgp q[X]{\logKX \otimes \mtidlof{\vphi}}$, in which $\vphi$
  is psh locally everywhere in $X$ and is smooth on a Zariski open set in $X$ (and can be
  more general than $\vphi_F+\phi_D$), one has
  \begin{equation*}
    *_{\clomega} \clt u \in \cohgp 0[X]{\holoform_X^{n-q} \otimes
      D\otimes F \otimes \mtidlof{\vphi}} \; ,
  \end{equation*}
  where $*_{\clomega}$ is the Hodge $*$-operator with respect to $\clomega$.
\end{thm}

In \cite{Matsumura_injectivity-lc}*{Thm.~3.3}, the statement is proved
for $\omega$ in place of $\clomega$.
The proof of the above statement is exactly the same as the one in
\cite{Matsumura_injectivity-lc}*{\S 3.1}.
A sketch of the proof is given below.

\begin{proof}[Sketch of proof]
  Let $\omega'$ be a complete \textde{Kähler} metric on $X^\circ
  \setminus \vphi^{-1}(-\infty)$ and set, for any number $\delta > 0$,
  \begin{equation*}
    \clomega_\delta :=\clomega +\delta \omega' \; ,
  \end{equation*}
  such that $\clomega_\delta$ has a locally bounded potential in $X$
  (so that the $L^2$ Dolbeault isomorphism holds when $\clomega_\delta$
  is considered).
  For any given $\clt u \in \Harm$ and for any $\delta > 0$, there
  exists $\clt u_\delta \in \Harm<\clomega_\delta>$ which represents
  the same cohomology class as $\clt u$ in $\cohgp q[X]{\logKX \otimes
    \mtidlof{\vphi}}$ and satisfies
  \begin{equation*} \tag{$*$} \label{eq:pf-u-and-clt-u-ineq}
    \norm{*_{\delta}\clt u_\delta}_{\vphilist<\alert{\omega}>}^2
    \leq \norm{*_{\delta}\clt u_\delta}_{\vphilist<\clomega_\delta>}^2
    =\norm{\clt u_\delta}_{\vphilist<\clomega_\delta>}^2
    \leq \norm{\clt u}_{\vphilist<\clomega_\delta>}^2
    \leq \norm{\clt u}_{\vphilist}^2 \; ,
  \end{equation*}
  where $*_\delta$ is the Hodge $*$-operator with respect to
  $\clomega_\delta$.
  As $\clomega_\delta$ is complete on $X^\circ \setminus
  \vphi^{-1}(-\infty)$, the Bochner--Kodaira formula
  $\BK_{\vphilist<\clomega_\delta>}$ (see Lemma \ref{lem:BK-formulas})
  is valid for all $L^2$ sections in the domains of $\dbar$ and its
  adjoint $\dbadj$ with respect to $\vphi$ and $\clomega_\delta$.
  In particular, it can be applied to $\clt u_\delta$ and, thanks to
  the fact that $\ibddbar\vphi \geq 0$, yields
  \begin{equation*}
    \norm{\dbar *_\delta \clt u_\delta}_{\vphilist<\clomega_\delta>}^2
    =\norm{*_\delta \dbar *_\delta \clt
      u_\delta}_{\vphilist<\clomega_\delta>}^2
    =\norm{\nabla^{(0,1)} \clt
      u_\delta}_{\vphilist<\clomega_\delta>}^2
    = 0
  \end{equation*}
  (see Remark \ref{remark:BK-BKN-compare}).
  From the inequality on the far left in
  \eqref{eq:pf-u-and-clt-u-ineq} and the fact that $e^{-\vphi} 
  \gtrsim 1$ locally on $X$, it follows that $*_\delta \clt u_\delta$
  is not only holomorphic on $X^\circ \setminus \vphi^{-1}(-\infty)$,
  but also on the whole of $X$.
  It also follows that the set of holomorphic $(n-q)$-forms
  $\set{*_\delta \clt u_\delta}_{\delta \in (0,1)}$ is locally
  uniformly bounded in $X$, thus exists a subsequence $\seq{*_{\delta_\nu}
    \clt u_{\delta_\nu}}_{\nu \in\Nnum}$ which converges locally
  uniformly to some holomorphic $(n-q)$-form $f$ on $X$ as $\delta_\nu
  \tendsto 0^+$.

  For any chosen $\delta' > 0$, one has $\norm{\clt
    u_\delta}_{\vphilist<\alert{\clomega_{\delta'}}>}^2 \leq
  \norm{\clt u_\delta}_{\vphilist<\clomega_{\delta}>}^2$ for all
  $\delta$ such that $\delta \leq \delta'$.
  Together with the inequalities on the right side of
  \eqref{eq:pf-u-and-clt-u-ineq}, it follows that, by passing to a
  further subsequence if necessary, $\seq{\clt
    u_{\delta_\nu}}_{\nu\in\Nnum}$ converges weakly to some $\clt u_0$
  in $\Ltwosp<\clomega_{\delta'}>$ as $\delta_\nu \tendsto
  0^+$.
  Via the use of Cantor's diagonal argument, the subsequence $\seq{\clt
    u_{\delta_\nu}}_{\nu\in\Nnum}$ can be chosen independent of
  $\delta' >0$ as $\delta'$ shrinks to $0$, and thus the weak limits
  of the sequence in $\Ltwosp<\clomega_{\delta'}>$ for various
  $\delta' > 0$ all coincide.
  It then follows from the property of weak limits $\norm{\clt
    u_0}_{\vphilist<\clomega_{\delta'}>}^2 \leq
  \varliminf_{\delta_{\nu} \tendsto 0^+} \norm{\clt
    u_{\delta_\nu}}_{\vphilist<\clomega_{\delta'}>}^2 \leq \norm{\clt
    u}_{\vphilist}^2$ and Fatou's lemma that $\norm{\clt
    u_0}_{\vphilist}^2 \leq \norm{\clt u}_{\vphilist}^2$.
  By considering the kernel of the composition of the maps 
  \begin{equation*}
    \renewcommand{\objectstyle}{\displaystyle}
    \xymatrix@R=0.5cm{
      {\paren{\ker\dbar}_{\vphilist}} \ar@{^(->}[r]
      &
      {\paren{\ker\dbar}_{\vphilist<\clomega_{\delta'}>}} \ar@{->>}[r]
      &
      {
        \frac{
          \paren{\ker\dbar}_{\vphilist<\clomega_{\delta'}>}
        }{
          \paren{\im\dbar}_{\vphilist<\clomega_{\delta'}>}
        }
        \isom
        \frac{
          \paren{\ker\dbar}_{\vphilist}
        }{
          \paren{\im\dbar}_{\vphilist}
        }
      }
    } \; ,
  \end{equation*}
  one obtains (cf.~\cite{Matsumura_injectivity-lc}*{Prop.~3.1})
  \begin{equation*}
    \paren{\im\dbar}_{\vphilist}
    =\paren{\ker\dbar}_{\vphilist} \cap
    \paren{\im\dbar}_{\vphilist<\clomega_{\delta'}>}
    =\Ltwosp \cap
    \paren{\im\dbar}_{\vphilist<\clomega_{\delta'}>} \; .
  \end{equation*}
  Since $\clt u_{\delta_{\nu}} -\clt u \in
  \paren{\im\dbar}_{\vphilist<\clomega_{\delta'}>}$ implies that $\clt
  u_0 -\clt u \in
  \cl{\paren{\im\dbar}}_{\vphilist<\clomega_{\delta'}>}
  =\paren{\im\dbar}_{\vphilist<\clomega_{\delta'}>}$,
  and since $\clt u_0 -\clt u \in \Ltwosp$, it follows that $\clt u_0
  -\clt u \in \paren{\im\dbar}_{\vphilist}$.
  As $\clt u \in \Harm \subset \paren{\im\dbar}_{\vphilist}^\perp$, it
  forces the relations
  \begin{equation*}
    \norm{\clt u_0 -\clt u}_{\vphilist}^2 +\norm{\clt u}_{\vphilist}^2
    =\norm{\clt u_0}_{\vphilist}^2 \leq \norm{\clt u}_{\vphilist}^2 \; ,
  \end{equation*}
  and thus $\clt u_0 =\clt u$ in $\Ltwosp$.

  It remains to check that $*_{\clomega} f = (-1)^{n-q} \:\clt u$, which
  then implies $f =*_{\clomega} \clt u$ and completes the proof.
  The fact that $*_{\delta_{\nu}} \clt u_{\delta_\nu}$ converges locally
  uniformly to $f$ in $X^\circ$ implies that $\abs{*_{\clomega}
    *_{\delta_{\mathrlap{\nu}}} \clt u_{\delta_\nu} -*_{\clomega}
    f}_{\clomega}^2
  =\abs{*_{\delta_{\nu}} \clt u_{\delta_\nu} -f}_{\clomega}^2$
  converges locally uniformly to $0$ in $X^\circ$.
  A direct computation also shows that $\abs{*_{\clomega}
    *_{\delta_{\mathrlap{\nu}}} \clt u_{\delta_\nu} -*_{\delta_{\mathrlap{\nu}}}
    *_{\delta_{\mathrlap{\nu}}} \clt u_{\delta_\nu}}_{\clomega}^2$ converges
  locally uniformly to $0$ in $X^\circ$ (see the proof of
  \cite{Matsumura_injectivity-lc}*{Prop.~3.8}), which then implies
  that, for any $\delta' > 0$,
  \begin{align*}
    \abs{(-1)^{n-q} \:\clt u_{\delta_\nu} -*_{\clomega}
      f}_{\alert{\clomega_{\delta'}}}^2
    &\leq 
    \abs{(-1)^{n-q} \:\clt u_{\delta_\nu} -*_{\clomega}
      f}_{\clomega}^2 \\
    &=
    \abs{*_{\delta_{\mathrlap{\nu}}} *_{\delta_{\mathrlap{\nu}}} \clt
      u_{\delta_\nu} -*_{\clomega} f}_{\clomega}^2
    \tendsto 0 \quad\text{ locally uniformly in } X^\circ \; .
  \end{align*}
  From the weak convergence of $\clt u_{\delta_\nu}$, it follows that,
  for any $\zeta \in \smform*/n,q/\paren{X^\circ ; D\otimes F}$ and for any
  $\delta' >0$,
  \begin{equation*}
    \iinner{(-1)^{n-q} \:\clt u - *_{\clomega} f}{\zeta}_{\vphilist<\clomega_{\delta'}>}
    =\lim_{\delta_{\nu} \tendsto 0^+}
    \iinner{(-1)^{n-q} \:\clt u_{\delta_\nu} -*_{\clomega}
    f}{\zeta}_{\vphilist<\clomega_{\delta'}>}
    =0 \; .
  \end{equation*}
  Therefore, one obtains $*_{\clomega} f =(-1)^{n-q} \:\clt u$ in
  $\Ltwosp<\clomega_{\delta'}>$, and then in $\Ltwosp$ by Fatou's
  lemma, as desired.
\end{proof}


\subsection{Residue functions and $\sigma$-lc-measures}
\label{sec:residue-functions}



{
  \setDefaultMetric{\omega}

  In this section, any claims made in terms of $(F,\vphi_F)$ hold true
  when $(F,\vphi_F)$ is replaced by $(F\otimes M, \vphi_F +\vphi_M)$.

  Let $V$ be an admissible open set in $X$ with respect to
  $(\vphi_F,\psi_D)$ (see Section \ref{sec:setup}).
  Let $\symmgp_m$ be the group of permutations on a set of $m$ elements
  for any $m\in\Nnum$ and set 
  \begin{equation*}
    \cbn := \symmgp_{\sigma_V} / \symmgp_{\sigma} \times
    \symmgp_{\sigma_V-\sigma} \; ,
  \end{equation*}
  which is the set of choices of $\sigma$ elements in a set of
  $\sigma_V$ elements.
  Any element $p \in \cbn$ is abused to mean a permutation on the set of
  integer $\set{1,\dots,\sigma_V}$ such that, if $p, p' \in \cbn$ and $p
  \neq p'$, then $p\paren{\set{1,\dots, \sigma}} \neq
  p'\paren{\set{1,\dots, \sigma}}$ (and one also has
  $p\paren{\set{\sigma+1,\dots, \sigma_V}} \neq
  p'\paren{\set{\sigma+1,\dots, \sigma_V}}$).
  Then, the set of subvarieties
  \begin{equation*}
    \lcS :=\set{z_{p(1)} = z_{p(2)} = \dotsm = z_{p(\sigma)} = 0}
    \quad\text{ for } p \in \cbn
  \end{equation*}
  are precisely the set of all of the \emph{lc centres of codimension
    $\sigma$} (or \emph{$\sigma$-lc centres} for short) of $(V,D \cap
  V)$, which is denoted by $\lcc<V>$ (see
  \cite{Kollar_Sing-of-MMP}*{Def.~4.15} for the definition of lc
  centres of lc pairs; see also
  \cite{Chan_adjoint-ideal-nas}*{\S 5.2} for $\sigma$-lc centres
  in a more general setting).
  Recall that $\sect_D$ is the canonical section of $D$ such that
  $\phi_D =\log\abs{\sect_D}^2$ and that $\sect_D =z_1 \dotsm
  z_{\sigma_V}$ on the admissible open set $V$.
  Let $(r_j,\theta_j)$ be the polar coordinate system of the
  $z_j$-plane for $j=1,\dots,n$.
  Define the ad hoc notations 
  \begin{equation*}
    \smooth_{X \,*}(V)
    :=\smooth_X(V)\left[e^{\pm \cplxi \theta_1}, \dots, e^{\pm \cplxi
        \theta_n} \right] \; ,
  \end{equation*}
  which is given as a $\smooth_X(V)$-algebra, and
  \begin{equation*}
    \smform/p,q/[X \,*](V) :=\set{\text{$(p,q)$-forms on $V$ with coefficients in }
      \smooth_{X\,*}(V)} 
  \end{equation*}
  (these notations are considered only on open sets in some given
  coordinate chart, like admissible open sets, in this paper).
  Note that $e^{m\cplxi \theta_j}$ for any $m \in \Znum$ and
  $j=1,\dots,n$ can be viewed as a bounded function on $V$ (smooth in
  the variable $\theta_j$) and thus integrable on $V$.
  Apparently, one has $\smform/p,q/\paren{V} \subset
  \smform/p,q/[X \,*]\paren{V}$ for any integers $p,q \geq 0$.
  For any smooth $(n,q)$-form $f \in \smform/n,q/\paren{X;D\otimes F}$
  and for any point $x \in X$, $f$ can be written locally on some
  admissible open set $V$ centred at $x$ as a finite sum
  \begin{align} \label{eq:decomposition-f}
    \res f_V
    &=\sum_{p \in \cbn} dz_{p(1)} \wedge \dotsm \wedge dz_{p(\sigma)}
      \wedge g_p \:z_{p(\sigma+1)} \dotsm z_{p(\sigma_V)} \\
    \notag
    &=\sum_{p \in \cbn} \frac{dz_{p(1)}}{z_{p(1)}} \wedge \dotsm
      \wedge \frac{dz_{p(\sigma)}}{z_{p(\sigma)}}
      \wedge g_p \:\sect_D
  \end{align}
  for some integer $\sigma \in [0,\sigma_V]$ and $g_p \in
  \smform/n-\sigma,q/[X \,*]\paren{V}$ \footnote{
    When $n=1$ (and $V \cap D =\set 0$) and $f$ is smooth on $V$
    with $f(0)=0$, one has $f = a\: z +b\: \conj z =
    \paren{a +b \:e^{-2\cplxi \theta}} z$, where $a$ and $b$ are
    smooth $(1,1)$-forms on $V$, and the corresponding function
    $g_p=g$ is given by $\frac{dz}z \wedge g =a +b \:e^{-2\cplxi
        \theta}$.
    Also notice that (for general $n$) if, for instance, $f=\rho h$
    for some smooth function $\rho$ and holomorphic
    $(n-\sigma,q)$-form $h$ on $V$, then the $(n-\sigma,q)$-forms
    $g_p$ can be chosen to be in (the ordinary)
    $\smform/n-\sigma,q/\paren{V}$.
  } for each $p \in \cbn$ (note that the set $\cbn$, and therefore $p$,
  depends on $\sigma$).
  Let $\sigma_f :=\sigma_{V,f}$ be the \emph{minimal} $\sigma \in
  [0,\sigma_V]$ such that $\res{g_p}_{\lcS|\sigma_f|} \not\equiv 0$ for some $p
  \in\cbn[\sigma_f]$, i.e.~$\sigma_f$ is the \emph{codimension of the mlc of
    $(V,D\cap V)$ with respect to $f$} (see
  \cite{Chan&Choi_ext-with-lcv-codim-1}*{Def.~2.2.5}; see also
  \cite{Chan_adjoint-ideal-nas}*{Thm.~4.1.2 and Remark 4.1.3}).
  Let $V =U \times W$ be the decomposition into
  a product of polydiscs such that $(z_1,\dots,z_{\sigma_V})$ and
  $(z_{\sigma_V+1}, \dots, z_n)$ are coordinate systems on
  $U=U^{\sigma_V}$ and $W=W^{n-\sigma_V}$ respectively.
  To allow more general $f$ for the applications in latter sections,
  the regularity on each $g_p$ ($p \in \cbn[\sigma_f]$) is relaxed
  such that
  \begin{subequations} \label{eq:regularity-on-g_p}
    \begin{equation} \label{eq:regularity-on-g_p-smooth-U}
      \parres{\text{coef.~of } \inner{g_p}{g_{p'}}_{\alert{\vphi_F},
          \omega}}_{U \times \set w} \in
      \smooth_{U\,*}\paren{U} \quad\text{ for a.e.~} w \in W 
    \end{equation}
    (thus all high order partial derivatives with respect to the
    radial coordinates $r_1, \dots, r_{\sigma_V}$ are well-defined for
    a.e.~fixed point in $W$) and
    \begin{equation} \label{eq:regularity-on-g_p-Lloc-W} 
      \frac{\diff^\sigma}{\diff r_{\sigma} \dotsm \diff r_{2} \diff
        r_{1}} \paren{\text{coef.~of } \inner{g_p}{g_{p'}}_{\alert{\vphi_F}, \omega}}
      \in \Lloc\paren{V} \quad
      \text{ for all } \sigma =0,1,\dots,\sigma_V 
    \end{equation}
  \end{subequations}
  for all $ p,p' \in \cbn[\sigma_f]$.
 
  With a suitable choice of orientation on the normal bundle of $\lcS$
  in $X$ which determines the sign on the \textfr{Poincaré} residue map
  $\PRes[\lcS]$ (see \cite{Kollar_Sing-of-MMP}*{para.~4.18}) and with a suitable
  extension of coefficients of the map $\PRes[\lcS]$, it follows that
  \begin{equation} \label{eq:g_p-from-Poincare-residue}
    \res{g_p}_{\lcS} \:
    = \PRes[\lcS](\frac{f}{\sect_D}) 
  \end{equation}
  for all $p \in \cbn$.
  Observe that each $\res{g_p}_{\lcS}$, as a section of $K_{\lcS} \otimes
  \parres{\bigwedge^{q}\conj{\cTgt_X} \otimes F}_{\lcS}$
  on $\lcS$ (cf.~the discussion in
  \cite{Chan_adjoint-ideal-nas}*{\S 4.2}), is uniquely determined by
  $f$ up to the choice of canonical section $\sect_D$ of $D$ and can
  be extended to the $\sigma$-lc centre $\lcS*[]$ of $(X,D)$ which
  contains $\lcS$.
  Note also that $\res{g_p}_{\lcS} = 0$ for all $p \in \cbn$ when
  $\sigma > \sigma_f$, as can be seen from \eqref{eq:decomposition-f}.

  In \cite{Chan&Choi_ext-with-lcv-codim-1}*{Prop.~2.2.1}, the following
  theorem is proved for $(n,0)$-forms via successive use of
  integration by parts.
  The proof for $(n,q)$-forms is essentially the same, which is given
  below.
  
  \begin{thm}[cf.~\cite{Chan&Choi_ext-with-lcv-codim-1}*{Prop.~2.2.1}]
    \label{thm:residue-fcts-and-norms}
    Recall that $\vphi =\vphi_F+\phi_D$.
    Let $\rho$ be a compactly supported smooth function on an admissible
    open set $V \subset X$ and let $f \in \Ltwosp/q/$ which admits the
    decomposition \eqref{eq:decomposition-f} with the coefficients
    satisfying \eqref{eq:regularity-on-g_p} on some neighbourhood of
    $\cl V$.
    Then, for any $\eps > 0$,
    \begin{equation*}
      \eps \int_V \frac{\rho \abs
        f_{\vphilist}^2}{\abs{\psi_D}^{\sigma+\eps}}
      =
      \begin{dcases}
        \frac{(-1)^{\sigma_f} \:\eps}{\prod_{j=1}^{\sigma_f} (\sigma -j+\eps)}
        \int_V \frac{G_{\sigma_f}}{\abs{\psi_D}^{\sigma -\sigma_f+\eps}}
        & \text{ when } \sigma \geq \sigma_f \; ,
        \\
        +\infty & \text{ when } 0 \leq \sigma < \sigma_f \text{ and }
        \eps < \sigma_f -\sigma \; ,
      \end{dcases}
    \end{equation*}
    where $G_{\sigma_f}$ is an $(n,n)$-form on $V$ independent of $\eps$
    with $L^1$ coefficients, which contain derivatives of $\rho \abs
    f_{\vphilist|\vphi_F|}^2$ of order at most $\sigma_f$ in the normal
    directions of the lc centres $\lcS$.
    Moreover, when $\sigma \geq \sigma_f$, one obtains the residue
    norm of $f$ on $\lcc<V>$ with respect to the $\sigma$-lc-measure,
    given by
    \begin{equation*}
      \lim_{\eps \tendsto 0^+} \eps \int_V \frac{\rho \abs
        f_{\vphilist}^2}{\abs{\psi_D}^{\sigma+\eps}} 
      =
      \begin{dcases}
        0 & \text{ when } \sigma > \sigma_f \text{ or }
        \sigma=\sigma_f=0 \; , \\
        \sum_{p \in \cbn} \frac{\pi^\sigma}{(\sigma -1)!} \int_{\lcS}
        \rho \abs{g_p}_{\vphilist|\vphi_F|}^2
        & \text{ when } \sigma \in [\sigma_f, \sigma_V] \cap \Nnum
        \; .
      \end{dcases}
    \end{equation*}
    (In particular, $\parres{\rho g_p}_{\lcS} \equiv 0$ for $\sigma >
    \sigma_f$, with the knowledge that the coefficients $g_p$ in
    \eqref{eq:decomposition-f} depend on $\sigma$.)
  \end{thm}

  \begin{remark} \label{rem:residue-with-clomega}
    As can be observed from the proof, Theorem
    \ref{thm:residue-fcts-and-norms} still holds true when $\omega$ in
    the claim and also in the assumptions in
    \eqref{eq:regularity-on-g_p} is replaced by $\clomega$.
  \end{remark}
  
  \begin{remark}
    The results in \cite{Chan_on-L2-ext-with-lc-measures}*{Prop.~2.2.1
      and Cor.~2.3.3}, which make similar claims as in Theorem
    \ref{thm:residue-fcts-and-norms} but for $(n,0)$-forms $f$ and for
    the residue function
    \begin{equation*}
      \RTF[\rho]|f|(\eps)[\sigma] :=\eps \int_V \frac{
        \rho \abs f_{\vphilist}^2
      }{\logpole}
    \end{equation*}
    instead of $\eps \int_V \frac{\rho \abs
      f_{\vphilist}^2}{\abs{\psi_D}^{\sigma+\eps}}$, can also be
    extended similarly to the current setting (where $f$ is an
    $(n,q)$-form with more relaxed regularity assumptions).
    Following the proof of
    \cite{Chan_on-L2-ext-with-lc-measures}*{Thm.~2.3.1}, one can also
    show that $\eps \mapsto \eps \int_V \frac{\rho \abs
      f_{\vphilist}^2}{\abs{\psi_D}^{\sigma+\eps}}$ can be continued
    analytically across the origin $\eps=0$, but it is not a priori an
    entire function as there are apparent poles at $\eps = -(\sigma
    -j)$ for $j=1, \dots, \sigma_f-1$ and also at $\eps
    =-(\sigma-\sigma_f)$ when $\sigma > \sigma_f$, as can be seen from
    the identity in Theorem \ref{thm:residue-fcts-and-norms}.
    For the purpose of this paper, it suffices to consider $\eps
    \int_V \frac{\rho \abs
      f_{\vphilist}^2}{\abs{\psi_D}^{\sigma+\eps}}$ instead of
    $\RTF[\rho]|f|(\eps)[\sigma]$.
  \end{remark}

  \begin{proof}
    Write
    \begin{equation*}
      g_p
      =\sgn{p} \:dz_{p(\sigma+1)} \wedge \dotsm \wedge
      dz_{p(\sigma_V)} \wedge dz_{\sigma_V+1} \wedge \dotsm \wedge
      dz_n \wedge \rs g_p \; ,
    \end{equation*}
    where $\rs g_p$ is a $(0,q)$-form and $\sgn{p}$ is
    the sign of the permutation representing the choice $p$.
    It follows from the decomposition \eqref{eq:decomposition-f} (with
    $\sigma=\sigma_f$) of $f$ that (note also Notation
    \ref{notation:norm-of-nq-form})
    \begin{equation} \label{eq:ptnorm-f-expansion}
      \rho\abs{f}_{\vphilist}^2
      =
        \smashoperator{\sum_{p,p' \in \cbn[\sigma_f]}} \;
        \overbrace{
          \rho \inner{\rs g_p}{\rs g_{p'}}_{\vphilist F}
        }^{F_{p,p'} \::=}
        \:
        \frac{
          \bigwedge_{j=1}^{\sigma_V} \pi\ibar \:dz_j \wedge
          d\conj{z_j}
        }{
          \prod_{j=1}^{\sigma_f} z_{p(j)} \conj{z_{p'(j)}}
        }
        \wedge \:\smashoperator{\bigwedge_{k=\sigma_V+1}^n}
        \pi\ibar\:dz_k \wedge d\conj{z_k} 
    \end{equation}
    and, following the notation in the proof of
    \cite{Chan_on-L2-ext-with-lc-measures}*{Prop.~2.2.1}, one can
    write accordingly
    \begin{equation*} \tag{$*$} \label{eq:RTF-decomposition}
      \eps \int_V \frac{\rho \abs
        f_{\vphilist}^2}{\abs{\psi_D}^{\sigma+\eps}}
      =\smashoperator[r]{\sum_{p,p' \in \cbn[\sigma_f]}} \: \RTI[F_{p,p'}](\eps)[\sigma] \; ,
    \end{equation*}
    where $\RTI[F_{p,p'}](\eps)[\sigma]$ is the summand containing the
    term $F_{p,p'}$.
    Let $(r_j,\theta_j)$ be the polar coordinate system in the
    $z_j$-plane.
    That $\inner{g_p}{g_{p'}}_{\vphilist F}$
    (more precisely, $\inner{\rs g_p}{\rs g_{p'}}_{\vphilist F}$)
    satisfying \eqref{eq:regularity-on-g_p} implies
    that the conditions in \eqref{eq:regularity-on-g_p} hold true with
    $F_{p,p'}$ in place of $\inner{g_p}{g_{p'}}_{\vphilist F}$ for any
    $p,p'\in \cbn[\sigma_f]$.
    Integration by parts with respect to the variables $r_1, \dots,
    r_{\sigma_V}$ can therefore be applied to the integral
    $\RTI[F_{p,p'}](\eps)[\sigma]$ without questions.
    The rest of the proof is proceeded as in
    \cite{Chan&Choi_ext-with-lcv-codim-1}*{Prop.~2.2.1} or
    \cite{Chan_on-L2-ext-with-lc-measures}*{Prop.~2.2.1}. 

    First consider the case when $\sigma \geq \sigma_f$.
    To compute the summands in \eqref{eq:RTF-decomposition} with
    $p=p'$, and it suffices to consider only the case where
    $p=p'=\id$, the identity permutation.
    Set 
    \begin{align*}
      F_0 &:=F_{\id,\id; 0} :=F_{\id,\id} \quad\text{ and } \\
      F_j &:=F_{\id,\id; j} :=\fdiff{r_j}
            \paren{\frac{F_{j-1}}{r_j^2 \fdiff{r_j^2}\psi_D}}
      =\fdiff{r_j} \paren{\frac{F_{j-1}}{1-\frac{r_j}{2}
          \fdiff{r_j} \sm\vphi_D}}
      \quad\text{ for } j=1,\dots,\sigma_V
    \end{align*}
    (when $F_{p,p';0} := F_{p,p'}$ is considered for general $p, p'
    \in\cbn[\sigma_f]$, all $r_j$ in the formula of $F_j
    (:=F_{p,p';j})$ should be replaced by $r_{p(j)}$).
    Note that $\frac{1}{r_j^2 \fdiff{r_j^2}\psi_D} > 0$ and is
    smooth on $V$ by the choice of the admissible set $V$ (see Section
    \ref{sec:setup}).
    In view of Fubini's theorem, it follows that
    \begin{align*}
      \RTI[F_{0}](\eps)[\sigma]
      &=\eps \int_{V} \frac{F_0}{\abs{\psi_D}^{\sigma+\eps}}
        \bigwedge_{j=1}^{\sigma_f}
        \frac{\pi\ibar \:dz_j \wedge d\conj{z_j}} {\abs{z_j}^2}
        \underbrace{
          \wedge \smashoperator{\bigwedge_{j=\sigma_f+1}^{\sigma_V}}
          \pi\ibar \:dz_j \wedge d\conj{z_j}
          \wedge \smashoperator{\bigwedge_{k=\sigma_V+1}^n}
          \pi\ibar \:dz_k \wedge d\conj{z_k}
        }_{\text{(made implicit)}} \\
      &=\eps \int_{V} \frac{F_0}{\abs{\psi_D}^{\sigma+\eps}}
        \prod_{j=1}^{\sigma_f} d\log r_j^2 
        \cdot \underbrace{\prod_{j=1}^{\sigma_f}
        \frac{d\theta_j}{2}}_{=:~\vect{d\theta}}
        =\eps \int_{V} \frac{-F_0}{r_1^2 \fdiff{r_1^2}\psi_D}
        \frac{d\abs{\psi_D}}{\abs{\psi_D}^{\sigma+\eps}}
        \prod_{j=\alert{2}}^{\sigma_f} d\log r_j^2 
        \cdot \vect{d\theta} \\
      &=\frac\eps{\sigma-1+\eps}
        \int_{V} \frac{F_0}{r_1^2 \fdiff{r_1^2}\psi_D}
        \:d\paren{\frac{1}{\abs{\psi_D}^{\sigma-1+\eps}}}
        \prod_{j=2}^{\sigma_f} d\log r_j^2 
        \cdot \vect{d\theta} \\
      &\overset{\mathclap{\text{int.~by parts}}}=
        \quad
        \frac{-\eps}{\sigma-1+\eps}
        \int_{V}
        \frac{\alert{F_1}}{\abs{\psi_D}^{\sigma-1+\eps}}
        \prod_{j=2}^{\sigma_f} d\log r_j^2 
        \cdot dr_1 \:\vect{d\theta} \\
      \label{eq:residue-fct-identity-RHS} \tag{$**$}
      &= \dotsm =
        \frac{(-1)^{\sigma_f} \:\eps} {\prod_{j=1}^{\sigma_f} \paren{\sigma-j+\eps}} 
        \int_{V}
        \frac{F_{\alert{\sigma_f}}}{\abs{\psi_D}^{\sigma-\sigma_f+\eps}}
        \prod_{j=1}^{\sigma_f} dr_j
        \cdot \vect{d\theta} \; .
    \end{align*}
    The boundary terms arising from integration by parts all vanish as
    one has $\res{\frac{1}{\abs{\psi_D}}}_{\mathrlap{\set{r_j=0}}}
    \quad = 0$ and $\res{F_\sigma}_{\set{r_j=1}} = 0$ for all
    $j=1, \dots, \sigma_V$ and $\sigma=0, \dots, \sigma_V$.
    The assumption \eqref{eq:regularity-on-g_p-Lloc-W} guarantees that
    $F_{\sigma_f}$ is integrable and so is the integral in
    \eqref{eq:residue-fct-identity-RHS} (note that
    $\abs{\psi_D}^{\sigma-\sigma_f+\eps} \geq 1$).
    This implies that all equalities above which invokes integration
    by parts are valid and the integral $\RTI[F_0](\eps)[\sigma]$ is
    convergent.
    Moreover, when $\sigma > \sigma_f$, as the integral in
    \eqref{eq:residue-fct-identity-RHS} remains convergent when
    $\eps=0$ and the coefficient of the integral is a multiple of
    $\eps$, it follows that
    \begin{equation*}
      \lim_{\eps \tendsto 0^+} \RTI[F_0](\eps)[\sigma] = 0
      \quad\text{ when } \sigma > \sigma_f \; .
    \end{equation*}
    When $\sigma =\sigma_f$, the fundamental theorem of calculus
    (which makes use of the assumption
    \eqref{eq:regularity-on-g_p-smooth-U}) yields
    \begin{align*}
      \lim_{\eps \tendsto 0^+} \RTI[F_0](\eps)[\sigma]
      &=\frac{(-1)^{\sigma}} {\paren{\sigma-1}!} 
        \int_{V} F_{\sigma}
        \prod_{j=1}^{\sigma} dr_j
        \cdot \vect{d\theta}
        =\frac{(-1)^{\sigma}} {\paren{\sigma-1}!} 
        \int_{V} \fdiff{r_\sigma} \paren{
          \frac{F_{\sigma-1}}{1-\frac{r_\sigma}{2} \fdiff{r_\sigma}\sm\vphi_D}
        }
        \prod_{j=1}^{\sigma} dr_j
        \cdot \vect{d\theta}
      \\
      &=\frac{(-1)^{\sigma-1}} {\paren{\sigma-1}!} 
        \int_{\mathrlap{\set{r_\sigma=0}}} \quad F_{\sigma-1}
        \prod_{j=1}^{\sigma-1} dr_j
        \cdot \vect{d\theta} \\
      &=\dotsm
        =\frac{1} {\paren{\sigma-1}!} 
        \;\;\int\limits_{\mathclap{\set{r_1=\dotsm=r_\sigma=0}}}
        F_{0}
        \:\alert{\vect{d\theta}}
        =\frac{\alert{\pi^\sigma}}{(\sigma-1)!} \int_{\lcS[\id]}
        \rho\abs{g_{\id}}_{\vphilist F}^2 \; .
    \end{align*}
    More generally, one has
    \begin{equation*}
      \lim_{\eps \tendsto 0^+} \RTI[F_{p,p}](\eps)[\sigma]
      =\frac{\alert{\pi^\sigma}}{(\sigma-1)!} \int_{\lcS}
      \rho\abs{g_p}_{\vphilist F}^2
      \quad\text{ for any } p \in \cbn[\sigma_f] \; .
    \end{equation*}
    Note that the above computation is still valid if the expansion of
    $\rho\abs f_{\vphilist}^2$ in \eqref{eq:ptnorm-f-expansion} is
    made into a sum of $p,p' \in \cbn[\alert{\sigma}]$ (instead of
    $\cbn[\sigma_f]$, thus coming with a different set of
    $\set{g_p}_p$) to start with.
    The above equation thus holds true for any $p \in
    \cbn[\alert{\sigma}]$ and all integers $\sigma$ such that
    $\sigma_V \geq \sigma \geq \sigma_f$.

    To compute the ``cross terms'' in \eqref{eq:RTF-decomposition},
    i.e.~summands $\RTI[F_{p,p'}](\eps)[\sigma]$ with $p \neq p'$,
    consider the special case where $p(j) = p'(j)$ for $j=1,\dots,
    \sigma_f-1$ but $p(\sigma_f) \neq p'(\sigma_f)$.
    In this case, one has
    \begin{equation*}
      \RTI[F_{p,p'}](\eps)[\sigma]
      =\eps \int_V \frac{F_{p,p'}}{\abs{\psi_D}^{\sigma+\eps}}
      \bigwedge_{j=1}^{\sigma_f-1} \frac{
        \pi\ibar \:dz_{p(j)} \wedge d\conj{z_{p(j)}}
      }{\abs{z_{p(j)}}^2}
      \underbrace{
        \wedge
        \smashoperator{\bigwedge_{j=\sigma_f}^{\sigma_V}} \frac{
          \pi\ibar \:dz_{p(j)} \wedge d\conj{z_{p(j)}}
        }{z_{p(\sigma_f)} \:\conj{z_{p'(\sigma_f)}}}
        \wedge
        \smashoperator{\bigwedge_{k=\sigma_V+1}^n}
        \pi\ibar \:dz_{k} \wedge d\conj{z_{k}}
      }_{\text{integrable}} \; .
    \end{equation*}
    In view of Fubini's theorem and the assumptions in
    \eqref{eq:regularity-on-g_p}, the computation for the case where
    $p=p' \in\cbn[\sigma_f]$ still applies and thus (with the same
    notation as before and making the irrelevant variables implicit)
    \begin{equation*}
      \RTI[F_{p,p'}](\eps)[\sigma]
      =\frac{(-1)^{\sigma_f-1} \:\eps}{
        \prod_{j=1}^{\sigma_f-1} (\sigma -j+\eps)
      } \int_V
      \frac{F_{p,p';\:\sigma_f-1}}{\abs{\psi_D}^{\sigma -\sigma_f+1+\eps}}
      \smashoperator{\prod_{j=1}^{\sigma_f-1}} dr_{p(j)} \cdot
      \frac{dr_{p(\sigma_f)}^2}{z_{p(\sigma_f)}} \:\vect{d\theta} 
    \end{equation*}
    (one can further lower the exponent on $\abs{\psi_D}$ by writing
    $\frac{F_{p,p'; \:\sigma_f -1}}{z_{p(\sigma_f)}}$ as
    $\frac{\conj{z_{p(\sigma_f)}} F_{p,p'; \:\sigma_f -1}
    }{r_{p(\sigma_f)}^2}$ and applying integration by parts with
    respect to the variable $r_{p(\sigma_f)}$ as before if one wants
    to fit this in the identity in the claim).
    As the integral is convergent even when $\eps =0$ and the
    coefficient of the integral is a multiple of $\eps$, one obtains
    \begin{equation*}
      \lim_{\eps \tendsto 0^+} \RTI[F_{p,p'}](\eps)[\sigma] = 0 \; .
    \end{equation*}
    Note that this holds true even for $\sigma > \sigma_f -1$.
    The computation for other summands $\RTI[F_{p,p'}](\eps)[\sigma]$
    with $p \neq p'$ is similar.
    
    Combining the above results, the claims for the case $\sigma \geq
    \sigma_f$ are proved.

    For the case $0 \leq \sigma < \sigma_f$ (which is to show divergence of the
    integral in question), as $\frac{1}{\abs{\psi_D}^{\sigma}} \leq
    \frac{1}{\abs{\psi_D}^{\sigma'}}$ for any $\sigma' \leq \sigma$, it
    suffices to consider $\sigma$ such that $\sigma_f -1 < \sigma <
    \sigma_f$.
    The above computation shows that the ``cross terms''
    $\RTI[F_{p,p'}](\eps)[\sigma]$ in \eqref{eq:ptnorm-f-expansion}
    are convergent, it remains to show that
    $\RTI[F_{p,p}](\eps)[\sigma] =+\infty$ for some $p
    \in\cbn[\sigma_f]$. 
    Without loss of generality, assume that
    $\res{F_{\id,\id}}_{\lcS|\sigma_f|} \not\equiv 0$ 
    (such assumption is valid by the definition of $\sigma_f$).
    Consider a change of coordinates on the
    admissible open set $V$ which changes the radial coordinates
    $(r_1,\dots,r_{\sigma_f})$ to $(\abs{\psi_D}, q_2, \dots,
    q_{\sigma_f})$, where
    \begin{equation*}
      q_j := \frac{\log r_j^2}{\psi_D}
      =\frac{\abs{\log r_j^2}}{\abs{\psi_D}}
      \quad\text{ for } j=2,\dots, \sigma_f \; ,
    \end{equation*}
    in which each $q_j$ varies within $[0,1]$ on $V$.
    Then, with the same notation as before and making the irrelevant
    variables implicit, one has, for $\eps < \sigma_f -\sigma$, 
    \begin{align*}
      \RTI[F_{\id,\id}](\eps)[\sigma]
      &=\eps \int_V \frac{
        F_{\id,\id}
        }{r_1^2 \fdiff{r_1^2} \psi_D}
        \frac{
          \abs{\psi_D}^{\sigma_f -1} \:d\abs{\psi_D}
        }{
          \abs{\psi_D}^{\sigma+\eps}
        } \prod_{j=2}^{\sigma_f} dq_j \cdot \vect{d\theta} \\
      &=\frac\eps{\sigma_f -\sigma -\eps} \int_V \frac{
        F_{\id,\id}
        }{r_1^2 \fdiff{r_1^2} \psi_D}
        \:d\paren{\abs{\psi_D}^{\sigma_f -\sigma -\eps}}
        \prod_{j=2}^{\sigma_f} dq_j \cdot \vect{d\theta} \; .
    \end{align*}
    Since $ \frac{F_{\id,\id}}{r_1^2 \fdiff{r_1^2} \psi_D} > 0$ on
    some open subset $V' \subset V$ such that $V' \cap \lcS|\sigma_f| \neq
    \emptyset$ and $d\paren{\abs{\psi_D}^{\sigma_f -\sigma -\eps}}$ is
    not integrable on $V'$, the above integral diverges.
    This completes the proof for the case $0 \leq \sigma < \sigma_f$.
  \end{proof}

  \begin{remark}
    When $\psi_D$ is replaced by a global function $\psi \leq -1$ of
    the form
    \begin{equation*}
      \res\psi_{V} = \sum_{j = 1}^{\sigma_V} \nu_j \log\abs{z_j}^2 
      +\smashoperator{\sum_{k=\sigma_V+1}^n} c_k \log\abs{z_k}^2 + \alpha 
    \end{equation*}
    where $\nu_j > 0$ for $j=1,\dots, \sigma_V$ and $c_k \geq 0$ for
    $k=\sigma_V, \dots, n$ are constants and $\alpha \in \smooth_X(V)$
    (as in \cite{Chan&Choi_ext-with-lcv-codim-1}*{Prop.~2.2.1}), the
    claims in Theorem \ref{thm:residue-fcts-and-norms} still hold true
    when the residue norm is given by
    \begin{equation*}
      \lim_{\eps \tendsto 0^+} \eps \int_V \frac{\rho \abs
        f_{\vphilist}^2}{\abs{\psi}^{\sigma+\eps}}
      =
      \begin{dcases}
        0 & \text{ when } \sigma > \sigma_f \text{ or }
          \sigma=\sigma_f=0 \; , \\
        \sum_{p \in \cbn} \frac{\pi^\sigma}{(\sigma -1)!
          \:\alert{\vect{\nu}_p}} \int_{\lcS} 
        \rho \abs{g_p}_{\vphilist|\vphi_F|}^2
        & \text{ when } \sigma \in [\sigma_f, \sigma_V] \cap \Nnum
        \; ,
      \end{dcases}
    \end{equation*}
    where $\vect{\nu}_p :=\prod_{j=1}^{\sigma} \nu_{p(j)}$, which can
    be checked easily by following the proof of Theorem
    \ref{thm:residue-fcts-and-norms}.
  \end{remark}

}


\section{Proof of Theorem \ref{thm:main-result}}
\label{sec:proof-main-result}

\subsection{Outline of the proof}
\label{sec:outline-of-pf}


An outline of the proof of Theorem \ref{thm:main-result} is given in
this section, with some detailed verifications for \emph{the smooth
  case} (i.e.~$\vphi_F$ and $\vphi_M$ being smooth and thus $X
=X^\circ$ and $\omega =\clomega$) given in Section
\ref{sec:pf-smooth-vphi_FM}.
Extra treatments to \emph{the singular case} (i.e.~$\vphi_F$ and $\vphi_M$
having neat analytic singularities as described in Section
\ref{sec:setup}) are given in Section \ref{sec:pf-singular-vphi_FM}.

Recall that $\vphi :=\vphi_F +\phi_D$.
Let $\sect_D$ be a canonical section of $D$ such that $\phi_D
=\log\abs{\sect_D}^2$.
Let $\theta \colon [0,\infty) \to [0,1]$ be a smooth non-decreasing
cut-off function such that $\res\theta_{[0,\frac 12]} \equiv 0$ and
$\res\theta_{[1,\infty)} \equiv 1$.
For $\eps \geq 0$, set $\theta_\eps := \theta \circ \frac
1{\abs{\psi_D}^\eps}$ and $\theta'_\eps := \theta' \circ \frac
1{\abs{\psi_D}^\eps}$ for convenience (where $\theta'$ is the
derivative of $\theta$).
Note that both $\theta_{\eps}$ and $\theta'_{\eps}$ have compact
supports inside $X \setminus D$ for $\eps > 0$.
One also has $\theta_\eps \nearrow 1$ pointwisely on $X\setminus D$ as
$\eps \searrow 0$.

\begin{enumerate}[label=\textbf{Step \Roman*:}, ref=\Roman*,
  leftmargin=0pt, labelsep=0pt, itemindent=*, align=left, itemsep=1.5ex]

\item \label{item:pick-u-decompose-v}
  Note that the space of harmonic forms $\Harm$ is isomorphic to
  $\cohgp q[X]{\logKX \otimes \mtidlof{\vphi}}$ (see Section
  \ref{sec:L2-Dolbeault-isom}).
  The homomorphism
  \begin{equation*}
    \iota_0 \colon
    \cohgp q[X]{\logKX \otimes \mtidlof{\vphi_F+\phi_D}} \to
    \cohgp q[X]{\logKX \otimes \mtidlof{\vphi_F}}
  \end{equation*}
  can be viewed as the homomorphism
  \begin{equation*}
    \iota_0 \colon \Harm \to
    \cohgp q[X]{\logKX \otimes \mtidlof{\vphi_F}} \; . 
  \end{equation*}
  Choose any $\clt u \in \paren{\ker\iota_0}^\perp 
  \subset \ker\iota_0 \oplus \paren{\ker\iota_0}^\perp =\Harm$
  such that $\clt u \in \ker\mu_0$, i.e.~$s\clt u$ represents the zero
  class in $\cohgp q[X]{\logKX M \otimes \mtidlof{\vphi_F+\vphi_M}}$.
  Then $s\clt u =\dbar \clt v$ for some $\clt v \in
  \Ltwosp/q-1/M|\vphi_F+\vphi_M+\sm\vphi_D|$ which can be decomposed into \mmark{$\clt v =
  \clt v_{(2)} + v_{(\infty)}$}{Need proof. \alert{Done.}} such that
  \begin{align*}
    \clt v_{(2)} &\in \Ltwo/n,q-1/[X^\circ]{D\otimes F\otimes
              M}_{\vphi_F+\vphi_M+\alert{\phi_D},\clomega}  \; , \\
    v_{(\infty)} &\in \smform/n,q-1/\paren{X; D\otimes F\otimes M
                   \otimes \mtidlof{\vphi_F+\vphi_M}} 
  \end{align*}
  (see Lemma \ref{lem:regularity-of-v} and Remark
  \ref{rem:regularity-v-for-sing-F-M}).
  The goal is to show that $\clt u = 0$, which will then imply that
  $\ker\mu_0 =\ker \iota_0$.
  
\item \label{item:results-from-harmonic-u-and-BK}
  Since $\vphi_F+\phi_D$ is psh and $\clt u$ is harmonic in the
  corresponding $L^2$ space, it follows from the refinement of the
  hard Lefschetz theorem proved in
  \cite{Matsumura_injectivity-lc}*{Thm.~3.3} (see also Theorem
  \ref{thm:refined-hard-Lefschetz}) that \mmark{$*_{\clomega} \clt
  u$ is holomorphic on $X^\circ$ (which also implies that $\frac{\clt
    u}{\sect_D}$ is smooth on $X^\circ$).
  One can derive from this fact}{Need proof, esp.~on why
  $*_{\clomega}\clt u$ is holomorphic when $P_F\cup P_M \neq \emptyset$
  and on how poles of $\clt u$
    along $P_F\cup P_M$ is controlled which leads to the claim on
    residue. \alert{Done.}} and
  Theorem \ref{thm:residue-fcts-and-norms} (or the computation
  corresponding to the $1$-lc-measure in 
  \cite{Chan&Choi_ext-with-lcv-codim-1}*{Prop.~2.2}) that
  \begin{equation*}
    \lim_{\eps \tendsto 0^+} \eps \int_{\mathrlap{X^\circ}\;\;}
    \frac{
      \abs{\dbar\psi_D \otimes \clt u}_{\vphilist}^2
    }{\abs{\psi_D}^{1+\eps}}
    < +\infty
  \end{equation*}
  and thus
  \begin{equation} \label{eq:limit-ctrt-u=0}
    \lim_{\eps \tendsto 0^+} \alert{\eps^2} \int_{\mathrlap{X^\circ}\;\;}
    \frac{
      \abs{\dbar\psi_D \otimes \clt u}_{\vphilist}^2
    }{\abs{\psi_D}^{1+\eps}} = 0
  \end{equation}
  (see Proposition \ref{prop:residue-of-dpsi-ctrt-u} for the smooth
  case and Proposition \ref{prop:residue-of-dpsi-ctrt-clt-u} for the
  singular case).
  Via a \mmark{careful use of the Bochner--Kodaira formula}{Need
    proof. \alert{Done.}} $\BK_{\vphilist}$ on $X^\circ$ in Lemma
  \ref{lem:BK-formulas} (in which the fact that the \text{Kähler}
  metric $\clomega$ being smooth along the general points of $D$ while
  the involving potential $\phi_D$ being singular along $D$ is taken
  care of using \eqref{eq:limit-ctrt-u=0}) and the fact that
  $\ibddbar\vphi_F \geq 0$ on $X$, \mmark{one obtains}{Need
    proof. (Validity of $\dfadj =\dfadj_{\vphi_F} +\idxup{\diff\phi_D}
    \ctrt $ on $X^\circ$ may be explained in the preliminaries.) \alert{Done.}} 
  \begin{equation} \label{eq:nable-u=0_ibddbar-vphi_F-u=0}
    \nabla^{(0,1)}\clt u = 0 \quad\text{ and }\quad
    \idxup{\ibddbar\vphi_F}\ptinner{\clt u}{\clt
      u}_{\clomega} = 0 
    \quad\text{ on } X^\circ \; ,
  \end{equation}
  (see Proposition \ref{prop:nabla-u=0_curv-u=0}, supplemented by
  Lemma \ref{lem:limits-of-d-of-cutoff-fcts} for the singular case)
  which consequently lead to
  \begin{equation} \label{eq:dfadj-su=0}
    \dfadj_{\vphi_M} \paren{s\clt u} = 0 \quad\text{ and }\quad
    s\clt u \in \Dom \dbadj_{\vphi_M} \; ,
  \end{equation}
  where $\dbadj_{\vphi_M}$ is the Hilbert space adjoint of $\dbar$ on
  the $L^2$ space $\Ltwosp/\bullet/M|\vphi+\vphi_M|$ and
  $\dfadj_{\vphi_M}$ is the corresponding formal adjoint (see
  Corollary \ref{cor:dfadj_M-su=0}, supplemented by Lemma
  \ref{lem:limits-of-d-of-cutoff-fcts} for the singular case).
  
  
\item \label{item:by-parts-on-norm-su}
  Now consider
  \begin{equation*}
    \norm{s\clt u}_{\vphilist M}^2
    =\iinner{s\clt u}{\dbar\clt v}_{\vphilist M}
    =\iinner{s\clt u}{\dbar\clt v_{(2)}}_{\vphilist M}
    +\iinner{s\clt u}{\dbar v_{(\infty)}}_{\vphilist M}
  \end{equation*}
  and apply integration by parts.
  The form $\clt v_{(2)}$ being $L^2$ with respect to the potential
  $\vphi_F+\phi_D +\vphi_M$, together with $\dfadj_{\vphi_M}\paren{s\clt u} =0$
  given in \eqref{eq:dfadj-su=0}, yields
  \begin{equation*}
    \iinner{s\clt u}{\dbar\clt v_{(2)}}_{\vphilist M}
    =\iinner{\dfadj_{\vphi_M} \paren{s\clt u}}{\clt
      v_{(2)}}_{\vphilist M} = 0 \; .
  \end{equation*}
  For the inner product involving $v_{(\infty)}$,
  as $v_{(\infty)}$ is $L^2$ with respect to the potential
  $\vphi_F+\vphi_M$ but not necessarily to $\vphi_F+\phi_D+\vphi_M$,
  the cut-off function $\theta_\eps$ is introduced to facilitate the
  integration by parts.
  As $\eps \tendsto 0^+$, one has
  \begin{align*}
    \iinner{s\clt u}{\dbar v_{(\infty)}}_{\vphilist M}
    &\leftarrow \iinner{s\clt u}{ \theta_\eps \dbar v_{(\infty)}}_{\vphilist M} \\
    &=
    \iinner{s\clt u}{ \dbar \paren{\theta_\eps v_{(\infty)}}}_{\vphilist M}
    -\iinner{s\clt u}{ \frac{\eps \theta'_\eps}{\abs{\psi_D}^{1+\eps}}
      \dbar\psi_D \wedge v_{(\infty)}}_{\vphilist M} \\
    &=\cancelto{0}{\iinner{\dfadj_{\vphi_M} \paren{s\clt u}}{ \theta_\eps  v_{(\infty)}}}_{\vphilist M}
      -\eps \:\iinner{\idxup{\diff\psi_D}[\clomega] \ctrt s\clt u}{\:
      \frac{\theta'_\eps v_{(\infty)}}{\abs{\psi_D}^{1+\eps}} 
      }_{\vphilist M} \; .
  \end{align*}
  It suffices to show that the inner product on the far
  right-hand-side converges to $0$ as $\eps \tendsto 0^+$, which then
  implies that $\norm{s\clt u}_{\vphilist M} = 0$, hence $s\clt u = 0$
  and the desired $\clt u = 0$ on $X^\circ$.
  
\item \label{item:Takegoshi-argument}
  Applying to $\clt u$ the \mmark{twisted Bochner--Kodaira
    formula}{Need proof. \alert{Done.}} $\tBK_{\eps,\vphilist}$ in Lemma
  \ref{lem:BK-formulas} with the twisting function $\eta_{\eps} :=
  \abs{\psi_D}^{1-\eps}$ for $\eps > 0$, together with the help from
  \eqref{eq:limit-ctrt-u=0} and
  \eqref{eq:nable-u=0_ibddbar-vphi_F-u=0}, yields
  \begin{equation*}
    \lim_{\eps \tendsto 0^+} \eps \int_{\mathrlap{X^\circ}\;\;}
    \frac{
      \abs{\idxup{\diff\psi_D}* \ctrt \clt u}_{\vphilist}^2
    }{\abs{\psi_D}^{1+\eps}}
    =\int_{X^\circ} \idxup{\ibddbar\sm\vphi_D} \ptinner{\clt u}{\clt
      u}_{\vphilist}
  \end{equation*}
  (see Proposition \ref{prop:tBK-valid-for-u}, supplemented by the
  treatment described in Section
  \ref{sec:justification-Takegoshi-argument-singular} for the singular
  case).
  Moreover, as $\clt u \in \paren{\ker\iota_0}^\perp$, one can apply
  the argument of Takegoshi
  (\cite{Takegoshi_cohomology-nef-line-bdl}*{Prop.~3.8}, see also
  \cite{Matsumura_injectivity-lc}*{Prop.~3.13}) to claim that the
  above expression vanishes as follows.
  \mmark{Note that}{Explain pointwise BK formula in preliminaries. \alert{Done.}}
  \begin{align*}
    \dbar\dfadj_{\vphi_F+\sm\vphi_D}\clt u
    &=\paren{\dbar\dfadj_{\vphi_F+\sm\vphi_D}
        +\dfadj_{\vphi_F+\sm\vphi_D}\dbar} \clt u
    \\
    &\overset{\mathclap{\text{\eqref{eq:pointwise-BK}}}}= \quad
      -\nabla^{(1,0)}_{\vphi_F+\sm\vphi_D} \cdot \nabla^{(0,1)} \clt u
      +\ibddbar\paren{\vphi_F+\sm\vphi_D}
      \Lambda_{\clomega}\clt u
    \\
    &\overset{\mathclap{\text{\eqref{eq:nable-u=0_ibddbar-vphi_F-u=0}}}}=
      \quad
      \ibddbar\sm\vphi_D \Lambda_{\clomega} \clt u
      \quad\text{ pointwisely on } X^\circ \; ,
  \end{align*}
  where $\dfadj_{\vphi_F+\sm\vphi_D}$
  (resp.~$-\nabla^{(1,0)}_{\vphi_F+\sm\vphi_D} \cdot$) is the formal adjoint
  of $\dbar$ (resp.~$\nabla^{(0,1)}$) with respect to the inner
  product $\iinner\cdot\cdot_{\vphilist|\vphi_F+\sm\vphi_D|}$.
  The form $\ibddbar\sm\vphi_D \Lambda_{\clomega}
  \clt u$ is in $\Ltwosp$ since so is $\clt u$ and
  $\ibddbar\sm\vphi_D$ is a smooth form on $X$.
  This, together with the equality above, therefore shows that
  \begin{equation*}
    \ibddbar\sm\vphi_D \Lambda_{\clomega} \clt u
    \in \paren{\ker\dbar}_{\vphilist} = \Harm \oplus
    \paren{\im\dbar}_{\vphilist}
    \subset \Ltwosp \; .
  \end{equation*}
  Recall that, on any admissible open set $V \subset X$ such that $V
  \cap D \neq \emptyset$, one has $\res{\phi_D}_V =\res{\log\abs{\sect_D}^2}_{V}
  =\sum_{j=1}^{\sigma_V} \log\abs{z_j}^2$ and, therefore,
  \begin{equation} \label{eq:diff-psi_D-on-open-set}
    \parres{\diff\psi_D}_V =\sum_{j=1}^{\sigma_V} \frac{dz_j}{z_j} -
    \parres{\diff\sm\vphi_D}_V \; .
  \end{equation}
  Since $\vphi_F +\sm\vphi_D =\vphi_F +\phi_D -\psi_D =\vphi -\psi_D$
  and since $\frac{\clt u}{\sect_D}$ is smooth on $X^\circ$ (knowing
  its singularities along $X \setminus X^\circ$ given by Proposition
  \ref{prop:regularity-of-clt-u}), 
  \mmark{it follows that}{May justify here the use of an
    $\clomega$ which is smooth on general points on $D$ instead of a
    complete metric on $X^\circ \setminus D$. \alert{Done.}}
  \begin{equation*}
    \dfadj_{\vphi_F+\sm\vphi_D}\clt u
    =\dfadj \clt u -\idxup{\diff\psi_D} \ctrt \clt u
    =-\idxup{\diff\psi_D} \ctrt \clt u \in \Ltwosp/q-1/|\vphi_F+\alert{\sm\vphi_D}|
  \end{equation*}
  and thus $\ibddbar\sm\vphi_D \Lambda_{\clomega} \clt
  u \in \paren{\im\dbar}_{\vphilist|\vphi_F+\alert{\sm\vphi_D}|}
  \subset \Ltwosp|\vphi_F+\alert{\sm\vphi_D}|$ (i.e.~the class of
  $\ibddbar\sm\vphi_D \Lambda_{\clomega} \clt u$ in $\cohgp
  q[X]{\logKX \otimes \mtidlof{\vphi}}$ is mapped to $0$ via
  $\iota_0$), which implies that 
  \begin{equation*}
    \ibddbar\sm\vphi_D \Lambda_{\clomega} \clt u 
    \in \ker\iota_0 \oplus \paren{\im\dbar}_{\vphilist} \subset \Ltwosp \; .
  \end{equation*}
  As a result, $\clt u \in \paren{\ker\iota_0}^\perp \subset \Harm$ implies that
  \begin{equation*}
    \lim_{\eps \tendsto 0^+} \eps \int_{\mathrlap{X^\circ}\;\;}
    \frac{
      \abs{\idxup{\diff\psi_D}* \ctrt \clt u}_{\vphilist}^2
    }{\abs{\psi_D}^{1+\eps}}
    =\int_{X^\circ} \idxup{\ibddbar\sm\vphi_D}
    \ptinner{\clt u}{\clt u}_{\vphilist}
    =\iinner{
      \ibddbar\sm\vphi_D \Lambda_{\clomega} \clt u
    }{\:\clt u}_{\vphilist} = 0 \; .
  \end{equation*}

  \begin{remark} \label{rem:reason-for-smooth-omega-along-D}
    The need to invoke the above argument of Takegoshi, or, more
    precisely, to make the claim $\idxup{\diff\psi_D} \ctrt \clt u \in
    \Ltwosp/q-1/|\vphi_F +\sm\vphi_D|$ legitimate, is the reason of
    the use of the metric $\clomega$ which is smooth on $D \cap X^\circ$
    instead of a complete metric on $X^\circ \setminus D$.
    If $\clomega$ were singular along $D$, the form $\frac{\clt u}{\sect_D}$
    would have poles along $D$ even though $*_{\clomega} \frac{\clt
      u}{\sect_D}$ is holomorphic on $X$ (cf.~Proposition
    \ref{prop:regularity-of-clt-u}, which describes the possible
    singularities of $\frac{\clt u}{\sect_D}$ along $X\setminus
    X^\circ$).
    The analysis to claim $\idxup{\diff\psi_D} \ctrt \clt u$
    being $L^2$ would then be complicated, if not impossible.
  \end{remark}

\item \label{item:residue-of-final-inner-prod}
  Let $V \Subset X$ be an admissible open subset of $X$ (not
  only of $X^\circ$) with respect to $(\vphi_F,\vphi_M,\psi_D)$ in the
  holomorphic coordinate system $(z_1, \dots, z_n)$, and
  let $\rho$ be any smooth cut-off function on $X$ compactly supported
  in $V$.
  Suppose that $V \cap D \neq \emptyset$.
  Taking into account the fact that $\frac{\clt u}{\sect_D}$ is smooth
  on $X^\circ$, one can derive from Theorem
  \ref{thm:residue-fcts-and-norms} (or the computation corresponding
  to the $1$-lc-measure in
  \cite{Chan&Choi_ext-with-lcv-codim-1}*{Prop.~2.2}) that
  \begin{equation*}
    0
    =\lim_{\eps \tendsto 0^+} \eps \int_{\mathrlap{V}\;\;}
    \frac{
      \rho \abs{\idxup{\diff\psi_D}* \ctrt \clt u}_{\vphilist|\vphi_F+\alert{\phi_D}|}^2
    }{\abs{\psi_D}^{1+\eps}}
    =\pi \sum_{j=1}^{\sigma_{V}} \int_{\lcS|1|[j]}
    \rho \abs{\PRes[\lcS|1|[j]](
      \!\idxup{\diff\psi_D}* \ctrt
      \frac{\clt u}{\sect_D}
      )}_{\vphilist F}^2
    \; , 
  \end{equation*}
  where $\lcS|1|[j] := \set{z_j =0}$ for $j=1,\dots,
  \sigma_{V}$ are the irreducible components ($1$-lc centres)
  of $V \cap D$.
  (This computation is explained in the justification for the
  finiteness of this integral in Step
  \ref{item:results-from-harmonic-u-and-BK}.)
  This implies that (see Remark \ref{rem:PRes-dpsi-ctrt-u-explain})
  \begin{equation} \label{eq:smooth-dz-over-z-ctrt-u-over-sect_D}
    \parres{\idxup{dz_j}* \ctrt
        \frac{\clt u}{\sect_D}}_{\lcS|1|[j]} 
    =0
    \quad\text{ and thus }\quad
    \idxup{\frac{dz_j}{z_j}}
    \ctrt \frac{\clt u}{\sect_D}
    \in \smform/n,q-1/[X\,*]\paren{V \cap X^\circ} 
  \end{equation}
  for all $j=1,\dots,\sigma_{V}$.
  Note that $\parres{\dotsm}_{\lcS|1|[j]}$ here means only the
  restriction on the coefficients of the $(n,q-1)$-form.
  Considering the expression of $\diff\psi_D$ on $V$ in
  \eqref{eq:diff-psi_D-on-open-set} and noticing that the above holds
  true for any admissible open subset $V \Subset X$, it follows that
  \begin{equation*} \tag{%
      \begin{NoHyper}%
        \ref{eq:smooth-dz-over-z-ctrt-u-over-sect_D}$'$%
      \end{NoHyper}%
    }
    \label{eq:smooth-dpsi-ctrt-u-over-sect_D}
    \text{coef.~of~ } \idxup{\diff\psi_D} \ctrt \frac{\clt u}{\sect_D}
    \text{ ~are locally in }
    \smooth_{X\,*} \text{ on }  X^\circ \; .
  \end{equation*}

      
  Now the last inner product in Step \ref{item:by-parts-on-norm-su} can be
  considered and its limit as $\eps \tendsto 0^+$ can be evaluated.
  The absolute value of the inner product satisfies
  \begin{align*}
    &~\abs{\eps \:\iinner{\idxup{\diff\psi_D} \ctrt s\clt u}{\:
        \frac{\theta'_\eps v_{(\infty)}}{\abs{\psi_D}^{1+\eps}} 
      }_{\vphilist|\vphi_F+\alert{\phi_D}| M}} \\
    \leq
    &~\eps \int_{\mathrlap{X^\circ}\;\;}
    \frac{
      \abs{
        \inner{
          \idxup{\diff\psi_D}* \ctrt
          \frac{\clt u}{\alert{\sect_D}}
        }{
          \:v_{(\infty)} e^{-\frac 12 \alert[Plum]{\sm\vphi_D} -\frac 12 \vphi_M}}
      }_{\vphilist F}
      \:\abs{s}_{\vphi_M}
      \:\abs{\theta'_\eps}
    }{\abs{\alert{\sect_D}}_{\alert[Plum]{\sm\vphi_D}} \:\abs{\psi_D}^{1+\eps}} \; .
  \end{align*}
  Note that both $\abs s_{\vphi_M}$ and $\abs{\theta'_\eps}$ are
  bounded uniformly in $\eps$ on $X$.
  The \mmark{regularity of $\clt u$ along $X
    \setminus X^\circ$ (see Proposition \ref{prop:regularity-of-clt-u}
    for the singular case) together with
    \eqref{eq:smooth-dpsi-ctrt-u-over-sect_D} shows that the
    pointwise inner product 
    $\abs{\inner{\dotsm}{\dotsm}_{\vphilist F}}$ is integrable on
    $X$}{Need clarification. \alert{Done.}}.
  Moreover, since $\abs{\inner{\dotsm}{\dotsm}_{\vphilist F}}$ is
  locally in $\smooth_{X \,*}$ on $X^\circ$ (with poles along $X
  \setminus X^\circ$ ``independent'' of those of $\psi_D$ along $D$) while
  $\frac{1}{\abs{\sect_D}_{\sm\vphi_D} \:\abs{\psi_D}^{1+\eps}}$ 
  is integrable on $X$ for all $\eps \geq 0$,
  in view of Fubini's theorem,
  $\frac{\abs{\inner{\dotsm}{\dotsm}_{\vphilist F}}}{\abs{\sect_D}_{\sm\vphi_D}
    \:\abs{\psi_D}^{1+\eps}}$ is also integrable on $X$ (see
  Proposition \ref{prop:justification-final-inner-prod} for the
  singular case) and, therefore, the
  term on the right-hand-side above is in the order of $\BigO(\eps)$
  as $\eps \tendsto 0^+$, that is,
  \begin{equation*}
    \lim_{\eps \tendsto 0^+}
    \eps \:\iinner{\idxup{\diff\psi_D} \ctrt s\clt u}{\:
      \frac{\theta'_\eps v_{(\infty)}}{\abs{\psi_D}^{1+\eps}} 
    }_{\vphilist M}
    =0 \; .
  \end{equation*}
  Therefore, as seen from Step \ref{item:by-parts-on-norm-su}, this implies
  that $s\clt u = 0$ and thus $\clt u = 0$.
  This completes the proof.
\end{enumerate}




\subsection{Proof for the case with smooth $\vphi_F$ and $\vphi_M$}
\label{sec:pf-smooth-vphi_FM}


{
  \setDefaultMetric{\omega}

  In this section, $\vphi_F$ and $\vphi_M$ are assumed to be smooth.
  Therefore, one has $X^\circ = X$ and $\clomega = \omega$.
  For consistency, $\clt u$, $\clt v$ and $\clt v_{(2)}$ in Section
  \ref{sec:outline-of-pf} are written as $u$, $v$ and $v_{(2)}$
  respectively.

  The claims made in each step in the outline of the proof in Section
  \ref{sec:outline-of-pf} are justified in this section.

  \subsubsection{Justification for Step \ref{item:pick-u-decompose-v}}

  It suffices to provide a proof of the following lemma.

  \begin{lemma} \label{lem:regularity-of-v}
    Suppose $u \in \Harm<\omega>$ and $s \in \cohgp 0[X]{M}$ such that
    $su$ represents the zero class in $\cohgp q[X]{K_X \otimes D \otimes
      F \otimes M \otimes \mtidlof{\vphi_F+\vphi_M}}$ under the $L^2$
    Dolbeault isomorphism.
    Then, there exists $v \in \Ltwosp/q-1/M|\vphi_F+\vphi_M+\sm\vphi_D|<\omega>$ such
    that $su =\dbar v$ and $v = v_{(2)} +v_{(\infty)}$, where
    \begin{align*}
      v_{(2)} &\in \Ltwo/n,q-1/[X]{D\otimes F\otimes
                M}_{\vphi_F+\vphi_M+\alert{\phi_D},\omega}  \; , \\
      v_{(\infty)} &\in \smform/n,q-1/\paren{X; D\otimes F\otimes M
                     \otimes \mtidlof{\vphi_F+\vphi_M}} \; .
    \end{align*}
  \end{lemma}

  \begin{proof}
    \newcommand{\cubrace}[2]{
      \begingroup
      \colorlet{currcolor}{.}
      \color{Gray}
      \underbrace{\color{currcolor}#1}_{\mathclap{#2}}
      \endgroup
    }

    \newcommand{\Lbdl}{L}

    \RenewDocumentCommand{\Ltwosp}{
      D//{\bullet}             
      o                        
      D||{\vphi_F+\vphi_M +\phi_D} 
      D<>{\omega}              
    }{\IfNoValueTF{#2}{\Ltwo/n,#1/*{\Lbdl}_{#3, #4}}{\Ltwo/n,#1/[#2]{\Lbdl}_{#3, #4}}}

    The lemma is proved via examining the $L^2$ Dolbeault isomorphism,
    which is explicitly described in \cite{Matsumura_injectivity}.
    Set $\Lbdl := D \otimes F \otimes M$
    for convenience.

    Let $\cvr U :=\set{U_i}_{i\in I}$ be a finite open Stein cover of $X$ with
    a smooth partition of unity $\set{\rho^i}_{i \in I}$ subordinate to it.
    Write $U_{i_0 i_1 \dots i_\nu} := U_{i_0} \cap U_{i_1} \cap \dots
    \cap U_{i_\nu}$ for any indices $i_0, i_1, \dots , i_\nu \in I$ as
    usual.
    Also let $\delta$ to be the coboundary operator on the \v Cech
    cochains. 
    Following the construction of the $L^2$ Dolbeault isomorphism (see
    \cite{Matsumura_injectivity}*{Prop.~5.5} and
    \cite{Matsumura_injectivity-lc}*{Prop,~2.8}), with the given
    $(n,q)$-form $su$, solve the equations
    \begin{alignat*}{2}
      \dbar\set{\beta_{i_0}}
      &=\set{\res{su}_{U_{i_0}}}
      &&\text{for } \beta_{i_0} \in
      \Ltwosp/q-1/[U_{i_0}]
      \quad\text{and}
      \\
      \dbar\set{\beta_{i_0 \dots i_\nu}}
      &=\underbrace{\delta\set{\beta_{i_0 \dots i_{\nu-1}}}}_{=:
        \:\set{\alpha_{i_0 \dots i_\nu}}}
      &\quad&\text{for } \beta_{i_0 \dots i_\nu} \in
      \Ltwosp/q-\nu-1/[U_{i_0 \dots i_\nu}]
    \end{alignat*}
    for $\nu = 0, \dots, q-1$ via $L^2$ method on relatively compact
    Stein subsets (see, for example, \cite{DHP}*{\S 4} for the
    procedures for smoothing out the metric before solving for the
    $\dbar$-equation on a relatively compact Stein subset).
    Note that $\delta\set{\beta_{i_0 \dots i_{q-1}}} =: \set{\alpha_{i_0
        \dots i_q}}$ has holomorphic components
    and is representing the same cohomology class as $su$ under the
    $L^2$ Dolbeault isomorphism.
    Using the fact that $\sum_{i\in I} \dbar\rho^i \equiv 0$ on $X$, the
    $(n,q)$-form $su$ can then be expressed as (under the Einstein
    summation convention)
    \begin{align*}
      su = \rho^{i_0} \dbar\beta_{i_0}
      &=\dbar\paren{\rho^{i_0} \beta_{i_0}}
        \cubrace{-\dbar\rho^{i_0} \wedge
        \beta_{i_0}}{=\:+\dbar\rho^{i_0} \wedge \rho^{i_1}
        \paren{\beta_{i_1} -\beta_{i_0}}} 
        =\dbar\paren{\rho^{i_0} \beta_{i_0}}
        +\dbar\rho^{i_0} \wedge \rho^{i_1} \dbar\beta_{i_0 i_1} \\
      &=\dbar\paren{\rho^{i_0} \beta_{i_0}}
        +\dbar\rho^{i_0} \wedge \dbar\paren{\rho^{i_1}\beta_{i_0 i_1}}
        \cubrace{-\dbar\rho^{i_0} \wedge \dbar\rho^{i_1} \wedge
        \beta_{i_0 i_1}}{=\: +\dbar\rho^{i_1} \wedge
        \dbar\rho^{i_0} \wedge \rho^{i_2} \paren{\beta_{i_0 i_1}
        -\beta_{i_0 i_2} +\beta_{i_1 i_2}} } \\
      &=\dbar\paren{\rho^{i_0} \beta_{i_0}
        -\dbar\rho^{i_0} \wedge
        \rho^{i_1}\beta_{i_0 i_1}}
        +\dbar\rho^{i_1} \wedge
        \dbar\rho^{i_0} \wedge \rho^{i_2} \dbar\beta_{i_0 i_1 i_2} \\
      &=
        \begin{aligned}[t]
          &\dbar\paren{ \rho^{i_0} \beta_{i_0}
            -\dbar\rho^{i_0} \wedge
            \rho^{i_1}\beta_{i_0 i_1}
            +\dbar\rho^{i_1} \wedge
            \dbar\rho^{i_0} \wedge  \rho^{i_2} \beta_{i_0
              i_1 i_2}} \\
          &\cubrace{-\dbar\rho^{i_2} \wedge \dbar\rho^{i_1}
            \wedge \dbar\rho^{i_0} \wedge \beta_{i_0 i_1 i_2}}{=
            \:+\dbar\rho^{i_2} \wedge \dbar\rho^{i_1}
            \wedge \dbar\rho^{i_0} \wedge \rho^{i_3} \:\alpha_{i_0 i_1 i_2 i_3}}
        \end{aligned} \\
      &=
        \begin{aligned}[t]
          &
          \begin{multlined}[t]
            \dbar\left(
              \alert{\rho^{i_0} \beta_{i_0}
              -\dbar\rho^{i_0} \wedge
              \rho^{i_1}\beta_{i_0 i_1}
              +\dbar\rho^{i_1} \wedge
              \dbar\rho^{i_0} \wedge \rho^{i_2} \beta_{i_0 i_1 i_2}}
            \right. \\
            \left. 
              \alert{-\dots
              +(-1)^{q-1} \:\dbar\rho^{i_{q-2}} \wedge \dots \wedge
              \dbar\rho^{i_0} 
              \cdot \rho^{i_{q-1}} \beta_{i_0\dots i_{q-1}}}
            \right)
          \end{multlined}
          \\
          &+\dbar\rho^{i_{q-1}} \wedge \dbar\rho^{i_{q-2}}
          \wedge \dots \wedge \dbar\rho^{i_0} \cdot \rho^{i_q}
          \alpha_{i_0 \dots i_{q}}
        \end{aligned}
      \\
      &=: \dbar \alert{v_{(2)}} +\dbar\rho^{i_{q-1}} \wedge \dbar\rho^{i_{q-2}}
        \wedge \dots \wedge \dbar\rho^{i_0} \cdot \rho^{i_q}
        \alpha_{i_0 \dots i_{q}} \; .
    \end{align*}
    Notice that, since all $\beta_{i_0,\dots,i_\nu}$'s are $L^2$ with respect
    to the weight $e^{-\vphi_F-\vphi_M -\phi_D}$, the $(n,q-1)$-form
    $v_{(2)}$
    can be chosen to be in $\Ltwosp/q-1/[X]$.

    As $su$ is representing the zero class in $\cohgp q[X]{K_X\otimes
      \Lbdl\otimes \mtidlof{\vphi_F+\vphi_M}}$, there exists
    a $(q-1)$-cochain $\set{\gamma_{i_0 \dots i_{q-1}}}$ such that
    $\gamma_{i_0 \dots i_{q-1}} \in K_X \otimes \Lbdl \otimes
    \mtidlof{\vphi_F+\vphi_M}\paren{U_{i_0\dots i_{q-1}}}$ (holomorphic) for
    all indices $i_0, \dots, i_{q-1} \in I$ and $\set{\alpha_{i_0 \dots
        i_q}} = \delta\set{\gamma_{i_0 \dots i_{q-1}}}$.
    Since
    \begin{align*}
      \dbar\rho^{i_{q-1}}
      \wedge \dots \wedge \dbar\rho^{i_0} \cdot \rho^{i_q}
      \alpha_{i_0 \dots i_{q}}
      &=\dbar\rho^{i_{q-1}}
        \wedge \dots \wedge \dbar\rho^{i_0} \cdot \rho^{i_q}
        \sum_{\nu=0}^q (-1)^\nu \gamma_{i_0 \dots \widehat{i_\nu} \dots
        i_{q}} \\
      &=(-1)^q \:\dbar\rho^{i_{q-1}}
        \wedge \dots \wedge \dbar\rho^{i_0} \cdot \gamma_{i_0 \dots
        i_{q-1}} \\
      &=\dbar\paren{(-1)^q \:\dbar\rho^{i_{q-2}}
        \wedge \dots \wedge \dbar\rho^{i_0}
        \cdot \rho^{i_{q-1}}\gamma_{i_0 \dots i_{q-1}}}
        =: \dbar v_{(\infty)}\; ,
    \end{align*}
    and since all $\rho^{i_\nu}$'s and $\gamma_{i_0\dots i_{q-1}}$'s are
    smooth functions, one sees that $v_{(\infty)}$ can be chosen to be
    smooth.
    One then obtains $v := v_{(2)} +v_{(\infty)} \in
    \Ltwosp/q-1/[X]|\vphi_F+\vphi_M+\sm\vphi_D|$ such that $su = \dbar v$
    with the acclaimed properties.
  \end{proof}

  \begin{remark} \label{rem:regularity-v-for-sing-F-M}
    \setDefaultMetric{\clomega}
    The proof of Lemma \ref{lem:regularity-of-v} is still valid when
    $\vphi_F$ and $\vphi_M$ have neat analytic singularities as
    described in Section \ref{sec:setup} and when $\clt u \in
    \Harm<\alert{\clomega}>$ instead of $u \in \Harm<\omega>$ is
    considered in the statement.
    This can be seen by noticing that the local $\dbar$-equations on
    $U_{i_0 \dots i_\nu}$ (for instance, $\dbar\set{\beta_{i_0}}
    =\set{\res{su}_{U_{i_0}}}$ for $\beta_{i_0} \in
    \Ltwo/n,q-1/[U_{i_0}]{D\otimes F\otimes
      M}_{\vphilist|\vphi_F+\vphi_M+\phi_D|}$) can be solved first on
    $U_{i_0 \dots i_\nu} \cap \paren{X^\circ \setminus D}$ with $L^2$
    estimate even though $\clomega$ is not complete there (see
    \cite{Demailly}*{Ch.~VIII, Thm.~(6.1)} or
    \cite{Matsumura_injectivity}*{Lemma 5.4}) and the solution can
    then be extended to $U_{i_0 \dots i_\nu}$ via the $L^2$ Riemann
    extension theorem (see \cite{Demailly_complete-Kahler}*{Lemma
      6.9}).
    The form $v_{(2)}$ obtained in the conclusion should then be
    replaced by
    \begin{equation*}
      \clt v_{(2)} \in \Ltwo/n,q-1/[X^\circ]{D\otimes F\otimes
        M}_{\vphi_F+\vphi_M+\phi_D,\alert{\clomega}}  \; .
    \end{equation*}
  \end{remark}
  

  \subsubsection{Justification for Step \ref{item:results-from-harmonic-u-and-BK}}
  \label{sec:justify-step-II-smooth}

  To justify the arguments in Step
  \ref{item:results-from-harmonic-u-and-BK}, one has to show that
  \begin{enumerate}
  \item \eqref{eq:limit-ctrt-u=0} holds true,
    
  \item $u$ satisfies the Bochner--Kodaira formula $\BK_{\vphilist}$,
    which yields \eqref{eq:nable-u=0_ibddbar-vphi_F-u=0}, and
    
  \item $su$ satisfies $\BK_{\vphilist M}$, which yields
    \eqref{eq:dfadj-su=0}. 
  \end{enumerate}

  It follows from the refined hard Lefschetz theorem in
  \cite{Matsumura_injectivity-lc}*{Thm.~3.3} (see also Theorem
  \ref{thm:refined-hard-Lefschetz}) that $*_\omega u$ is holomorphic
  on $X$, which implies that $u =(-1)^{n-q} *_\omega *_\omega u$ is
  smooth on $X$.
  Recall that $\sect_D$ is a canonical section of $D$ such that
  $\phi_D =\log\abs{\sect_D}^2$.
  Since $\abs u^2 = \abs{*_\omega u}_\omega^2 \dvol_{X,\omega}$ is
  $L^1$ with respect to the weight $e^{-\vphi_F-\phi_D}$, it follows
  that $\frac{*_\omega u}{\sect_D}$, hence $\frac{u}{\sect_D}$, is
  smooth on $X$.
  This piece of information is sufficient to give the following
  proposition.

  \begin{prop} \label{prop:residue-of-dpsi-ctrt-u}
    With the knowledge that $\frac{u}{\sect_D}$ being smooth on $X$,
    one has
    \begin{equation*}
      \lim_{\eps \tendsto 0^+} \eps
      \int_{\mathrlap{X}\;\;}
      \frac{
        \abs{\dbar\psi_D \otimes u}_{\vphilist}^2
      }{\abs{\psi_D}^{1+\eps}}
      =\pi \sum_{i \in I_D} \int_{D_i} \abs{
        \PRes[D_i](\dbar\psi_D \otimes \frac{u}{\sect_D})
      }_{\vphilist F}^2 < +\infty \; ,
    \end{equation*}
    where $D =\sum_{i \in I_D} D_i$ is the decomposition of $D$ into
    irreducible components in $X$ and $\PRes[D_i]$ is the
    \textfr{Poincaré} residue map corresponding to the restriction
    from $X$ to $D_i$ (see \cite{Kollar_Sing-of-MMP}*{Def.~4.1 and
      para.~4.18}; see also Section \ref{sec:residue-functions}).
    In particular, the equation \eqref{eq:limit-ctrt-u=0} holds true, that is,
    \begin{equation*}
      \lim_{\eps \tendsto 0^+} \alert{\eps^2}
      \int_{\mathrlap{X}\;\;}
      \frac{
        \abs{\dbar\psi_D \otimes u}_{\vphilist}^2
      }{\abs{\psi_D}^{1+\eps}} = 0 \; .
    \end{equation*}
  \end{prop}
  
  \begin{proof}
    Let $\set{V_\gamma}_{\gamma \in \Gamma}$ be a \emph{finite} open
    cover of $X$ such that each $V_\gamma$ is an admissible open set
    with respect to $\psi_D$ in the holomorphic coordinate system
    $(z^\gamma_1, \dots, z^\gamma_n)$, and let
    $\set{\rho_\gamma}_{\gamma\in \Gamma}$ be a partition of unity
    subordinate to this cover.
    On any $V_\gamma$ such that $V_\gamma \cap D \neq \emptyset$, one has
    \begin{equation*}
      \parres{\diff\psi_D}_{V_\gamma}
      =\sum_{j=1}^{\sigma_{V_\gamma}} \frac{dz_j^\gamma}{z_j^\gamma}
      -\parres{\diff\sm\vphi_D}_{V_\gamma} \; ,
    \end{equation*}
    and, therefore, on $V_\gamma$,
    \begin{equation*}
      \dbar\psi_D \otimes u
      =\sum_{j=1}^{\sigma_{V_\gamma}} \paren{
        d\conj{z_j^\gamma}
        -\frac{\conj{z_j^\gamma} \:\dbar\sm\vphi_D}{\sigma_{V_\gamma}}
      } \otimes \frac u{\conj{z_j^\gamma}}
      =\sum_{j=1}^{\sigma_{V_\gamma}}
      \underbrace{
        \paren{
          d\conj{z_j^\gamma}
          -\frac{\conj{z_j^\gamma} \:\dbar\sm\vphi_D}{\sigma_{V_\gamma}}
        } \otimes \frac u{\sect_D} \:\frac{z_j^\gamma}{\conj{z_j^\gamma}}
      }_{\textstyle =: \: dz_j^\gamma \wedge g_j^\gamma \;\;\footnotemark}
      \:\prod_{\substack{k=1 \\ k\neq j}}^{\sigma_{V_\gamma}}
      z_k^\gamma \; .
    \end{equation*}\footnotetext{\label{fn:definition-of-g_j}%
      More precisely, $g_j^\gamma$ is chosen to be $\fdiff{z_j^\gamma}
      \ctrt \paren{d\conj{z_j^\gamma} -\frac{\conj{z_j^\gamma}
          \:\dbar\sm\vphi_D}{\sigma_{V_\gamma}}} \otimes \frac
      u{\sect_D} \:\frac{z_j^\gamma}{\conj{z_j^\gamma}} \;$, so
      $g_j^\gamma$ contains no terms of the holomorphic differential
      form $dz_j^\gamma$.  
    }%
    Compare this expansion with the one in \eqref{eq:decomposition-f}.
    As $\frac u{\sect_D}$ is smooth on $X$, it follows that $g_j^\gamma$
    is in $\smooth_{X\,*}$ on $V_\gamma$ for $j=1,\dots, \sigma_{V_\gamma}$ and
    thus satisfies the conditions in \eqref{eq:regularity-on-g_p}.
    Note that, for every $\lcS|1,\gamma|[j] :=\set{z_j^\gamma =
      0}$ where $j=1,\dots,\sigma_{V_\gamma}$ (using the notation in Section
    \ref{sec:residue-functions}), there exists a unique $i \in
    I_D$ such that $\lcS|1,\gamma|[j] = V_\gamma \cap D_i$ and 
    \begin{equation*}
      \res{g_j^\gamma}_{\lcS|1,\gamma|[j]}
      =\PRes[\lcS|1,\gamma|[j]](\dbar\psi_D \otimes \frac u{\sect_D})
      =\parres{\PRes[D_i](\dbar\psi_D \otimes \frac
        u{\sect_D})}_{V_\gamma \cap D_i} \; .
    \end{equation*}
    It, therefore, follows from Theorem
    \ref{thm:residue-fcts-and-norms} (or the computation corresponding
    to the $1$-lc-measure in
    \cite{Chan&Choi_ext-with-lcv-codim-1}*{Prop.~2.2.1}) that
    \begin{equation*}
      \lim_{\eps \tendsto 0^+} \eps \int_{\mathrlap{V_\gamma}\;\;}
      \frac{
        \rho_\gamma \abs{\dbar\psi_D \otimes u}_{\vphilist}^2
      }{\abs{\psi_D}^{1+\eps}}
      =\pi \sum_{j=1}^{\sigma_{V_\gamma}} \int_{\lcS|1,\gamma|[j]}
      \rho_\gamma \abs{g_j^\gamma}_{\vphilist F}^2 <+\infty \; .
    \end{equation*}
    The claim thus follows after summing up over $\gamma \in \Gamma$,
    where $\Gamma$ is just a finite set (thus finiteness of the
    integral is guaranteed).
    \qedhere

  \end{proof}

  \begin{remark} \label{rem:PRes-dpsi-ctrt-u-explain}
    The computation and results in Proposition
    \ref{prop:residue-of-dpsi-ctrt-u} 
    still hold true when $\dbar\psi_D \otimes u$ is replaced by
    $\idxup{\diff\psi_D} \ctrt u$ or $\dbar\psi_D \wedge u$.
    For fear of being confused by the notation and also for the use in
    Step \ref{item:residue-of-final-inner-prod} of Section
    \ref{sec:outline-of-pf},
    notice that, if $
    \res{g_j^\gamma}_{\lcS|1,\gamma|[j]}
    =\PRes[\lcS|1,\gamma|[j]](\!\idxup{\diff\psi_D} \ctrt \frac 
    u{\sect_D}) = 0$ under the notation in the proof above, one
    actually has $\parres{dz_j^\gamma \wedge
      g_j^\gamma}_{\lcS|1,\gamma|[j]} =\parres{\!\idxup{dz_j^\gamma}
      \ctrt \frac{u}{\sect_D}}_{\lcS|1,\gamma|[j]} = 0$,
    where $\parres{\dotsm}_{\lcS|1,\gamma|[j]}$ here is the restriction
    only on the coefficients of the $(n,q-1)$-form.
  \end{remark}

  It can now be shown that $u$ satisfies the Bochner--Kodaira formula
  $\BK_{\vphilist}$ in Lemma \ref{lem:BK-formulas} and therefore
  \eqref{eq:nable-u=0_ibddbar-vphi_F-u=0}.
  For that purpose, let $\Dom \dbar \subset \Ltwosp$ be the domain of
  $\dbar$ on $\Ltwosp$ and $\Dom \dbadj , \Dom \dfadj \subset \Ltwosp$
  be the domains of respectively the Hilbert space adjoint and the
  formal adjoint of $\dbar \colon \Ltwosp/q+1/ \birat \Ltwosp$ (see
  Lemma \ref{lem:formal-adjoint-densely-defined}).
  Note that one has $\Dom\dbadj \subset \Dom\dfadj$ and $\dbadj =
  \dfadj$ on $\Dom\dbadj$ (see Section \ref{sec:BK-formulas}).
  Let again $\theta_\eps :=\theta \circ \frac{1}{\abs{\psi_D}^{\eps}}$
  be the smooth cut-off function described in Section
  \ref{sec:outline-of-pf}. 

  \begin{prop} \label{prop:nabla-u=0_curv-u=0}
    Given $u \in \Harm$, which is smooth on $X$ and satisfies
    \eqref{eq:limit-ctrt-u=0} according to Proposition
    \ref{prop:residue-of-dpsi-ctrt-u}, it follows that $u$ satisfies
    the Bochner--Kodaira formula $\BK_{\vphilist}$ in Lemma
    \ref{lem:BK-formulas} and, consequently,
    \eqref{eq:nable-u=0_ibddbar-vphi_F-u=0}, i.e.~
    \begin{equation*}
      \nabla^{(0,1)} u = 0 \quad\text{ and }\quad
      \idxup{\ibddbar\vphi_F}\ptinner{u}{u}_{\omega} = 0 
      \quad\text{ on } X \; .
    \end{equation*}
  \end{prop}

  \begin{proof}
    Since $u$ is smooth on $X$, it follows that $\theta_{\eps} u$
    satisfies $\BK_{\vphilist}$ by Lemma \ref{lem:BK-formulas}.
    Moreover, $u \in \ker\dbadj \subset \Dom\dbadj$ implies that
    $\dfadj u = 0$.
    Noting the pointwise identities $\dbar\paren{\theta_\eps u}
    =\dbar\theta_\eps \wedge u$, \;
    $\dfadj\paren{\theta_\eps u} =-\idxup{\diff\theta_\eps} \ctrt u$, \;
    $\nabla^{(0,1)}\paren{\theta_\eps u} =\theta_\eps \nabla^{(0,1)}u
    +\dbar\theta_\eps \otimes u$ and 
    \begin{equation*}
      \abs{\dbar\theta_\eps \wedge u}_{\vphilist}^2
      +\abs{\idxup{\diff\theta_\eps} \ctrt u}_{\vphilist}^2
      =\abs{\dbar\theta_\eps \otimes u}_{\vphilist}^2
      \quad\text{ on } X\;\footnotemark
    \end{equation*}%
    \footnotetext{The identity can be obtained by applying
      $\inner{u}{\cdot}_{\vphilist}$ to both sides of the identity (see,
      for example, \cite{Federer}*{1.5.3})
      \begin{equation*}
        u = \frac{
          \idxup{\diff\theta_\eps} \ctrt \paren{\dbar\theta_\eps \wedge u}
        }{\abs{\dbar\theta_\eps}_\omega^2}
        +\frac{
          \dbar\theta_\eps \wedge \paren{\idxup{\diff\theta_\eps} \ctrt  u}
        }{\abs{\dbar\theta_\eps}_\omega^2} \; .
      \end{equation*}
      Note that $\abs{\dbar\theta_\eps \otimes u}_{\vphilist}^2
      =\abs{\dbar\theta_\eps}_\omega^2 \: \abs{u}_{\vphilist}^2$ and
      $\dbar\theta_\eps \wedge \cdot$ and $\idxup{\diff\theta_\eps}
      \ctrt \cdot$ are adjoints of each other.
    }%
    and $\ibddbar\phi_D = 0$
    on $X \setminus D$ as well, $\BK_{\vphilist}$ (with $\theta_\eps
    u$ in place of $\zeta$) yields
    \begin{equation*}
      0 = \norm{\theta_\eps \nabla^{(0,1)} u}_{\vphilist}^2
      +2\Re\iinner{\theta_\eps \nabla^{(0,1)} u}{
        \frac{\eps \theta'_\eps}{\abs{\psi_D}^{1+\eps}}
        \dbar\psi_D \otimes u
      }_{\vphilist}
      +\int_X \idxup{\ibddbar\vphi_F}
      \ptinner{\theta_\eps u}{\theta_\eps u}_{\vphilist} \; .
    \end{equation*}
    Noting that $\theta'_\eps$ is bounded uniformly in $\eps$ on $X$,
    a use of the Cauchy--Schwarz inequality followed by the AM-GM
    inequality for any fixed constant $\alpha \in (0,1)$ on the inner
    product above yields
    \begin{equation*}
      (1-\alpha) \norm{\theta_\eps \nabla^{(0,1)} u}_{\vphilist}^2
      +\int_X \idxup{\ibddbar\vphi_F}
      \ptinner{\theta_\eps u}{\theta_\eps u}_{\vphilist}
      \lesssim
      \frac{\eps^2}{\alpha} \int_X
      \frac{
        \abs{\dbar\psi_D \otimes u}_{\vphilist}^2
      }{\abs{\psi_D}^{2+2\eps}} \; ,
    \end{equation*}
    where the constant involved in $\lesssim$ is independent of
    $\eps$.
    As $\frac{1}{\abs{\psi_D}^{2+2\eps}} \leq
    \frac{1}{\abs{\psi_D}^{1+\eps}}$, Proposition
    \ref{prop:residue-of-dpsi-ctrt-u} guarantees that the
    right-hand-side above has its limit equal $0$ as $\eps \tendsto 0^+$.
    Using the assumption that $\ibddbar\vphi_F \geq 0$, one can apply
    Fatou's lemma on the left-hand-side of the inequality above and
    obtain
    \begin{equation*}
      0 \leq (1-\alpha) \norm{\nabla^{(0,1)} u}_{\vphilist}^2
      +\int_X \idxup{\ibddbar\vphi_F}
      \ptinner{u}{u}_{\vphilist}
      \lesssim 0 \; .
    \end{equation*}
    The desired equalities thus follow.
  \end{proof}


  A similar argument can show that the Bochner--Kodaira formula
  $\BK_{\vphilist M}$ is valid for $su$ and the desired claims in
  \eqref{eq:dfadj-su=0} can then be proved.

  \begin{cor} \label{cor:dfadj_M-su=0}
    Given $u \in \Harm$, which satisfies the hypotheses and claims in
    Propositions \ref{prop:residue-of-dpsi-ctrt-u} and
    \ref{prop:nabla-u=0_curv-u=0}, it follows that $su$ satisfies the
    Bochner--Kodaira formula $\BK_{\vphilist M}$ and, consequently,
    claims in \eqref{eq:dfadj-su=0} hold, that is,
    \begin{equation*}
      \dfadj_{\vphi_M}\paren{su} = 0 \quad\text{ and }\quad
      su \in \Dom\dbadj_{\vphi_M} \; ,
    \end{equation*}
    where $\dbadj_{\vphi_M}$ and $\dfadj_{\vphi_M}$ are respectively
    the Hilbert space and formal adjoints of $\dbar$ with respect to
    $\iinner{\cdot}{\cdot}_{\vphilist M}$.
  \end{cor}

  \begin{proof}
    Since $su$ is smooth on $X$, $\BK_{\vphilist M}$ is valid for
    $\theta_\eps su$.
    Noting that $\dbar\paren{su} =s \dbar u =0$, \;
    $\nabla^{(0,1)}\paren{su} =s \nabla^{(0,1)} u = 0$ and
    $\dfadj_{\vphi_M}\paren{\theta_\eps su} = \theta_\eps
    \dfadj_{\vphi_M}\paren{su} -\idxup{\diff\theta_\eps} \ctrt su$, together with
    the assumption $\ibddbar\vphi_M \leq C \ibddbar\vphi_F$ for some
    constant $C > 0$ (see Section \ref{sec:setup}), the argument as in
    Proposition \ref{prop:nabla-u=0_curv-u=0} turns
    $\BK_{\vphilist M}$ (with $\theta_\eps s u$ in place of $\zeta$) into
    \begin{align*}
      &~\norm{\theta_\eps \dfadj_{\vphi_M} \paren{su}}_{\vphilist M}^2
      -2\Re\iinner{
        \theta_\eps \dfadj_{\vphi_M} \paren{su}
      }{
        \frac{\eps \theta'_\eps}{\abs{\psi_D}^{1+\eps}}
        \idxup{\diff\psi_D} \ctrt su
      }_{\vphilist M} \\
      =
      &~\int_X \idxup{\ibddbar\paren{\vphi_F +\vphi_M}}
      \ptinner{\theta_\eps su}{\theta_\eps su}_{\vphilist M} \\
        \leq 
      &~\paren{1+C} \int_X \abs s_{\vphi_M}^2 \theta_\eps^2
      \idxup{\ibddbar\vphi_F} \ptinner u u_{\vphilist}
      =0 \; .
    \end{align*}
    The use of the Cauchy--Schwarz and AM-GM inequalities for some
    constant $\alpha \in (0,1)$ then yields
    \begin{align*}
      \paren{1-\alpha} \norm{\theta_\eps \dfadj_{\vphi_M}
        \paren{su}}_{\vphilist M}^2
      &\lesssim \frac{\eps^2}{\alpha} \int_X
      \frac{\abs{\idxup{\diff\psi_D} \ctrt su}_{\vphilist M}^2}{
        \abs{\psi_D}^{2+2\eps}
      } \\
      &\leq \sup_X\abs s_{\vphi_M}^2 \:
       \frac{\eps^2}{\alpha} \int_X
      \frac{\abs{\idxup{\diff\psi_D} \ctrt u}_{\vphilist}^2}{
        \abs{\psi_D}^{2+2\eps}
      } \; ,
    \end{align*}
    where the constant involved in $\lesssim$ is from the estimate of
    $\theta'_\eps$ which is independent of $\eps$.
    Applying Fatou's lemma on the left-hand-side and Proposition
    \ref{prop:residue-of-dpsi-ctrt-u} (see also Remark
    \ref{rem:PRes-dpsi-ctrt-u-explain}) on the right-hand-side while
    taking the limit $\eps \tendsto 0^+$ gives
    \begin{equation*}
      \dfadj_{\vphi_M} \paren{su} = 0
    \end{equation*}
    as desired.

    To see that $su\in \Dom\dbadj_{\vphi_M}$ (see, for example,
    \cite{Demailly}*{Ch.~VIII, \S 1} for the definition of the domain
    of an Hilbert space adjoint), notice that, for every
    $\zeta \in \Dom\dbar \subset \Ltwosp/q-1/M|\vphi+\vphi_M|$, by considering a
    partition of unity and the convolution with smoothing kernels on
    local coordinate charts, one obtains a sequence
    $\seq{\zeta_{\eps,\nu}}_{\nu\in\Nnum}$ of smooth $D\otimes
    F\otimes M$-valued $(n,q-1)$-forms \emph{compactly supported in $X
      \setminus D$} such that $\zeta_{\eps,\nu} \tendsto \theta_\eps
    \zeta$ in the graph norm $\paren{\norm\cdot_{\vphilist M}^2
      +\norm{\dbar \:\cdot}_{\vphilist M}^2}^{\frac 12}$ of $\dbar$
    according to the lemma of Friedrichs (see \cite{Friedrichs} or
    \cite{Demailly}*{Ch.~VIII, Thm.~(3.2)}).
    Therefore, one sees that
    \begin{align*}
      \iinner{su}{\dbar\zeta}_{\vphilist M} 
      \xleftarrow{\eps \tendsto 0^+}
      &~\iinner{su}{\theta_\eps \dbar\zeta}_{\vphilist M} \\
      =&~\iinner{su}{\dbar\paren{\theta_{\eps}\zeta}}_{\vphilist M}
         -\iinner{su}{\dbar\theta_\eps \wedge \zeta}_{\vphilist M} \\
      \xleftarrow{\nu \tendsto \infty}
      &~\iinner{su}{\dbar\zeta_{\eps,\nu}}_{\vphilist M}
        -\iinner{su}{\frac{\eps \theta'_\eps}{\abs{\psi_D}^{1+\eps}}
        \dbar\psi_D \wedge \zeta}_{\vphilist M} \\
      =&~\iinner{\dfadj_{\vphi_M} \paren{su}}{\zeta_{\eps,\nu}}_{\vphilist M}
         -\iinner{\frac{\eps \theta'_\eps}{\abs{\psi_D}^{1+\eps}}
         \idxup{\diff\psi_D} \ctrt su}{\zeta}_{\vphilist M} \; .
    \end{align*}
    As $\abs{\idxup{\diff\psi_D} \ctrt su}_{\vphilist M}^2 
    \leq \sup_X \abs s_{\vphi_M}^2 \: 
    \abs{\idxup{\diff\psi_D} \ctrt u}_{\vphilist}^2 $, Proposition
    \ref{prop:residue-of-dpsi-ctrt-u} (and Remark
    \ref{rem:PRes-dpsi-ctrt-u-explain}) guarantees that the inner
    product on the far right-hand-side converges to $0$ when $\eps
    \tendsto 0^+$.
    Since $\iinner{\dfadj_{\vphi_M}
      \paren{su}}{\zeta_{\eps,\nu}}_{\vphilist M} = 0$ (or one has
    $\lim_{\eps \tendsto 0^+} \lim_{\nu \tendsto \infty} 
    \iinner{\dfadj_{\vphi_M} \paren{su}}{\zeta_{\eps,\nu}}_{\vphilist
      M}
    =\lim_{\eps \tendsto 0^+} \iinner{\dfadj_{\vphi_M}
      \paren{su}}{\theta_\eps \zeta}_{\vphilist M} 
    =\iinner{\dfadj_{\vphi_M} \paren{su}}{\zeta}_{\vphilist M}$ if
    $\dfadj_{\vphi_M}\paren{su}$ were not $0$), it can be seen that $\Dom\dbar
    \ni \zeta \mapsto \iinner{su}{\dbar\zeta}_{\vphilist M}$ is a
    bounded linear functional, so $su \in \Dom\dbadj_{\vphi_M}$.
  \end{proof}

  \begin{remark}
    The argument for the claim $su \in \Dom\dbadj_{\vphi_M}$ indeed
    shows that
    \begin{align*}
      \Dom \dbadj \cap \smform/n,q/\paren{X;D\otimes F}
      &=\Dom \dfadj \cap \smform/n,q/\paren{X;D\otimes F}
        \quad\text{ and}
      \\
      \Dom \dbadj_{\vphi_M} \cap \smform/n,q/\paren{X;D\otimes
      F\otimes M}
      &=\Dom \dfadj_{\vphi_M} \cap \smform/n,q/\paren{X;D\otimes
      F\otimes M} \; .
    \end{align*}
  \end{remark}

  \subsubsection{Justification for Step \ref{item:Takegoshi-argument}}

  There is no need of extra clarification for Step
  \ref{item:by-parts-on-norm-su}.
  To justify Step \ref{item:Takegoshi-argument} in the current
  situation (where both $\vphi_F$ and $\vphi_M$ are smooth and
  $\frac{u}{\sect_D}$ is smooth on the whole of $X$), it
  suffices to show the following.

  \begin{prop} \label{prop:tBK-valid-for-u}
    Given $u \in \Harm$, which satisfies the hypotheses and claims in
    Propositions \ref{prop:residue-of-dpsi-ctrt-u} and
    \ref{prop:nabla-u=0_curv-u=0}, it follows that $u$ satisfies the
    twisted Bochner--Kodaira formula $\tBK_{\eps,\vphilist}$ with
    $\eps > 0$ and, consequently,
    \begin{equation*}
      \lim_{\eps \tendsto 0^+} \eps \int_{\mathrlap{X}\;\;}
      \frac{
        \abs{\idxup{\diff\psi_D}* \ctrt u}_{\vphilist}^2
      }{\abs{\psi_D}^{1+\eps}}
      =\int_{X} \idxup{\ibddbar\sm\vphi_D}
      \ptinner{u}{u}_{\vphilist} \; .
    \end{equation*}
  \end{prop}

  \begin{proof}
    The proof is similar to the one of Proposition
    \ref{prop:nabla-u=0_curv-u=0}.
    Note that $\theta_{\alert{\eps'}} u$ satisfies $\tBK_{\eps,\vphilist}$
    for any $\eps, \eps' > 0$.
    Taking into account of the vanishing results $\dbar u =0$, $\dfadj
    u =0$, $\nabla^{(0,1)}u =0$ and $\idxup{\ibddbar\vphi_F}\ptinner u
    u_{\omega} = 0$ on $X$, together with $\ibddbar\phi_D = 0$ on
    $X\setminus D$, the argument as in Proposition
    \ref{prop:nabla-u=0_curv-u=0} turns $\tBK_{\eps,\vphilist}$
    into
    \begin{equation*}
      0
      =
      \begin{aligned}[t]
        &-\int_X \frac{1-\eps}{\abs{\psi_D}} \idxup{\ibddbar\sm\vphi_D}
        \ptinner{\theta_{\eps'} u}{\theta_{\eps'} u}_{\vphilist}
        \:\abs{\psi_D}^{1-\eps} \\
        &+2(1-\eps) \Re \int_X
        \inner{\idxup{\diff\theta_{\eps'}} \ctrt u}{
          \frac{\idxup{\diff\psi_D} \ctrt \theta_{\eps'} u}{\abs{\psi_D}}
        }_{\vphilist} \abs{\psi_D}^{1-\eps} \\
        &+\eps \int_X \frac{1-\eps}{\abs{\psi_D}^2}
        \abs{\idxup{\diff\psi_D} \ctrt \theta_{\eps'} u}_{\vphilist}^2
        \:\abs{\psi_D}^{1-\eps} \; .
      \end{aligned}
    \end{equation*}
    After dividing the factor $1-\eps$, expanding $\diff\theta_{\eps'}$
    and rearranging terms, it becomes
    \begin{equation*}
      \int_X \frac{\theta_{\eps'}^2}{\abs{\psi_D}^\eps}
      \idxup{\ibddbar\sm\vphi_D}
      \ptinner{u}{u}_{\vphilist}
      =
      \int_X \paren{
        \frac{\eps \theta_{\eps'}^2}{\abs{\psi_D}^{1+\eps}}
        +\frac{2\eps' \:\theta'_{\eps'} \theta_{\eps'}}{\abs{\psi_D}^{1+\eps+\eps'}}
      }
      \abs{\idxup{\diff\psi_D} \ctrt  u}_{\vphilist}^2 \; .
    \end{equation*}
    On the left-hand-side, since $\ibddbar\sm\vphi_D$ is a smooth form
    on $X$ and $u$ is $L^2$ with respect to $\vphi$ and $\omega$, it
    follows that the limit as $\eps' \tendsto 0^+$ followed by $\eps
    \tendsto 0^+$ is finite and equal to $\int_X
    \idxup{\ibddbar\sm\vphi_D} \ptinner u u_{\vphilist}$.
    On the right-hand-side, notice that $\frac{2\eps'
      \:\theta'_{\eps'} \theta_{\eps'}}{\abs{\psi_D}^{1+\eps+\eps'}}
    \geq 0$.
    As $\theta_{\eps'} \ascendsto 1$ as $\eps' \descendsto 0$, it follows
    from the monotone convergence theorem that $\eps \int_X 
    \frac{\abs{\idxup{\diff\psi_D} \ctrt
        u}_{\vphilist}^2}{\abs{\psi_D}^{1+\eps}} < \infty$ for any
    $\eps >0$.
    This in turn implies that $\lim_{\eps' \tendsto 0^+} 2\eps' \int_X
    \frac{\theta'_{\eps'} \theta_{\eps'}}{\abs{\psi_D}^{1+\eps+\eps'}}
    \abs{\idxup{\diff\psi_D} \ctrt  u}_{\vphilist}^2 = 0$, as
    $\theta'_{\eps'} \theta_{\eps'}$ is bounded uniformly in $\eps'$.
    Therefore, the limit of the above equality as $\eps' \tendsto
    0^+$ followed by $\eps \tendsto 0^+$ yields the desired result.
  \end{proof}

  There is no need of further clarification for Step
  \ref{item:residue-of-final-inner-prod} in the case where $\vphi_F$
  and $\vphi_M$ are smooth (thus $X\setminus X^\circ =\emptyset$).
  The proof of Theorem \ref{thm:main-result} is therefore completed
  for this case.
  
}


\subsection{Proof for the case with singular $\vphi_F$ and $\vphi_M$}
\label{sec:pf-singular-vphi_FM}


In this section, $\vphi_F$ and $\vphi_M$ possess neat analytic
singularities described as in Section \ref{sec:setup}.
Recall that $X^\circ = X \setminus \paren{P_F \cup P_M}$, where $P_F =
\vphi_F^{-1}(-\infty)$ and $P_M = \vphi_M^{-1}(-\infty)$ are the polar
sets, and $\clomega$ is a complete \textde{Kähler} metric on $X^\circ$
given by the formula in item \eqref{item:setup-clomega} in Section
\ref{sec:setup} such that $\clomega \geq \omega$.
Write
\begin{equation*}
  \psi_{FM} :=\psi_{P_F \cup P_M}
  =\phi_{P_F\cup P_M} -\sm\vphi_{P_F \cup P_M}
  =:\phi_{FM} -\sm\vphi_{FM}
\end{equation*}
for convenience in what follows.
Note that $\log\paren{e\log\abs{\ell\psi_{FM}}}$ is a smooth
exhaustive function on $X^\circ$ with
$\abs{d\paren{\log\paren{e\log\abs{\ell\psi_{FM}}}}}_{\clomega}^2
\lesssim 1$ on $X^\circ$ by the choice of $\clomega$ (and the constant
$\ell \gg e$). Let $\chi \colon [0,\infty) \to [0,1]$ be a smooth
non-increasing cut-off function such that $\res{\chi}_{[0,1]} \equiv
1$ and $\res{\chi}_{[2,\infty)} \equiv 0$.
For every $\nu \in \Nnum$, set $\chi_{\nu} := \chi \circ
\paren{\frac{1}{2^\nu} \log\paren{e\log\abs{\ell\psi_{FM}}}}$ and
$\chi'_{\nu} :=\chi' \circ \paren{\frac{1}{2^\nu}
  \log\paren{e\log\abs{\ell\psi_{FM}}}}$ for convenience (where
$\chi'$ is the derivative of $\chi$).
Notice that both $\chi_\nu$ and $\chi'_\nu$ are compactly supported in
$X^\circ$ for $\nu < +\infty$ and $\chi_\nu \ascendsto 1$ pointwisely
on $X^\circ$ as $\nu \ascendsto +\infty$.
As $\abs{\chi'_\nu}$ is bounded uniformly in $\nu$, it follows that
\begin{equation} \label{eq:d-chi-bound}
  \abs{d\chi_\nu}_{\clomega}^2 \lesssim \frac{1}{2^{2\nu}} \; ,
\end{equation}
in which the constant involved in $\lesssim$ is independent of $\nu$
(so $\seq{\chi_\nu}_{\nu\in\Nnum}$ is just the exhaustive sequence of
cut-off functions introduced in \cite{Demailly}*{Ch.~VIII, Lemma (2.4
  c)}).

\subsubsection{Justification for Step \ref{item:pick-u-decompose-v}}

The justification for $\clt v = \clt v_{(2)} + v_{(\infty)}$ is
already given by Lemma \ref{lem:regularity-of-v} together with Remark
\ref{rem:regularity-v-for-sing-F-M}.

\subsubsection{Justification for Step \ref{item:results-from-harmonic-u-and-BK}}

To justify the arguments in Step
\ref{item:results-from-harmonic-u-and-BK} for the case where $\vphi_F$
and $\vphi_M$ have neat analytic singularities, it suffices to
supplement the proof of Proposition \ref{prop:residue-of-dpsi-ctrt-u}
with the clarification of the regularity of $\clt u \in \Harm$ along
$X \setminus X^\circ =P_F \cup P_M$ which guarantees that the proof
remains valid (thus \eqref{eq:limit-ctrt-u=0} holds true).
Moreover, a simple adjustment, making use of the completeness of
$\clomega$ on $X^\circ$, to the proofs of Proposition
\ref{prop:nabla-u=0_curv-u=0} and Corollary \ref{cor:dfadj_M-su=0} is
made to guarantee that $\clt u$ and $s \clt u$ still satisfy the
Bochner--Kodaira formulas $\BK_{\vphilist}$ and $\tBK_{\vphilist M}$
respectively.

The regularity of $\clt u$ along $X \setminus X^\circ$ is controlled
as follows.

\begin{prop} \label{prop:regularity-of-clt-u}
  For any $\clt u \in \Harm$, one has $\frac{\clt u}{\sect_D}$ being
  smooth on $X^\circ$.
  Furthermore, 
  on an admissible open set $V \subset X$ with respect to $(\vphi_F,
  \vphi_M, \psi_D)$ in the holomorphic coordinate system $(z_1,\dots,
  z_n)$ such that $D \cap V = \set{z_1 \dotsm z_{\sigma_V} = 0}$
  and
  \begin{equation*}
    \res{\phi_{FM}}_{V}
    =\smashoperator[r]{\sum_{j=\sigma_V+1}^{\sigma_V+\mu}}
    \:\log\abs{z_j}^2
    =\sum_{k=1}^{\mu} \log\abs{w_k}^2
    \qquad\paren{\text{setting }\, w_k := z_{\sigma_V+k}}
  \end{equation*}
  for some integer $\mu = \mu_V \in [0 , n-\sigma_V]$ ($\mu=0$ when $V
  \cap \paren{X\setminus X^\circ} =\emptyset$),
  one also has
  \begin{equation*}
    \res{\frac{\clt u}{\sect_D}}_{V}
    \in \smform/n,q/\paren{V} \left[
      \frac{1}{\abs{\psi_{FM}}} ,\: \frac{1}{\log\abs{\ell\psi_{FM}}} ,\:
      \frac{1}{w_1},\: \frac{1}{\conj{w_1}}, \dots,\: \frac{1}{w_{\mu}},\:
      \frac{1}{\conj{w_{\mu}}} 
    \right] \; ,
  \end{equation*}
  where the right-hand-side is generated as an
  $\smform/n,q/\paren{V}$-algebra.
\end{prop}

\begin{proof}
  Following the argument as in Section
  \ref{sec:justify-step-II-smooth}, the refined hard Lefschetz theorem
  in \cite{Matsumura_injectivity-lc}*{Thm.~3.3} (see also Theorem
  \ref{thm:refined-hard-Lefschetz}) implies that $*_{\clomega} \clt u$
  is holomorphic on $X$, which, together with the fact that $u$ being
  $L^2$ with respect to $e^{-\phi_D}$, in turn implies that
  $\frac{\clt u}{\sect_D}$ is smooth on $X^\circ$.

  For the remaining claim, note that, as $\frac{*_{\clomega} \clt
    u}{\sect_D}$ is smooth (indeed holomorphic) on $X$, it follows
  from the formula of the Hodge $*$-operator or 
  \begin{equation*} 
    \frac{\clt u}{\sect_D} \wedge \conj{\paren{\frac{*_{\clomega} \clt u}{\sect_D}}}
    =\abs{\frac{\clt u}{\sect_D}}_{\clomega}^2
    =\abs{\frac{*_{\clomega} \clt u}{\sect_D}}_{\clomega}^2
    \dvol_{X^\circ, \clomega}  
  \end{equation*}
  that the singularities of $\frac{\clt u}{\sect_D}$
  along $X \setminus X^\circ$ is determined by those of
  $d\vol_{X^\circ, \clomega} = \frac{\clomega^{\wedge n}}{n!}$.
  From the definition of $\clomega$ in item \eqref{item:setup-clomega}
  in Section \ref{sec:setup} and the expression $\psi_{FM} =\phi_{FM}
  -\sm\vphi_{FM} =\sum_{k=1}^\mu \log\abs{w_k}^2 -\sm\vphi_{FM}$,
  one can express $\frac{\clomega^{\wedge n}}{n!}$ on $V$ as 
  \begin{align*}
    \frac{\clomega^{\wedge n}}{n!}
    &=\frac{1}{n!} \paren{\alert{
      2\omega +
      \frac{ \ibddbar{\psi_{FM}}}{\abs{\psi_{FM}}
      \paren{\log\abs{\ell{\psi_{FM}}}}^{2}}
      }
      +\frac{1 +\frac{2}{\log\abs{\ell{\psi_{FM}}}}}{\abs{\psi_{FM}}^{2}
      \paren{\log\abs{\ell{\psi_{FM}}}}^2}
      \ibar \diff{\psi_{FM}} 
      \wedge \dbar{\psi_{FM}} }^{\wedge n} \\
    &=\frac{1}{n!} \paren{\alert{\alpha'^{\,1,1}}
      +\frac{
        1 +\frac{2}{\log\abs{\ell{\psi_{FM}}}} 
      } {\abs{\psi_{FM}}^{2}
      \paren{\log\abs{\ell{\psi_{FM}}}}^2}
      \:\ibar \paren{\sum_{j=1}^\mu \frac{dw_j}{w_j} -\diff\sm\vphi_{FM}}
      \!\!\wedge\!\! \paren{\sum_{j=1}^\mu \frac{d\conj w_j}{\conj w_j}
      -\dbar\sm\vphi_{FM}}
      }^{\wedge n} \\
    &=
      \begin{aligned}[t]
        &\alpha^{n,n} 
        +2\Re\paren{\sum_{j=1}^\mu
          \frac{dw_j \;\wedge \alpha^{n-1,n}_{j}}{w_j \abs{\psi_{FM}}^{2}
            \paren{\log\abs{\ell{\psi_{FM}}}}^2} 
          \;+ \;\sum_{\mathclap{1\leq j<k\leq \mu}} \;\;\:
          \frac{dw_j \wedge dw_k \;\wedge \alpha^{n-2,n}_{jk}}{w_j
            w_k  \abs{\psi_{FM}}^{4} \paren{\log\abs{\ell{\psi_{FM}}}}^4} 
        +\dots } \\
        &+\sum_{j,k=1}^\mu \frac{dw_j \wedge d\conj w_k \;\wedge
          \alpha^{n-1,\mathrlap{n-1}}_{j\conj k}} {w_j \conj w_k
          \abs{\psi_{FM}}^{2} \paren{\log\abs{\ell{\psi_{FM}}}}^2} \quad
        +\;\; \sum_{\mathclap{\substack{1\leq j_1 < j_2 \leq \mu \\ 1\leq k_1 < k_2
            \leq \mu}}} \;
        \frac{dw_{j_1} \wedge dw_{j_2} \wedge d\conj w_{k_1} \wedge
          d\conj w_{k_2} \;\wedge \alpha^{n-2,n-2}_{j_1 j_2 \conj
          k_1 \conj k_2}}{w_{j_1} w_{j_2} \conj w_{k_1} \conj w_{k_2}
          \abs{\psi_{FM}}^{4} \paren{\log\abs{\ell{\psi_{FM}}}}^4}  \\
        &+ \dots +\frac{\bigwedge_{j=1}^\mu dw_j \wedge d\conj w_j}
        {\prod_{j=1}^\mu \abs{w_j}^2 \cdot \abs{\psi_{FM}}^{2\mu}
          \paren{\log\abs{\ell{\psi_{FM}}}}^{2\mu}} \wedge \alpha^{n-\mu,n-\mu} \; ,
      \end{aligned} 
  \end{align*}
  where each $\alpha^{p',q'}_{j_1 j_2 \dots \conj k_1 \conj k_2 \dots}$
  is a $(p',q')$-form with \emph{continuous coefficients on $V$} such
  that
  \begin{equation*}
    \alpha^{p',q'}_{j_1 j_2 \dots \conj k_1 \conj k_2 \dots}
    \in \smform/p',q'/\paren{V} \left[
      \frac{1}{\abs{\psi_{FM}}} , \: \frac{1}{\log\abs{\ell\psi_{FM}}}
    \right] \; .
  \end{equation*}
  Note that $\alpha'^{\,1,1}$ can also be viewed as an $(1,1)$-form
  with continuous coefficients on $V$ since
  \begin{equation*}
    \frac{\ibddbar\phi_{FM}}{\logpole|\psi_{FM}|//[\ell]<2>} = 0
    \quad\text{ on } X \text{ as a current.}
  \end{equation*}
  (Recall the \textfr{Poincaré}--Lelong formula which states that
  $\ibddbar\phi_{FM} = [P_F \cup P_M]$, the current of integration
  along $P_F \cup P_M$.)
  Since $\frac{\clt u}{\sect_D}$ has at worst the singularities of
  $\frac{\clomega^{\wedge n}}{n!}$, the last claim follows by
  observing the formula of $\frac{\clomega^{\wedge n}}{n!}$.
\end{proof}

Recall that $D$ and $X \setminus X^\circ =P_F \cup P_M$ have no common
components and intersect each other with snc.
Recall also that $\vphi_F$ has only neat analytic singularities along
$P_F$.
One implication of Proposition \ref{prop:regularity-of-clt-u} is that,
in view of Fubini's theorem, the derivatives of $\frac{\clt
  u}{\sect_D}$ in the normal directions of any components of $D$ (more
precisely, derivatives with respect to $r_1, \dots, r_{\sigma_V}$ on
an admissible set $V$) are locally $L^2$ with respect to $\vphi_F$ and
$\clomega$ in $X$ (not only in $X^\circ$), because $\frac{\clt
  u}{\sect_D}$ is also $L^2$ with respect to $\vphi_F$ and $\clomega$
and the derivatives of $\frac{1}{\abs{\psi_{FM}}}$ and
$\frac{1}{\log\abs{\ell\psi_{FM}}}$ with respect to $r_1, \dots,
r_{\sigma_V}$ are more readily integrable than their primitives
(note also that $\frac{1}{w_k}$ and $\frac{1}{\conj{w_k}}$ are simply
constant functions with respect to these derivatives).
The derivatives of coefficients of $\clomega$ and its inverse (as well
as $\frac{1}{\det\clomega}$) with respect to $r_1, \dots,
r_{\sigma_V}$ are also more readily integrable than their primitives
for the same reason (as they live in the algebra generated over
bounded functions on $V$ with generators given by all
$\frac{1}{w_k}$'s, $\frac{1}{\conj{w_k}}$'s,
$\frac{1}{\abs{\psi_{FM}}}$, $\frac{1}{\log\abs{\ell\psi_{FM}}}$ and
the derivatives of the generators with respect to $r_1, \dots,
r_{\sigma_V}$).
Such argument is used in the following proposition.

\begin{prop}[cf.~Proposition \ref{prop:residue-of-dpsi-ctrt-u}]
  \label{prop:residue-of-dpsi-ctrt-clt-u}
  With $\frac{\clt u}{\sect_D}$ satisfying the conclusion in
  Proposition \ref{prop:regularity-of-clt-u} on any admissible open
  sets $V \subset X$, one has
  \begin{equation*}
    \lim_{\eps \tendsto 0^+} \eps
    \int_{\mathrlap{X^\circ}\;\;}
    \frac{
      \abs{\dbar\psi_D \otimes \clt u}_{\vphilist}^2
    }{\abs{\psi_D}^{1+\eps}}
    =\pi \sum_{i \in I_D} \int_{D_i} \abs{
      \PRes[D_i](\dbar\psi_D \otimes \frac{\clt u}{\sect_D})
    }_{\vphilist F}^2 < +\infty \; ,
  \end{equation*}
  where $D =\sum_{i \in I_D} D_i$ is the decomposition of $D$ into
  irreducible components in $X$ and $\PRes[D_i]$ is the
  \textfr{Poincaré} residue map corresponding to the restriction
  from $X$ to $D_i$.
  The equation \eqref{eq:limit-ctrt-u=0} therefore also holds true.
\end{prop}

\begin{proof}
  The proof is exactly the same as the proof of Proposition
  \ref{prop:residue-of-dpsi-ctrt-u}, provided that the singularities
  of $\frac{\clt u}{\sect_D}$ along $X\setminus X^\circ$ is taken care
  of.

  Using the same notation as in the proof of Proposition
  \ref{prop:residue-of-dpsi-ctrt-u}, consider again a finite cover
  $\set{V_\gamma}_{\gamma \in \Gamma}$ of $X$ by admissible open sets
  $V_\gamma$ and write
  \begin{equation*}
    \dbar\psi_D \otimes \clt u
    =\sum_{j=1}^{\sigma_{V_\gamma}} \paren{
      d\conj{z_j^\gamma}
      -\frac{\conj{z_j^\gamma} \:\dbar\sm\vphi_D}{\sigma_{V_\gamma}}
    } \otimes \frac{\clt u}{\conj{z_j^\gamma}}
    =\sum_{j=1}^{\sigma_{V_\gamma}}
    \underbrace{
      \paren{
        d\conj{z_j^\gamma}
        -\frac{\conj{z_j^\gamma} \:\dbar\sm\vphi_D}{\sigma_{V_\gamma}}
      } \otimes \frac{\clt u}{\sect_D} \:\frac{z_j^\gamma}{\conj{z_j^\gamma}}
    }_{\textstyle =: \: dz_j^\gamma \wedge g_j^\gamma}
    \:\prod_{\substack{k=1 \\ k\neq j}}^{\sigma_{V_\gamma}}
    z_k^\gamma 
  \end{equation*}
  on each $V_\gamma$ for $\gamma \in \Gamma$ (where $g_j^\gamma$ is
  chosen as in footnote \ref{fn:definition-of-g_j} on page
  \pageref{fn:definition-of-g_j} such that $g_j^\gamma
  =\fdiff{z_j^\gamma} \ctrt \paren{dz_j^\gamma \wedge g_j^\gamma}$).
  With the conclusion of Proposition \ref{prop:regularity-of-clt-u}
  and the fact that $\norm{\frac{\clt u}{\sect_D}}_{\vphilist F}^2
  =\norm{\clt u}_{\vphilist}^2 < +\infty$, one can check readily (in
  view of Fubini's theorem) that both conditions in
  \eqref{eq:regularity-on-g_p} are satisfied (with
  $\inner{g_j^\gamma}{g_{j'}^\gamma}_{\vphilist F}$ in place of
  $\inner{g_p}{g_{p'}}_{\vphilist F<\omega>}$ for any $j,j' =1,\dots,
  \sigma_{V_\gamma}$).
  Indeed, this can be seen from the following facts:
  \begin{enumerate}[label=(\alph*), ref=\alph*]
  \item it follows from the definition of $g_j^\gamma$ that
    $\abs{dz_j^\gamma \wedge g_j^\gamma}_{\vphilist F}^2 \lesssim \abs{\frac{\clt
        u}{\sect_D}}_{\vphilist F}^2$ (thus the coefficient of the $(n-1,n-1)$-form
    $\abs{g_j^\gamma}_{\vphilist F}^2$ is in $\Lloc[1](V_\gamma)$,
    and, from the openness property of multiplier ideal sheaves and
    the fact that $\frac{*_{\clomega} \clt u}{\sect_D}$ is
    holomorphic, $\abs{g_j^\gamma}_{(1+\eps)\vphilist F}^2$ is also in
    $\Lloc[1](V_\gamma)$ for sufficiently small $\eps > 0$);

  \item $e^{-\vphi_F} \sim \prod_{k=1}^\mu \frac{1}{\abs{w_k}^{2
        b_k}}$ in the notation in Proposition
    \ref{prop:regularity-of-clt-u}, where $b_k \geq 0$ (but may not be
    integers);
    
  \item Proposition
    \ref{prop:regularity-of-clt-u}
    implies that $g_j^\gamma$, as well as
    $\abs{g_j^\gamma}_{\clomega}^2$ (without $e^{-\vphi_F}$), is a
    \emph{polynomial} in $\frac 1{\abs{w_1}}, \dots, \frac
    1{\abs{w_\mu}}$, $\frac 1{\abs{\psi_{FM}}}$, $\frac
    1{\slog[\ell]|\psi_{FM}|}$ and $\frac{1}{\det\clomega}$ with
    coefficients in $\smform/\bullet,\bullet/[X\,*]({V_\gamma})$,
    which are therefore \emph{smooth in the variables $r_1, \dots,
      r_n$ on ${V_\gamma}$}, where $r_j :=\abs{z_j}$, when the other
    variables are fixed (note that $\frac{1}{\det\clomega}$ is
    involved as $\abs{g_j^\gamma}_{\clomega}^2$ involves the
    coefficients of the inverse of $\clomega$);

  \item \label{item:det-omega-factorised}
    the computation of $\frac{\clomega^{\wedge n}}{n!}$ in the
    proof of Proposition \ref{prop:regularity-of-clt-u} implies that
    $\frac{1}{\det\clomega} =\abs{w_1 \dotsm w_\mu}^2
    \logpole|\psi_{FM}|/2\mu/[\ell]<2\mu> \: B$ for some positive
    (nowhere zero) continuous function $B$ on $V_\gamma$ which is
    smooth in the variables $r_1, \dots, r_{\sigma_{V_\gamma}}$.
  \end{enumerate}
  Taking into account the factorisation of $\frac 1{\det\clomega}$ in
  item \eqref{item:det-omega-factorised} in the polynomial expression of
  $\abs{g_j^\gamma}_{\clomega}^2$, it follows that
  $\abs{g_j^\gamma}_{\clomega}^2$ is a polynomial in variables
  $\abs{w_1}^{\pm 1}, \dots, \abs{w_\mu}^{\pm 1}$,
  $\abs{\psi_{FM}}^{\pm 1}$ and $\paren{\log\abs{\ell\psi_{FM}}}^{\pm
    1}$ over the algebra $\smform/n-1,n-1/[X\,*]\paren{V_\gamma}[B]$.
  By factoring out powers of the variables $\abs{w_1}, \dots,
  \abs{w_\mu}, \frac{1}{\abs{\psi_{FM}}}$ and
  $\frac{1}{\log\abs{\psi_{FM}}}$ (or possibly their reciprocals)
  without modifying the coefficients, one can write
  \begin{equation*} \label{eq:g_j-factorised} \tag{$*$}
    \abs{g_j^\gamma}_{\clomega}^2 = G \: \frac{
      \abs{w_1}^{m_1} \dotsm \abs{w_\mu}^{m_\mu}
    }{
      \logpole|\psi_{FM}|/p/[\ell]<p'>
    } \; ,
  \end{equation*}
  where $m_1, \dots, m_\mu, p, p' \in \Znum$ are the maximal (possibly
  negative) integers such that the function $G$ is a polynomial in
  $\abs{w_1}, \dots, \abs{w_\mu}$, $\abs{\psi_{FM}}$ and
  $\log\abs{\psi_{FM}}$ with coefficients in
  $\smform/n-1,n-1/[X\,*]\paren{V_\gamma}[B]$ and that $G$ is
  not divisible by $\abs{\psi_{FM}}$ and $\log\abs{\psi_{FM}}$
  over $\smform/n-1,n-1/[X\,*]\paren{V_\gamma}[B]$.
  With the fact that $\abs{g_j^\gamma}_{(1+\eps)\vphilist F}^2$ is integrable
  on $V_\gamma$, one can adjust the exponents $m_1, \dots, m_\mu$ and
  modify the coefficients of $G$ until the factor $Q :=\frac{
    \abs{w_1}^{m_1} \dotsm \abs{w_\mu}^{m_\mu}
  }{
    \logpole|\psi_{FM}|/p/[\ell]<p'>
  }$ is integrable with respect to $e^{-(1+\eps)\vphi_F}$ on
  $V_\gamma$ for all sufficiently small $\eps > 0$.
  Indeed, if $Q e^{-(1+\eps)\vphi_F}$ is not integrable for all $\eps > 0$,
  then $G$ has zeros along some divisor in $\psi_{FM}^{-1}(-\infty)$,
  say, $\set{w_k = 0}$, and can be factored into $\abs{w_k}^{m_k'}
  G'$, where $G' \in
  \smform/n-1,n-1/[X\,*]\paren{V_\gamma}[B, \abs{\psi_{FM}},
  \log\abs{\psi_{FM}}]$ ($\abs{w_k}$'s are absorbed into
  $\smform/n-1,n-1/[X\,*]\paren{V_\gamma}$) and $m_k' \in
  \Nnum$, and $Q \abs{w_k}^{m_k'} e^{-(1+\eps)\vphi_F}$ is then
  integrable for sufficiently small $\eps >0$ around general points of
  $\set{w_k=0}$ in $V_\gamma$.
  Therefore, by induction, one can assume that $Q$ satisfies the
  acclaimed integrability condition.

  There is a minimal exponent $\rs p \in \Nnum_{\geq 0}$ such that
  $\frac{G}{\abs{\psi_{FM}}^{\rs p}}$ is bounded on $V_\gamma$ (note
  that
  $\frac{\paren{\log\abs{\psi_{FM}}}^{\beta'}}{\abs{\psi_{FM}}^\beta}$
  is bounded on $V_\gamma$ for any $\beta > 0$ and $\beta' \geq 0$).
  As $Q\:\abs{\psi_{FM}}^{\rs p} e^{-(1+\eps)\vphi_F}$ is still integrable
  for all sufficiently small $\eps > 0$, one can then assume that
  $G$ in \eqref{eq:g_j-factorised} is bounded on $V_\gamma$
  (but now $G \in \smform/n-1,n-1/[X\,*]\paren{V_\gamma}[B,
  \abs{\psi_{FM}}^{\alert{\pm 1}}, \log\abs{\psi_{FM}}]$).
  Note that $G$ is smooth in the variables $r_1, \dots,
  r_{\sigma_{V_\gamma}}$ (when the other variables are fixed).
  (The boundedness of $G$ is used in the proof of Proposition
  \ref{prop:justification-final-inner-prod}.)

  By a similar analysis, $g_j^\gamma$ takes the form analogous to
  \eqref{eq:g_j-factorised} (but the corresponding factor $G$ is a
  polynomial over $\paren{\smform/0,1/[X\,*] \otimes
    \smform/n-1,q/[X\,*]}\paren{V_\gamma}[B]$) and satisfies the
  condition (which is stronger than \eqref{eq:regularity-on-g_p}) 
  \begin{equation*}
    \paren{\fdiff{r_{\sigma_{V_\gamma}}}}^{\alpha_{\sigma_{V_\gamma}}} \dotsm
    \paren{\fdiff{r_{1}}}^{\alpha_1} g_j^\gamma
    \in \Lloc[2]\paren{{V_\gamma}} \quad
    \text{ for any } \alpha_1, \dots,\alpha_{\sigma_{V_\gamma}} \in \Nnum_{\geq 0}
  \end{equation*}
  when the discussion before the statement of this
  proposition (on the integrability of the derivatives of $\frac
  1{\abs{\psi_{FM}}}$, $\frac 1{\slog[\ell]|\psi_{FM}|}$ and $\frac
  1{\det\clomega}$, which applies also to multiples of
  $\abs{\psi_{FM}}$ and $\log\abs{\psi_{FM}}$) is taken into account.
  The fact that $\inner{g_j^\gamma}{g_{j'}^\gamma}_{\vphilist F}$
  satisfies the conditions in \eqref{eq:regularity-on-g_p} (even when
  $j \neq j'$) then follows from the Cauchy--Schwarz inequality and
  the Leibniz rule for differentiation.
  Therefore, Proposition \ref{thm:residue-fcts-and-norms} with Remark
  \ref{rem:residue-with-clomega} can be applied.
  The rest of the argument is the same as in the proof of Proposition
  \ref{prop:residue-of-dpsi-ctrt-u}.  
\end{proof}

The proofs of Proposition \ref{prop:nabla-u=0_curv-u=0} and Corollary
\ref{cor:dfadj_M-su=0}, which confirm the claims in
\eqref{eq:nable-u=0_ibddbar-vphi_F-u=0} and \eqref{eq:dfadj-su=0},
need adjustments only in replacing $\omega$ by $\clomega$ and
$\theta_\eps$ by $\chi_\nu \theta_\eps$ (where $\chi_\nu$ is defined
at the beginning of Section \ref{sec:pf-singular-vphi_FM}).
They remain valid without further changes once the following lemma is
observed.

\begin{lemma} \label{lem:limits-of-d-of-cutoff-fcts}
  Given $\frac{\clt u}{\sect_D}$ satisfying the conclusion in
  Proposition \ref{prop:regularity-of-clt-u}, one has
  \begin{equation*}
    \lim_{\eps \tendsto 0^+} \lim_{\nu \tendsto \infty}
    \int_{X^\circ} \abs{\dbar\paren{\chi_\nu \theta_\eps} \otimes
      \clt u}_{\vphilist}^2 = 0 \; .
  \end{equation*}
  Note that the statement with $\dbar\paren{\chi_\nu \theta_\eps}
  \otimes \clt u$ replaced by $\idxup{\dbar\paren{\chi_\nu
      \theta_\eps}} \ctrt \clt u$ also holds true since
  $\abs{\idxup{\dbar\paren{\chi_\nu \theta_\eps}} \ctrt \clt
    u}_{\vphilist}^2 \leq \abs{\dbar\paren{\chi_\nu \theta_\eps}
    \otimes \clt u}_{\vphilist}^2$.
\end{lemma}

\begin{proof}
  Noting that $\chi_\nu$, $\theta_\eps$ and $\theta'_\eps$ are bounded
  uniformly in $\nu$ and $\eps$, and recalling the bound on
  $d\chi_\nu$ (thus on $\dbar\chi_\nu$) in \eqref{eq:d-chi-bound}, a
  direct computation yields
  \begin{align*}
    \abs{\dbar\paren{\chi_\nu \theta_\eps} \otimes
      \clt u}_{\vphilist}^2
    &=\abs{\paren{\theta_\eps \dbar\chi_\nu   +\chi_\nu \dbar\theta_\eps} \otimes
      \clt u}_{\vphilist}^2 \\
    &\leq 2 \abs{\theta_\eps \dbar\chi_\nu \otimes \clt
      u}_{\vphilist}^2
    +2\eps^2 \abs{\chi_\nu \theta'_\eps}^2
    \frac{ \abs{\dbar\psi_D \otimes \clt u
      }_{\vphilist}^2}{\abs{\psi_D}^{2+2\eps}} \\
    &\lesssim \frac{1}{2^{2\nu}} \abs{\clt u}_{\vphilist}^2
      +\eps^2 \frac{ \abs{\dbar\psi_D \otimes \clt u
      }_{\vphilist}^2}{\abs{\psi_D}^{2+2\eps}} \; ,
  \end{align*}
  where the constant involved in $\lesssim$ is independent of $\nu$
  and $\eps$.
  Noting also the fact $\abs{\psi_D} \geq 1$, the claim then follows
  immediately from Proposition \ref{prop:residue-of-dpsi-ctrt-clt-u}.
\end{proof}

This completes the justification for Step \ref{item:results-from-harmonic-u-and-BK}.

\subsubsection{Justification for Step \ref{item:Takegoshi-argument}}
\label{sec:justification-Takegoshi-argument-singular}

There is no extra clarification needed for Step
\ref{item:by-parts-on-norm-su}. 
To justify the argument in Step \ref{item:Takegoshi-argument}, it
suffices to reprove Proposition \ref{prop:tBK-valid-for-u} under the
assumption that $\clt u$ has singularities along $X \setminus X^\circ$
described as in Proposition \ref{prop:regularity-of-clt-u}.
This can be achieved by replacing $\theta_{\eps'}$ in the proof of
Proposition \ref{prop:tBK-valid-for-u} by $\chi_{\nu} \theta_{\eps'}$
and noticing the identity
\begin{align*}
  &~2\Re \int_X
  \inner{
    \idxup{\diff\paren{\chi_\nu \theta_{\eps'}}} \ctrt \clt u
  }{
    \frac{
      \idxup{\diff\psi_D} \ctrt \chi_\nu \theta_{\eps'} \clt u
    }{\abs{\psi_D}}
  }_{\vphilist} \abs{\psi_D}^{1-\eps} \\
  =
  &~2\Re\int_X 
    \inner{
      \idxup{\diff\chi_\nu} \ctrt \theta_{\eps'}\clt u
    }{
      \frac{
        \idxup{\diff\psi_D} \ctrt \chi_\nu \theta_{\eps'} \clt u
      }{\abs{\psi_D}^\eps}
    }_{\vphilist}
    +
    \int_X \frac{
      2\eps' \chi_\nu^2 \theta'_{\eps'} \theta_{\eps'}
    }{\abs{\psi_D}^{1+\eps+\eps'}}
    \abs{\idxup{\diff\psi_D} \ctrt \clt u}_{\vphilist}^2 \; .
\end{align*}
The first integral on the right-hand-side converges to $0$ as $\nu
\tendsto +\infty$ thanks to the bound on $\diff\chi_\nu$ in
\eqref{eq:d-chi-bound}.
The desired result is obtained by following the remaining arguments in
the proof of Proposition \ref{prop:tBK-valid-for-u} and taking the
limits $\nu \ascendsto +\infty$, $\eps' \descendsto 0^+$ and $\eps
\tendsto 0^+$ in order.

\subsubsection{Justification for Step \ref{item:residue-of-final-inner-prod}}

The justification is done with the following proposition.
\begin{prop} \label{prop:justification-final-inner-prod}
  Under the assumption $\vphi_F$ and $\vphi_M$ having neat analytic
  singularities with snc along $X \setminus X^\circ$, together with
  the regularity of $\frac{\clt u}{\sect_D}$ along $X \setminus
  X^\circ$ described in Proposition \ref{prop:regularity-of-clt-u},
  one has
  \begin{equation*}
    \int_X \abs{
      \inner{
        \idxup{\diff\psi_D}* \ctrt
        \frac{\clt u}{\sect_D}
      }{
        \:v_{(\infty)} e^{-\frac 12 \sm\vphi_D -\frac 12 \vphi_M}}
    }_{\vphilist F} < +\infty \; ,
  \end{equation*}
  which also implies that 
  \begin{equation*}
    \int_X \frac{
      \abs{
        \inner{
          \idxup{\diff\psi_D}* \ctrt
          \frac{\clt u}{\sect_D}
        }{
          \:v_{(\infty)} e^{-\frac 12 \sm\vphi_D -\frac 12 \vphi_M}}
      }_{\vphilist F}
    }{\abs{\sect_D}_{\sm\vphi_D} \:\abs{\psi_D}^{1+\eps}}
    < +\infty \quad \text{ for any } \eps \in \fieldR \; .
  \end{equation*}
\end{prop}

\begin{proof}
  It suffices to prove both claims on an arbitrary admissible open set
  $V \Subset X$ with respect to $(\vphi_F, \vphi_M, \psi_D)$.
  Let $(z_1, \dots, z_n)$ be a holomorphic coordinate system on $V$
  such that $D \cap V =\set{z_1 \dotsm z_{\sigma_V} =0}$.
  Decompose $V$ into the product $U \times W$ of polydiscs $U =
  U^{\sigma_V}$ and $W =W^{n-\sigma_V}$ such that
  $(z_1,\dots,z_{\sigma_V})$ and $(z_{\sigma_V+1},\dots,z_n) =(w_1,
  \dots, w_\mu, z_{\sigma_V +\mu+1}, \dots, z_n)$ are
  coordinate systems on $U$ and $W$ respectively (in the notation in
  Proposition \ref{prop:regularity-of-clt-u}).

  First note that $\int_X
  \abs{v_{(\infty)}}_{\vphilist|\sm\vphi_D+\vphi_F|M}^2 \leq \int_X
  \abs{v_{(\infty)}}_{\vphilist|\sm\vphi_D+\vphi_F|M<\alert{\omega}>}^2 < +\infty$
  according to Lemma \ref{lem:regularity-of-v} (together with Remark
  \ref{rem:regularity-v-for-sing-F-M}).
  Furthermore, by writing
  \begin{equation*}
    \idxup{\diff\psi_D} \ctrt \frac{\clt u}{\sect_D}
    =\sum_{j =1}^{\sigma_V} \underbrace{
      \idxup{dz_j -\frac{z_j \:\diff\sm\vphi_D}{\sigma_V}} \ctrt
      \frac{\clt u}{\sect_D}
    }_{\textstyle =: \: dz_j \wedge g_j} \: \frac{1}{z_j}
  \end{equation*}
  and following the analysis (as well as the notation) in the proof of
  Proposition \ref{prop:residue-of-dpsi-ctrt-clt-u}, one has
  \begin{equation*}
    \abs{dz_j \wedge g_j}_{\vphilist F}^2
    =G_j \:\underbrace{
      \frac{\abs{w_1}^{m_1} \dotsm \abs{w_\mu}^{m_\mu}}{\logpole|\psi_{FM}|/p/[\ell]<p'>}
    }_{\textstyle =: \: Q_j} \:e^{-\vphi_F} \; ,
  \end{equation*}
  where $G_j \in \smform/n,n/[X\,*]\paren{V}[B,
  \abs{\psi_{FM}}^{\pm 1}, \log\abs{\psi_{FM}}]$ is \emph{bounded}
  on $V$ and $Q_j \:e^{-\vphi_F}$ is integrable on $V$.
  The fact that $\idxup{\frac{dz_j}{z_j}} \ctrt \frac{\clt
    u}{\sect_D} \in \smform/n,q-1/[X\, *]\paren{V \cap X^\circ}$ for
  all $j=1,\dots,\sigma_V$ according to
  \eqref{eq:smooth-dz-over-z-ctrt-u-over-sect_D} (thus having 
  no poles along $D \cap V\cap X^\circ$)
  implies that $G_j=\abs{z_j}^2 G_j'$ for some $G_j' \in
  \smform/n,n/[X\,*]\paren{V}[B, \abs{\psi_{FM}}^{\pm 1},
  \log\abs{\psi_{FM}}]$ such that $G_j'$ is also \emph{bounded} on
  $V$.
  This implies immediately the convergence of
  $\int_V \abs{\idxup{\diff\psi_D} \ctrt \frac{\clt
      u}{\sect_D}}_{\vphilist F}^2$.
  The first claim then follows.

  In view of Fubini's theorem, the first claim
  implies that the function
  \begin{equation*}
    U \ni x \mapsto
    \int_{\set{x} \times W}
    \abs{
      \inner{
        \idxup{\diff\psi_D}* \ctrt
        \frac{\clt u}{\sect_D}
      }{
        \:v_{(\infty)} e^{-\frac 12 \sm\vphi_D -\frac 12 \vphi_M}}
    }_{\vphilist F}
  \end{equation*}
  is bounded (indeed even continuous) on $U$, as the integrand above
  is continuous in $z_1, \dots, z_{\sigma_V}$ (when the other
  variables are fixed at almost every point in $W$) and it has an
  estimate 
  \begin{align*}
    \abs{
      \inner{
        \dotsm
      }{
        \dotsm
      }
    }_{\vphilist F}
    &\leq
    \paren{\abs{
        \idxup{\diff\psi_D}* \ctrt
        \frac{\clt u}{\sect_D}
      }_{\vphilist F}^2}^{\frac 12} \:
    \paren{\abs{
        v_{(\infty)} 
      }_{\vphilist|\sm\vphi_D+\vphi_F|M<\alert{\omega}>}^2}^{\frac 12}
    \\
    &\lesssim
    \paren{\sum_{j=1}^{\sigma_V} Q_j \:e^{-\vphi_F}\dvol_X}^{\frac 12} \:
    \paren{\abs{
        v_{(\infty)}
      }_{\vphilist|\sm\vphi_D+\vphi_F|M<\omega>}^2}^{\frac 12}
  \end{align*}
  which can be seen easily that the right-hand-side is dominated by some
  $(n,n)$-form $\Phi \dvol_X$ on $V$ such that $\Phi =\Phi(w_1, \dots, w_\mu)$
  depends only on the variables $w_1, \dots, w_\mu$ and is integrable
  on $W$ (see, for example, \cite{Apostol_analysis}*{Thm.~10.38}).
  Since $\frac{1}{\abs{\sect_D}_{\sm\vphi_D} \:
    \abs{\psi_D}^{1+\eps}}$ is integrable on $X$ for any $\eps \in
  \fieldR$ and has poles only along $D$, the second claim then follows
  again in view of Fubini's theorem.
\end{proof}

This completes the justification of Step
\ref{item:residue-of-final-inner-prod} and thus the proof of Theorem
\ref{thm:main-result}.




\end{document}